\documentclass[final,3p,times]{elsarticle}
\usepackage{amsfonts,amsmath,amsthm,amssymb}
\usepackage{verbatim}
\usepackage{color}
\usepackage{moreverb}
\usepackage{mathrsfs}
\usepackage{caption}
\usepackage{multirow}

\usepackage{algorithm}
\usepackage{algpseudocode}
\usepackage{fixltx2e}

\usepackage{subfigure}
\usepackage{caption}
\usepackage{wrapfig}
\usepackage{epsfig}
\usepackage{subfig}

\usepackage{tikz}
\usetikzlibrary{arrows}

\def\NN{\hbox{I\kern-.2em\hbox{N}}}
\def\RR{{\mathop{{\rm I}\kern-.2em{\rm R}}\nolimits}}

\def\f{{\boldsymbol F}}

\def\n{{\boldsymbol n}}
\def\x{{\boldsymbol x}}
\def\y{{\boldsymbol y}}

\def\dx{\mathop{dx}}
\def\dy{\mathop{dy}}

\def\dz{\mathop{d \boldsymbol z}}
\def\dt{\mathop{d \boldsymbol t}}
\def\dtg{\mathop{d \Gamma_{\y}}}
\def\mg{{\mathcal{G}}}
\def\mh{{\mathcal{H}}}

\def\mk{{\mathcal{K}}}

\def\d{{\boldsymbol D}}
\def\s{{\boldsymbol s}}
\def\t{{\boldsymbol t}}
\def\z{{\boldsymbol z}}
\def\a{{\boldsymbol a}}

\definecolor{orange}{rgb}{1,0.5,0}
\definecolor{gray}{rgb}{0.5,0.5,0.5}

\def\cc{}

\newcommand{\vertiii}[1]{{\left\vert\kern-0.25ex\left\vert\kern-0.25ex\left\vert #1
    \right\vert\kern-0.25ex\right\vert\kern-0.25ex\right\vert}}

\newcommand{\be}{\begin{equation}}
\newcommand{\ee}{\end{equation}}
\newcommand{\ba}{\begin{eqnarray}}
\newcommand{\ea}{\end{eqnarray}}

\newtheorem{thm}{Theorem}
\newtheorem{lma}[thm]{Lemma}
\newtheorem{rmk}[thm]{Remark}
\newtheorem{prn}[thm]{Proposition}
\newtheorem{cor}[thm]{Corollary}

\begin{document}

\begin{frontmatter}

\title{Isoparametric singularity extraction technique\\ for 3D potential
problems in BEM 
}


%
\author[label1,label2]{Tadej Kandu\v{c}}
\cortext[cor1]{Corresponding Author}
\ead{tadej.kanduc@fmf.uni-lj.si}
%
\address[label1]{Faculty of Mathematics and Physics and Faculty of Mechanical Engineering, University of Ljubljana, Ljubljana, Slovenia}

\begin{abstract}
To solve boundary integral equations for potential problems using collocation Boundary Element Method (BEM) on smooth curved 3D geometries,  an analytical singularity extraction technique is employed. By adopting the isoparametric approach, curved geometries that are represented by mapped rectangles or triangles from the parametric domain are considered. 
The singularity extraction on the governing singular integrals can be performed either as an operation of subtraction or division, each having some advantages.

A particular series expansion of a singular kernel about a source point is investigated. The series in the intrinsic coordinates consists of functions of a type $R^p x^q y^r$, where $R$ is a square root of a quadratic bivariate homogeneous polynomial, corresponding to the first fundamental form of a smooth surface, and $p,q,r$ are integers, satisfying $p\leq -1$ and $q,r \geq 0$. By extracting more terms from the series expansion of the singular kernel, the smoothness of the regularized kernel at the source point can be increased. Analytical formulae for integrals of such terms are obtained from antiderivatives of $R^p x^q y^r$, using recurrence formulae, and by evaluating them at the edges of rectangular or triangular parametric domains.

Numerical tests demonstrate that the singularity extraction technique can be a useful prerequisite for a numerical quadrature scheme to obtain accurate evaluations of the governing singular integrals in 3D collocation BEM.

\end{abstract}
\begin{keyword}
Analytical integration \sep Singular integral \sep Singularity extraction \sep Boundary element method\sep Isogeometric analysis
\end{keyword}

\end{frontmatter}


\section{Introduction}


Isogeometric Analysis ({IGA}) \cite{Hughes_2005,LibroHughes} has in recent years brought a renewed interest in developing highly accurate simulation models in Boundary Element Method (BEM); e.g., see advances in local and adaptive mesh refinements \cite{GINNIS2014425,Dirk1,Dirk2},  
trimming \cite{BEER2015776},
fast matrix formation \cite{Dolz18}, 
fluid–structure interaction \cite{HELTAI2017522}, 
fracture \cite{PENG2017151},
shape optimization \cite{KOSTAS2018600}, 
acoustics and linear elasticity problems \cite{SimpsonScott,VENAS2020112670,TAUS2019112591},
and references therein. A predominant setting in IGA combines exact representation of the geometry using CAD standards such as polynomial and rational splines, and the same type of basis functions to represent the approximate solutions of PDEs. 
To fully profit by  high convergence orders of approximate solutions of isogeometric methods, the construction of efficient numerical rules for (singular) integrals is one of the most crucial challenges to be addressed in IGA-BEM \cite{ACDS3,AIMI2020113441,FGKSS_2019,heltai2014nonsingular,TauRodHug}. Although the Gaussian quadratures are usually considered as the optimal choice for smooth integrands with polynomial bases, there are better alternatives for integrals with B-splines bases, which exploit the interelement continuity of splines; see the review paper \cite{IGA_quads}. 

A special care is needed to treat singular integrals with numerical schemes -- a continuous research in the last decades (mostly outside isogeometric context) to design accurate, efficient and easy to implement integration rules produced a vast number of approaches, including explicit formulae for flat domains \cite{Fata2009ExplicitEF,jarvenpaa}, radial integration \cite{GAO2002905,gao_cmame}, (near) singularity smoothing \cite{JOHNSTON2013148,Monegato_airfol,SCUDERI2009406,telles1987self,ZHANG201557},  true  desingularisation \cite{KLASEBOER2009796} and other techniques \cite{Chen_Hong_Review_dualBEM}; see also references therein. Coordinate transformation, such as polar \cite{Guiggiani_generalAlg_HBIE,Khayat_polarSingInteg,Rong_polarTransform,TauRodHug} and Duffy transformation \cite{Duffy,Duffy4,Duffy2,Duffy3}, can reduce and simplify the nature of the singularity. 
An established technique, used in various forms, is the singularity extraction, in which the singular part of the integral is evaluated analytically \cite{AD2002,Fata2009ExplicitEF,gao_cmame,Guiggiani_generalAlg_HBIE,jarvenpaa}. 
%
%
In the preceding research to this paper, a promising approach for curved isoparametric boundaries for 2D problems combines an elegant singularity extraction and a local spline quasi-interpolation operators \cite{AIMI2020113441,CFSS18,FGKSS_2019}. For this type of quadrature schemes the optimal convergence orders of the approximate solution can be recovered with a small numbers of quadrature nodes \cite{FKdreams}. 

In this paper we study a singularity extraction for singular integrals for 3D potential problems. Focusing on isogeometric collocation BEM, curved boundary domains are represented as a set of rectangles or triangles in the parametric domain, mapped to the physical space using B-spline (or NURBS) functions. As a proof of concept, the main focus is on weakly and nearly singular integrals, appearing for Laplace problems, although the ideas could be applied to other types of integrals as well. 
The singularity extraction utilizes a particular series expansion of the singular kernel about a source point in the parametric domain -- to the best of author's knowledge this type of series expansions have not been studied before in this form in BEM.  No distortion of the integration domain occurs since the extraction is applied directly in the initial parametric domain; this is an important accuracy and efficiency feature since the same quadrature nodes can be used for regular and singular integrals and for several basis functions and source points. The expansion depends on the (higher) derivatives of the geometry parameterization, which are considered to be easily computable in IGA setting. For each summand in the series, multiplied with a polynomial basis function, recursive formulae to evaluate its double integral analytically are provided. The smoothness of the regularized kernel at the source point is controlled by the number of considered terms in the series. An additional constant cost to evaluate more terms can be outweighed by the higher accuracy of a quadrature for regularized integrals -- this is especially a preferred trade-off for isogeometric methods with high convergence orders.
The developed singularity extraction has already been applied in BEM for Laplace problems \cite{nash20,FGKSS_2021} and Helmholtz equation can be dealt similarly \cite{FKSS_2021}, since the appearing singularities are essentially the same.

The structure of the paper is as follows.  In Section {\cc 2.1} we recall a basic setting of integral formulation for potential problems -- it serves as a short introduction and motivation to Section {\cc 2.2}, where series expansions of kernels for the single and double layer potentials are derived. For completeness, we briefly outline the considered numerical integration schemes for the governing integrals in BIE in Section {\cc 2.3} and we provide some directions for implementation speedups for the analytical evaluation of the studied singular integrals in {\cc 2.4}. In {\cc 2.5}, smoothness of a more general regularized kernel is studied with respect to the two regularization techniques: subtraction and division. In Section {\cc 3} we perform tests of numerical integrations for the two common singular kernels for Laplace problems and for the two regularization techniques. In {\cc appendix} we {\cc collect} basic properties of functions $R^p x^q y^r$ and derive the recursive formulae to compute their antiderivatives.


\section{Singularity extraction}\label{singExtr}

In this section we demonstrate how to apply a singularity extraction technique on weakly singular kernels for the single and double layer potential for the 3D Laplace problems (see \cite{FKSS_2021} how to rewrite the kernels in Helmholtz equation). As a first step, a series expansion of the singular kernels about the source point in intrinsic coordinates is derived. A truncated series is an approximation of the singular kernel and it is used either in a singularity subtraction or in a singularity division to obtain regularized kernels. The remaining regular parts of the integrands can be locally approximated with a suitable polynomial. Therefore, the singular integrals in the underlying Boundary Integral Equations (BIE) are approximated as a finite sum of integrals, described in Section~\ref{basicInt}. At the end of the section we also analyze the smoothness of the derived regularized kernels.

\subsection{Integral formulation of the potential problem}

To be concise, let us consider 3D potential problems described by the Laplace equation with Dirichlet boundary conditions on finite volumes $\Omega$ with closed smooth boundary surfaces $\Gamma$,
\begin{align} \label{BVP_esse}
\left\{ \begin{array}{ll}
\Delta u=0& {\rm in}\; \Omega,\\
u=u_{D}& {\rm on}\; \Gamma,
\end{array} \right.
\end{align}
where the solution $u$ belongs to the Sobolev space $H^1(\Omega)$ and the Dirichlet boundary datum $u_D$ is in $H^{1/2}(\Gamma)$, the trace space of $H^1(\Omega)$. For more details we refer to \cite{Atkinson_2009,costabel1986principles,BEMbook}.

To solve \eqref{BVP_esse} using the so-called direct approch we rewrite the problem into the following BIE (the Symm's integral equation)
%
%
\begin{align}\label{eqn:BIEdirect}
 \int_{\Gamma} \mg(\x, \y) \phi (\y) \dtg =  \frac 1 2 u_D(\x) + \int_{\Gamma} \mh(\x, \y) u_D(\y) \dtg, \qquad \x\in\Gamma,
\end{align}
where the function $\phi$ denotes the unknown flux of $u$ and belongs to $H^{-1/2}(\Gamma)$, the dual space of $H^{1/2}(\Gamma)$ (duality is defined with respect to the usual $L^2(\Gamma)$-scalar product). 
The weakly singular kernels for the single and the double layer {\cc potentials} are
\begin{align}\label{eqn:kernelsGH}
\mg(\x,\y) = \frac{1}{4\pi} \frac{1}{\| \x-\y \|_2}, \qquad \mh(\x, \y)=\frac{\partial\mg(\x,\y)}{\partial \n_{\y}} = \frac{1}{4\pi} \frac{(\x - \y)^\top}{\| \x - \y \|_2^3} \n(\y),
\end{align}
and $\n$ denotes the outward unit normal vector.

Following the isoparametric approach, the boundary $\Gamma$ is represented as a mapped domain $\hat \Gamma$ in $\RR^2$,
\begin{align*}
\f:  \hat \Gamma  \to \RR^3, \; \t\mapsto \f(\t),
\end{align*}
where the geometry mapping $\f$ is described in terms of bivariate splines or NURBS \cite{librodeBoor01,Farin-02-CAGD,schumaker2007spline} in parametric coordinates $\t:=(t_1,t_2)$.
By applying the collocation discretization, we search for the approximate solution $\phi_h$ of \eqref{eqn:BIEdirect} in the discretization space, spanned by basis functions $B_j:\hat \Gamma \to \RR$, $j=1,2,\dots,N_{\rm DOF}$,
\begin{align*}
\phi_h(\x)  := \sum_{j=1}^{N_{\rm DOF}} \alpha_j (B_j \circ \f^{-1})(\x).
\end{align*}
By fixing suitable collocation sites $\x_i$, for $i=1,2,\dots, N_{\rm DOF}$, the unknown entries in the vector $\boldsymbol \alpha = (\alpha_1,\dots,\alpha_{N_{\rm DOF}})^{\top}$ are obtained by solving a square linear system,
\begin{align*}
A\boldsymbol{\alpha} = \boldsymbol{\beta}.
\end{align*}
The system matrix $A$ is fully populated and its entries $a_{ij}$   are
\begin{align*}
a_{ij} &= \frac{1}{4\pi} \int_{\Gamma} \frac{1}{\| \x_i-\y \|_2} (B_j \circ \f^{-1}) (\y) \dtg
= \frac{1}{4\pi} \int_{\hat \Gamma}   \frac{1}{\| \f(\s_i)-\f(\t) \|_2} B_{j}(\t)
J(\t) \dt,
\end{align*}
where $J := \| D_1 \f \times D_2 \f  \|_2$ 
is the area of the infinitesimal surface element, $D_1:=\partial \f /{\partial t_1}, D_2:=\partial \f /{\partial t_2}$, and $\f(\s_i)=\x_i$, $\f(\t)=\y$. 
Entries of the right-hand side vector $\boldsymbol \beta$ are
\begin{align*}
\beta_i &= \frac 1 2  (u_D \circ \f) (\s_i)
+ \frac 1 {4 \pi} \int_{\hat \Gamma}  \mh (\f(\s_i), \f(\t))\, J(\t)\,  (u_D \circ \f) (\t) \dt.
\end{align*}
%
%
%
%
%
Integral in $a_{ij}$ is classified singular if the source point $\s_i$ is inside the support of $B_j$. If $\s_i$ is outside the support but close to it, we declare the integral nearly singular. In the remaining case we say the integral is regular.

\subsection{Series expansion of a singular kernel}\label{kernelSeries}

In this subsection we derive a series expansion of kernels $\mg$ and $\mh$ for fixed source point $\s$. By restricting the kernel to a line segment, the singular part of the kernel is easily decoupled from the regular part -- for the latter part we then apply the Taylor series expansion. The procedure shares some similarity to polar and Duffy transformation and to radial integration approach by Gao \cite{GAO2002905}, with the difference that in our case the final expression of the series is written again in the starting intrinsic coordinates. Each summand in the series is a function of a type, described in Section~\ref{sec:recursiveFormulae}, thus a closed form expression for its integral exists. 

\subsubsection{Series expansion of $\mg$}

As a proof of concept we focus on the simplest kernel, $\mg$ in \eqref{eqn:kernelsGH}.
Let us write $\x = \f(\s)$, $\y = \f(\t)$ and let us assume $\f$ to be analytical (or sufficiently smooth) near fixed $\x = \f(\s)$. For a fixed $\s$ we can express kernel $\mg$, $\mg_\s (\z) := \mg(\f(\s), \f(\s-\z))$, in local coordinates $\z := (z_1, z_2) := \s - \t$. Taylor expansion of $\f$ about $\s$ can be compactly written as
\begin{align*}
\f(\t) =  \f(\s) + \sum_{i=1}^\infty \d_i \f (\t),
\qquad \d_i \f (\t):= \sum_{|\alpha|=i} \a_\alpha \z^{\alpha},
\qquad \a_\alpha :=  \frac{(-1)^{|\alpha|}}{\alpha!} D^{\alpha} \f(\s).
\end{align*}
%
%
Here we adopt the multi-index notation. Specifically, $\alpha = (\alpha_1, \alpha_2)$, $\alpha_1, \alpha_2 \in \mathbb N_0$ and $|\alpha|=\alpha_1+\alpha_2$, $\alpha!=\alpha_1! \alpha_2!$. The partial derivative operator $D^\alpha$  acts on each component of a vector function $\f$ separately and on both variables,
\begin{align*}
D^\alpha = D_1^{\alpha_1} D_2^{\alpha_2} = \frac{\partial^{|\alpha|}}{\partial t_1^{\alpha_1} \partial t_2^{\alpha_2}}.
\end{align*}
Also we compactly write the scalar $\z^\alpha = z_1^{\alpha_1} z_2^{\alpha_2}$. 

From the expression
\begin{align*}
\f(\s) -  \f(\t) = - \sum_{|\beta| \geq 1} \a_\beta \z^{\beta}
\end{align*}
we can easily express the kernel
$\mg_{\s}$ as a power series of $\z$ about $(0,0)$,
\begin{align}\label{eqn:g_expansion}
4\pi \mg(\f(\s),\f(\t)) &= \frac{1}{\| \f(\s) - \f(\t) \|_2} 
= \left( \sum_{{|\alpha|\geq 2}} c_{\alpha} \z^{\alpha}  \right)^{-\frac{1}{2}},
\end{align}
where 
\begin{align}\label{eqn:coefs_c}
c_{\alpha} = \sum_{\substack{\beta + \gamma = \alpha\\ |\beta|\geq1, |\gamma|\geq1}} \a_\beta^\top \a_\gamma
\end{align}
are the coefficients in the series expansion of $\|\f(\s) - \f(\t)\|_2^2$ about $\s$ in local variable $\z$. From  \eqref{eqn:coefs_c} it is easy to check that coefficient $c_{\alpha}$ is a linear combination of products of derivatives of $\f$, where the total degree of the derivatives is at most $|\alpha|-1$.

In an ideal world one would probably just truncate the series in \eqref{eqn:g_expansion} and integrate the obtained function, involving the Taylor polynomial. Unfortunately, closed forms for these types of integrals are not known in general, therefore a further transformation of the expression is needed. In the next step we derive a Taylor expansion of the regular part of  \eqref{eqn:g_expansion}, restricted to an arbitrary line $(z_1, v z_1)$, $v\in\mathbb R$. From the parameter $v$ depended linear expansion it is relatively straightforward to derive a truncated series expansion of $\mg_\s$ that can be integrated analytically.

Let $L_v :=\{(z_1,v z_1), z_1\in (-\delta, \delta) \}$ be a line segment passing the origin $(0,0)$, where $v\in \mathbb R$ is a fixed parameter and $\delta$ is sufficiently small.
With $\mg_v$ let us write the kernel $\mg$, restricted on $L_v$ and written as a function of $z_1$.
By multiplying the expression \eqref{eqn:g_expansion} with $|z_1|$ we obtain
\begin{align}\label{eqn:rv}
r_v(z_1) := |z_1|\, 4\pi \mg_{v}(z_1) = |z_1| \left( \sum_{|\alpha|\geq 2} c_{\alpha} z_1^{\alpha_1} (v z_1)^{\alpha_2}  \right)^{-\frac{1}{2}}
 = \left( \sum_{|\alpha|\geq 2} c_{\alpha} v^{\alpha_2} z_1^{|\alpha|-2}   \right)^{-\frac{1}{2}},
\end{align}
and the newly defined function $r_v$ is regular (it holds $  \sum_{|\alpha|= 2} c_{\alpha} v^{\alpha_2}  >0$ due to the first fundamental form at a regular point on the surface).
That way we separate the singular part $|z_1|^{-1}$ of the kernel $\mg_{v}$ from the regular part $r_v(z_1)$ ,
\begin{align}\label{eqn:1DApprox}
  4\pi \mg_{v}(z_1)  = |z_1|^{-1} r_v(z_1).
\end{align}
 Since $r_v$ is a regular function in $z_1$, we can replace it with its Taylor series expansion in $z_1$ about origin 0 using expression \eqref{eqn:rv},
\begin{align}\label{eqn:rvTaylor}
r_v(z_1)  = \sum_{\ell=0}^\infty \frac{r_v^{(\ell)}(0)}{\ell!} z_1^\ell  
 = \big( P_{v,2} \big)^{-\frac{1}{2}}
-  \frac 1 2 \big( P_{v,2} \big)^{-\frac{3}{2}}  P_{v,3}\, z_1
+ \bigg( -  \big( P_{v,2} \big)^{-\frac{3}{2}}  P_{v,4}
+ \frac{3}{4} \big( P_{v,2} \big)^{-\frac{5}{2}}  \big( P_{v,3} \big)^2 \bigg) \frac{z_1^2}{2} + \dots
\end{align}
where $P_{v,\ell}:=\sum_{|\alpha|=\ell} c_{\alpha} v^{\alpha_2}$ is a polynomial of degree $\ell$ in $v$. 

The obtained expression for $4\pi \mg_{v}$ is valid for any $v$, we insert back $\displaystyle{v = \frac{z_2}{z_1}}$ into \eqref{eqn:1DApprox} using the expansion \eqref{eqn:rvTaylor}
 and get the sought series expansion of $\mg_\s$,
\begin{align}\label{eqn:Gseries}
\mg_\s (\z)  = \frac{1}{4\pi} \left[ R(\z)^{-1}  - \frac{1}{2}  R(\z)^{-3}
P_3(\z) 
+ \left( -  \frac{1}{2} R(\z)^{-3} P_4(\z) + \frac{3}{8} R(\z) ^{-{5}} P_3(\z)^2  \right)  + \dots \right],
\end{align}
where $R(\z):=R(c_{(2,0)},c_{(1,1)},c_{(0,2)},\z) = (\sum_{|\alpha|=2} c_{\alpha}  \z^{\alpha})^{1/2}$, $P_\ell (\z):=\sum_{|\alpha|=\ell} c_{\alpha}  \z^{\alpha}$ is the bivariate homogeneous polynomial of degree $\ell$, and coefficients $c_{\alpha}$ are defined in \eqref{eqn:coefs_c}. Observe that in \eqref{eqn:Gseries} there is no problem in the definition when $z_1=0$ and $z_2 \neq 0$; it is only a limitation of the definition of $\mg_v$.

For the sake of easier notation and vocabulary, we consider an $\ell$-th term of \eqref{eqn:Gseries} a weighted sum of all entries $R^p z_1^q z_2^r$, that satisfy $\zeta(R^p z_1^q z_2^r) = \ell-2$. Indeed, compactly we could express \eqref{eqn:Gseries} as
\begin{align}\label{eqn:Gseries_v2}
\mg_\s  = \frac{1}{4\pi} \sum_{\ell=1,2,\dots} R^{-2 \ell+1} P_{3\ell-3}^{[\ell]},
\end{align}
where $P_{3\ell-3}^{[\ell]}$ are appropriate homogeneous polynomials of degree $3\ell-3$ and clearly for the $\ell$-th term it holds $\zeta (R^{-2 \ell+1} P_{3\ell-3}^{[\ell]}) = \ell-2$. In this notation we assume that there are no common factors between polynomials $R^2$ and $P_{3\ell-3}^{[\ell]}$.

With $\mg_{\s,n}$ we denote a kernel, obtained from $\mg_\s$ by truncating the series in \eqref{eqn:Gseries} after the $n$-th term. The kernel $\mg_{\s,n}$ is a sum of $n$ addends and double integral of each of the terms is a sum of functions of a type $I_{p,q,r}(z_1,z_2)$, defined in \eqref{eqn:fullInt}. As we {\cc show} in Section~\ref{sec:recursiveFormulae}, these type of integrals have closed form expressions.


\begin{rmk}
If the geometry is locally flat and parameterized by a linear map, then the cancellation of the singularity with the first term in the series expansion is exact, $\rho:= \mg_\s - \mg_{\s,1} \equiv 0$ and $\rho:= \mg_\s / \mg_{\s,1} \equiv 1$ , for singularity subtraction and division, respectively (see Section~\ref{sec:num_int} for more details about utilization of the singularity extraction in numerical integration). Thus, the procedure generalizes the common singularity extraction formulae for flat surfaces. 
\end{rmk}

\bigskip


\subsubsection{Series expansion of $\mh$}

We proceed similarly as for the kernel $\mg$. Since we want to express the kernel $\mh$ in the intrinsic coordinates from the information of the geometry mapping $\f$,  it is more convenient to consider the kernel $\bar \mh := \mh J$, where $J(\t) = \| D_1 \f(\t) \times D_2 \f (\t) \|_2$,
\begin{align}\label{eqn:kernelH}
\bar \mh(\f(\s), \f(\t)) = \mh(\f(\s), \f(\t))\, J(\t) = \left( \frac{\f(\s)-\f(\t)}{\| \f(\s)-\f(\t) \|^3} \right)^\top (D_1 \f(\t) \times D_2 \f (\t)).
\end{align}
The inner product of the two vectors in \eqref{eqn:kernelH} can be compactly written as a linear combination of its components. Therefore
\begin{align}\label{eqn:kernelH2}
\bar \mh(\f(\s), \f(\t))  
= \frac{\varepsilon_{\ell_1 \ell_2 \ell_3} \big( \f(\s)-\f(\t) \big)_{\ell_1} \big( D_1 \f(\t) \big)_{\ell_2} \big( D_2 \f (\t) \big)_{\ell_3}}{\| \f(\s)-\f(\t) \|^3},
\end{align}
where we use Levi-Civita symbol $\varepsilon$, components of vectors $(\bullet)_{\ell_i}$, for ${\ell_i}=1,2,3$, and Einstein summation convention.

By fixing the  source points $\s$ we define $\bar \mh_\s (\z) := \bar \mh(\f(\s), \f(\s-\z))$ and the geometric quantities in the kernel can be written in the Taylor series about $\s$ in the local variable $\z := \s - \t$,
\begin{align*}
\f(\s) -  \f(\t) &= - \sum_{|\beta | \geq 1} \a_\beta \z^{\beta},\\
D_1 \f(\t) &=  D_1 \bigg(\f(\s) + \sum_{|\gamma| \geq 1} \a_\gamma \z^{\gamma} \bigg)
=- \sum_{|\gamma| \geq 0} (\gamma_1+1)  \a_{(\gamma_1+1,\gamma_2)} \z^{\gamma},\\
D_2 \f(\t) &=  D_2 \bigg(\f(\s) + \sum_{|\delta| \geq 1} \a_\delta \z^{\delta} \bigg)
=- \sum_{|\delta| \geq 0} (\delta_2+1)  \a_{(\delta_1,\delta_2+1)} \z^{\delta},
\end{align*}
%
%
%
%
%
%
where $\gamma= (\gamma_1,\gamma_2)$, $\delta= (\delta_1,\delta_2)$. Therefore the kernel $\bar \mh_\s$ can be compactly written as
\begin{align*}
\bar \mh_\s = \frac{\displaystyle \sum_{{|\alpha|\geq 2}} d_{\alpha} \z^{\alpha} }
{\displaystyle \left( \sum_{{|\alpha|\geq 2}} c_{\alpha} \z^{\alpha} \right)^{\frac{3}{2}}},
\end{align*}
where $c_{\alpha}$ are defined as in \eqref{eqn:coefs_c} and
\begin{align*}
d_{\alpha} = - \sum_{\substack{\beta + \gamma + \delta = \alpha\\ |\beta|\geq1, |\gamma|\geq0, |\delta|\geq0}} \varepsilon_{\ell_1 \ell_2 \ell_3} \big( \a_\beta \big)_{\ell_1} \big( (\gamma_1+1) \a_{(\gamma_1+1,\gamma_2)} \big)_{\ell_2} \big( (\delta_2+1)\a_{(\delta_1,\delta_2+1)} \big)_{\ell_3}.
\end{align*}
Note that $d_\alpha=0$ for $|\alpha|=1$; namely
\begin{align*}
d_{(1,0)} &= \varepsilon_{\ell_1 \ell_2 \ell_3} \big( \a_{(1,0)} \big)_{\ell_1} \big( \a_{(1,0)} \big)_{\ell_2} \big( \a_{(0,1)} \big)_{\ell_3}= 
\left|\begin{matrix}
\a_{(1,0)} & \a_{(1,0)} & \a_{(0,1)}
\end{matrix}\right|
= 0,\\
d_{(0,1)} &= \varepsilon_{\ell_1 \ell_2 \ell_3} \big( \a_{(0,1)} \big)_{\ell_1} \big( \a_{(1,0)} \big)_{\ell_2} \big( \a_{(0,1)} \big)_{\ell_3}=
\left|\begin{matrix}
\a_{(0,1)} & \a_{(1,0)} & \a_{(0,1)}
\end{matrix}\right|
=0,
\end{align*}
since in both cases the involved three vectors are linearly dependent.

Again, let $L_v :=\{(z_1,v z_1), z_1\in (-\delta, \delta) \}$ for fixed $v\in \mathbb R$, $\delta$ is sufficiently small and let $\bar \mh_v$ be the kernel $\bar \mh$, restricted on $L_v$ and written as a function of $z_1$.
Then
\begin{align}\label{eqn:rv2}
r_v(z_1) := |z_1|\,  4\pi \bar \mh_{v}(z_1) =
 |z_1|  \frac{\displaystyle \sum_{{|\alpha|\geq 2}} d_{\alpha} z_1^{\alpha_1} (v z_1)^{\alpha_2} }
{\displaystyle \left( \sum_{{|\alpha|\geq 2}} c_{\alpha} z_1^{\alpha_1} (v z_1)^{\alpha_2}  \right)^{\frac{3}{2}}}=
%
 \frac{ \displaystyle \sum_{{|\alpha|\geq 2}} d_{\alpha} v^{\alpha_2} z_1^{|\alpha|-2} }
{ \displaystyle \left( \sum_{{|\alpha|\geq 2}} c_{\alpha} v^{\alpha_2}  z_1^{|\alpha|-2}  \right)^{\frac{3}{2}}}.
\end{align}
Its Taylor expansion reads
\begin{align}\label{eqn:rvTaylor2}
r_v(z_1) & = \sum_{\ell=0}^\infty \frac{r_v^{(\ell)}(0)}{\ell!} z_1^\ell  \nonumber\\
& =   \big( P_{v,2} \big)^{-\frac{3}{2}} Q_{v,2}
+ \left(  \big( P_{v,2} \big)^{-\frac{3}{2}} Q_{v,3}  -\frac{3}{2} \big( P_{v,2} \big)^{-\frac{5}{2}}  Q_{v,2}  P_{v,3} \right) \, z_1 \nonumber \\
&+ \left( 2 ( P_{v,2} )^{-\frac{3}{2}} Q_{v,4}  - 3 ( P_{v,2} )^{-\frac{5}{2}}  Q_{v,3}  P_{v,3}  
- 3 \big( P_{v,2} \big)^{-\frac{5}{2}} Q_{v,2}  P_{v,4} 
+ \frac{15}{4} \big( P_{v,2} \big)^{-\frac{7}{2}} Q_{v,2}  (P_{v,3})^2 
\right)  \frac{1}{2}z_1^2 + \dots
\end{align}
where $P_{v,\ell}=\sum_{|\alpha|=\ell} c_{\alpha} v^{\alpha_2}$, $Q_{v,\ell}:=\sum_{|\alpha|=\ell} d_{\alpha} v^{\alpha_2}$ are polynomials of degree $\ell$ in $v$.

Since the obtained expression for $4\pi \bar \mh_{v}$ is valid for any $v$, 
\begin{align}\label{eqn:1DApprox2}
  4\pi \bar \mh_{v}(z_1)  = |z_1|^{-1} r_v(z_1),
\end{align}
we insert back $\displaystyle{v = \frac{z_2}{z_1}}$ into \eqref{eqn:1DApprox2} using the expansion \eqref{eqn:rvTaylor2}
 and get the sought series expansion of $\bar \mh_\s$,
\begin{align}\label{eqn:Hseries}
\bar \mh_\s (\z)  = &\frac{1}{4\pi} \Bigg[ R(\z)^{-3} Q_2(\z)  
 + \left( R(\z)^{-3} Q_3(\z)  -\frac{3}{2} R(\z)^{-5} Q_{2}(\z)   P_{3}(\z)   \right) \nonumber \\
&+ \left( R(\z)^{-3} Q_4(\z) - \frac{3}{2} R(\z)^{-5} Q_3(\z) P_3(\z)  
- \frac{3}{2} R(\z)^{-5} Q_2(\z)  P_4(\z)  + \frac{15}{8}  R(\z)^{-7} Q_2(\z)  P_{3}(\z)^2
\right) + \dots \Bigg],
\end{align}
where 
$P_\ell (\z)=\sum_{|\alpha|=\ell} c_{\alpha} \z^{\alpha}$, $Q_\ell (\z):=\sum_{|\alpha|=\ell} d_{\alpha} \z^{\alpha}$ are two homogeneous polynomials.
Compactly we can write \eqref{eqn:Hseries} as
\begin{align}\label{eqn:Hseries_v2}
\bar \mh_\s  = \frac{1}{4\pi} \sum_{\ell=1,2,\dots} R^{-2 \ell-1} P_{3\ell-1}^{[\ell]},
\end{align}
where $P_{3\ell-1}^{[\ell]}$ are appropriate homogeneous polynomials of degree $3\ell-1$, there are no common factors between $R^2$ and $P_{3\ell-1}^{[\ell]}$ and the continuity of the $\ell$-th term is characterized by $\zeta (R^{-2 \ell-1} P_{3\ell-1}^{[\ell]}) = \ell-2$.

\subsection{Numerical integration scheme} \label{sec:num_int}

By utilizing the kernel expansion from the previous subsection we can now focus on numerical integration of integrals in the Symm's integral equation \eqref{eqn:BIEdirect}.
In this section we consider a more general singular kernel $\mk_\s$ that can be written as a similar series expansion and two regularization approaches are introduced.

%

Let $\mk_\s = \mk_\s(\z)$ be a singular kernel, expressed about a source point $\s$ in local intrinsic coordinates $\z$,
\begin{align}\label{eqn:kernel}
\mk_\s  = 
\sum_{\ell=1,2,\dots} R^{-2 (\ell+m_1)+1} P_{3\ell +m_2}^{[\ell]},
\end{align}
with regularity $m := \zeta(\mk_\s) = -2m_1+m_2+2$ and integers $m_2 \geq -3$, $m_1 \geq 0$. We assume $m<0$ and smaller $m$ corresponds to a stronger type of singularity. 
Here $P_{3\ell +m_2}^{[\ell]}$ are again suitable homogeneous polynomials of degree $3\ell+m_2$. 
Additionally, let $\mk_\s$ be sufficiently smooth for $\z\in D_{\varepsilon}$, where
\begin{align*}
D_\varepsilon := \{ \z: 0 < \| \z \|_2 \leq \varepsilon\}
\end{align*}
is a punctured disk for some fixed $\varepsilon>0$. 
Let  $\mk_{\s,n}$ be an approximation of $\mk_\s$ by truncating the infinite series \eqref{eqn:kernel} after the $n$-th term,
\begin{align}\label{eqn:kernelK}
\mk_{\s,n}  = 
\sum_{\ell=1,2,\dots,n} R^{-2 (\ell+m_1)+1} P_{3\ell +m_2}^{[\ell]}.
\end{align}

The singular kernel $\mk_\s$ is either regularized by the subtraction, i.e., $\rho := \mk_\s - \mk_{\s,n}$, or the division, i.e., $\rho := \mk_\s / \mk_{\s,n}$. Our goal is to numerically compute integrals 
\begin{align}\label{eqn:initInt}
\frac{1}{4\pi} \int_{I_1} \int_{I_2}  \mk_\s(\z)  B_{j}(\s - \z) v(\s - \z) \dz,
\end{align}
where $B_{j}$ is a polynomial basis function and $v$ is an auxiliary smooth function that can include for example the area of the infinitesimal surface element $J$ or a boundary datum.

\begin{enumerate}
\item[Case 1:]
If the integral \eqref{eqn:initInt} is regular, we can apply one of the common numerical integration schemes since the integrand is a smooth function. 

\item[Case 2:]
If the integral \eqref{eqn:initInt} is singular or nearly singular, the following transformation needs to be applied first.

\begin{enumerate}
\item
The first strategy, which is more common in literature, is the singularity subtraction. By writing $\rho := \mk_\s - \mk_{\s,n}$, the integral \eqref{eqn:initInt} is expressed as a sum of two integrals
\begin{align}\label{eqn:subInt}
\frac{1}{4\pi} \int_{I_1} \int_{I_2}  \rho(\z)  B_{j}(\s - \z) v(\s - \z) + 
\frac{1}{4\pi} \int_{I_1} \int_{I_2}  \mk_{\s,n}(\z)  B_{j}(\s - \z) v(\s - \z) \dz.
\end{align}
The integrand in the first integral in \eqref{eqn:subInt} is sufficiently regular if the geometry parameterization (hence $v$) and $\rho$ are sufficiently regular. The smoothness of $\rho$ is studied in the next subsection. Thus, for this type of integrals we can apply a quadrature scheme for regular integrals. In the second integral in  \eqref{eqn:subInt} the regular part $B_{j}(\s - \z) v(\s - \z)$ can be approximated sufficiently well with an appropriate polynomial $P$ via least square fitting, quasi-interpolation etc.

\item
In the regularization process via division step we consider $\rho := \mk_\s / \mk_{\s,n}$ and express the integral \eqref{eqn:initInt} as 
\begin{align}\label{eqn:divInt}
\frac{1}{4\pi} \int_{I_1} \int_{I_2}  \mk_{\s,n}(\z)  B_{j}(\s - \z) v(\s - \z) \rho(\z) \dz.
\end{align}
The regular part of the integrand, $B_{j}(\s - \z) v(\s - \z) \rho(\z)$, is again approximated with a suitable polynomial $P$.
\end{enumerate}

In both strategies (a), (b), the singular part of the integral \eqref{eqn:initInt} is approximated with
\begin{align}\label{eqn:simplifiedInt}
\frac{1}{4\pi} \int_{I_1} \int_{I_2}  \mk_{\s,n}(\z) P(\z) \dz
\end{align}
for a suitable polynomial $P$. The integral \eqref{eqn:simplifiedInt} can be expressed as a sum of fundamental integrals, presented in Section~\ref{sec:recursiveFormulae}.

\end{enumerate}

\begin{rmk}
In \eqref{eqn:initInt}, instead of a polynomial basis we can proceed similarly if we consider $B_j$ to play a role of a polynomial spline basis function. In that case we can decompose the basis function into a sum of polynomial contributions (for example using a B\'ezier extraction technique) and from thereon follow a similar procedure as before.

Alternatively, we can approximate the smooth part of the integrand of \eqref{eqn:initInt} with a spline function and write the product of the two splines as a new spline function, for example by using formulae in \cite{Morken91} for the tensor product basis functions. For a successful implementation of this idea to solve Laplace and Helmholtz boundary integral equations in BEM see \cite{FGKSS_2021,nash20,FKSS_2021}. The model combines a spline quasi-interpolation operator \cite{MSbit09, MSJcam12}, together with a spline recurrence relation formula \cite{FKSS_2021} to compute the integrals \eqref{eqn:simplifiedInt}, where $P$ plays a role of a polynomial spline function.
\end{rmk}

\subsection{Implementation speedup}

The presented singularity extraction technique can only be considered useful in practice if it can be efficiently implemented in software. For example, to perform a simulation using BEM, the matrix formation can induce a computation of thousands or millions of integrals of the type  \eqref{eqn:simplifiedInt} and it can represent the main bottleneck in the system matrix formation. However, since the evaluation of these integrals for different source points is embarrassingly parallelizable, the time complexity of this step can be partially alleviated. In this subsection we propose an additional reduction of the computational complexity by accessing specific precomputed integrals from lookup tables.

By exchanging the order of summation and integration in \eqref{eqn:simplifiedInt}, and writing each polynomial $P$ as a weighted sum of monomials (the weights depend on the local geometry around the point $\s$ and can be efficiently computed in a preprocess step), we are left to analyze how to efficiently integrate functions $R^p x^q y^r$. Unfortunately, analytical formulae for the indefinite integrals \eqref{eqn:fullInt} can have very long symbolic expressions that can result in undesired lengthy evaluation, especially for higher values of $-p+q+r$. To overcome this problem, we apply several simplifications of the expressions to reduce the number of parameters, influencing the integral values, and precompute values of particular simpler integrals. We can later on {\cc access} these values in the stored lookup tables to compute the definite integrals.

\subsubsection{Rectangles}
The goal is to derive to an efficient procedure to compute the building blocks of the definite integral of \eqref{eqn:kernelK},
\begin{align*}
I := \int_{y_0}^{y_1} \int_{x_0}^{x_1} R^p x^q y^r \dx \dy,
\end{align*}
on aligned rectangular domains $D = [x_0, x_1]\times[y_0, y_1]$, where $R(a,b,c,(x,y)) = \sqrt{a x^2 + b x y + c y^2}$, as in \eqref{eqn:fullInt}. Unfortunately, the integrals depend on too many ``continuous'' variables (e.g., $a,b,c,x_0,x_1,y_0,y_1$) to be storable in lookup tables for a dense enough grid of values of variables. Therefore, we need to simplify the integrals by exploiting all the correlations between the variables.

\begin{figure}[t!]
\centering
\subfigure[Interior source point: \newline ``$D = D_1 + D_2 + D_3 + D_4$'']{
\begin{tikzpicture}[->,>=stealth',auto,node distance=2cm,
  main node/.style={circle,draw,font=\sffamily\Large\bfseries}]
\node[inner sep=0pt] (IG) at (0,0)
 {\includegraphics[trim = 6cm 11.5cm 5.5cm 10.5cm, clip = true, height=3cm]{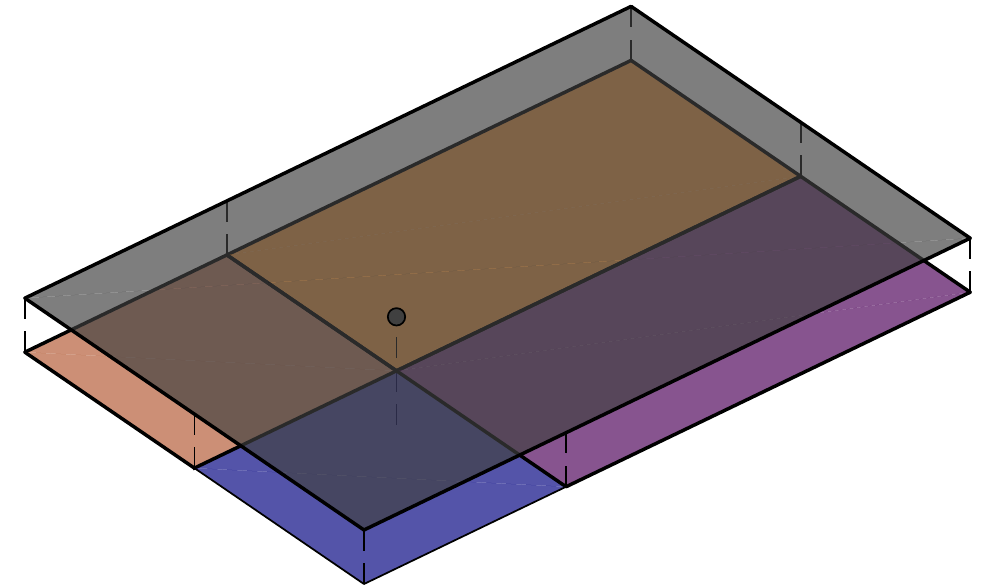}};

\draw[-] (-2.25,.5) -- (-1.3,.5);
\draw[->] [dashed](-1.3,.5) -- (0,.5);
\node[text width=0cm] at (-2.75,.5) {$D_3$};

\draw[-] (-2.25,-.5) -- (-1.7,-.5);
\draw[->] [dashed](-1.7,-.5) -- (-1.3,-.5);
\node[text width=0cm] at (-2.75,-.5) {$D_4$};

\draw[-] (-2.25,-1) -- (-0.85,-1);
\draw[->][dashed] (-.8,-1) -- (-0.5,-1);
\node[text width=0cm] at (-2.75,-1) {$D_1$};

\draw[->] (2.,-.5) -- (1.,-.5);
\node[text width=0cm] at (2.1,-.5) {$D_2$};
\end{tikzpicture}
}
\subfigure[Edge source point: \newline ``$D = D_1 + D_2$'']{
\begin{tikzpicture}[->,>=stealth',auto,node distance=2cm,
  main node/.style={circle,draw,font=\sffamily\Large\bfseries}]
\node[inner sep=0pt] (IG) at (0,0)
{\includegraphics[trim = 6cm 11.5cm 5.5cm 10.5cm, clip = true, height=3cm]{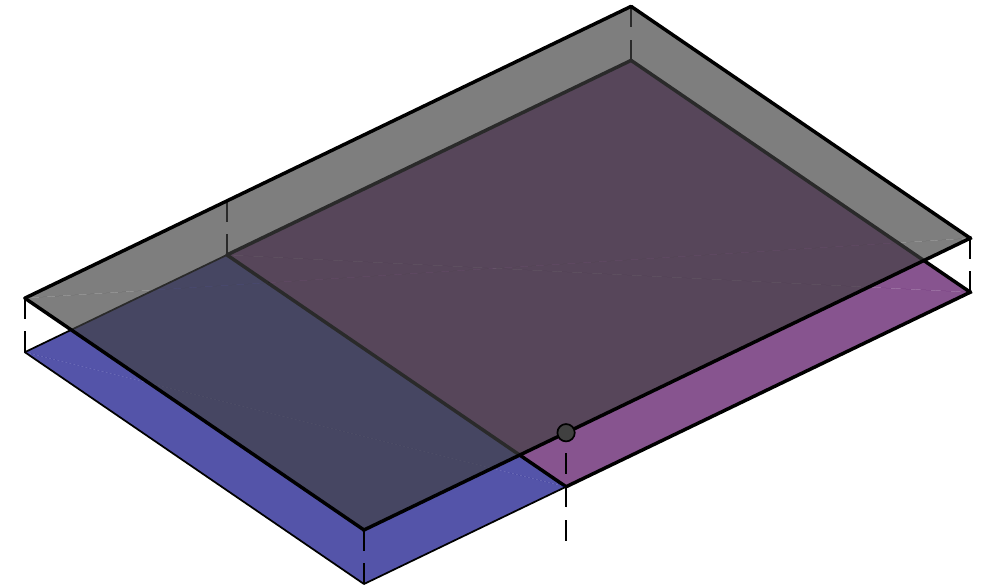}};


\draw[-] (-2.25,-.75) -- (-1.4,-.75);
\draw[->] [dashed](-1.4,-.75) -- (-0.75,-.75);
\node[text width=0cm] at (-2.75,-.75) {$D_1$};

\draw[->] (2.,-.5) -- (1.,-.5);
\node[text width=0cm] at (2.1,-.5) {$D_2$};
\end{tikzpicture}
}
\subfigure[Corner source point]{
{\includegraphics[trim = 6cm 11.cm 5.5cm 10.75cm, clip = true, height=3cm]{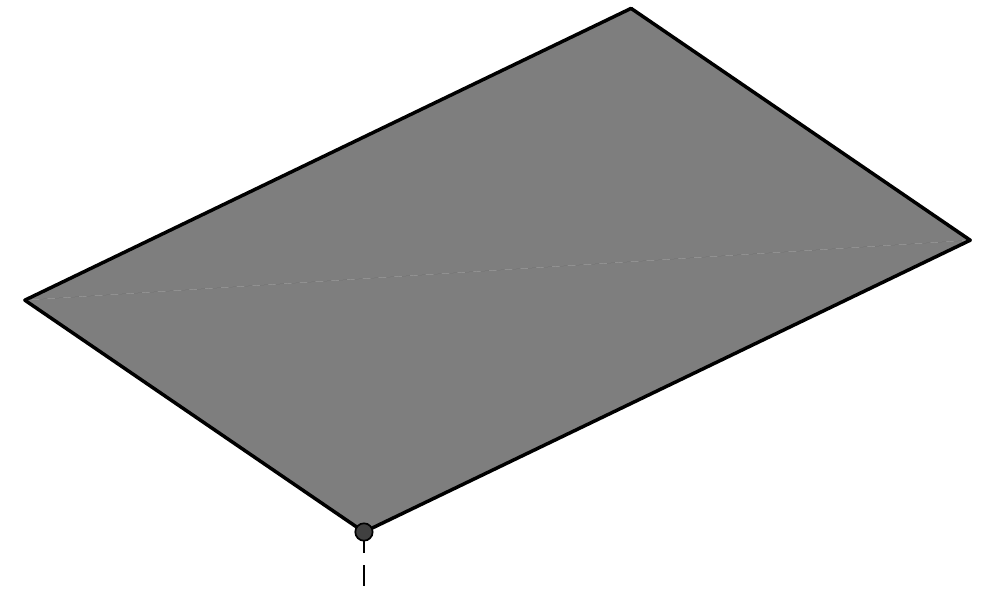}}
}
\subfigure[Aligned exterior source point: \newline ``$D = D_1 - D_2$'']{
\begin{tikzpicture}[->,>=stealth',auto,node distance=2cm,
  main node/.style={circle,draw,font=\sffamily\Large\bfseries}]
\node[inner sep=0pt] (IG) at (0,0)
{\includegraphics[trim = 6cm 10.5cm 5.5cm 10.25cm, clip = true, height=3cm]{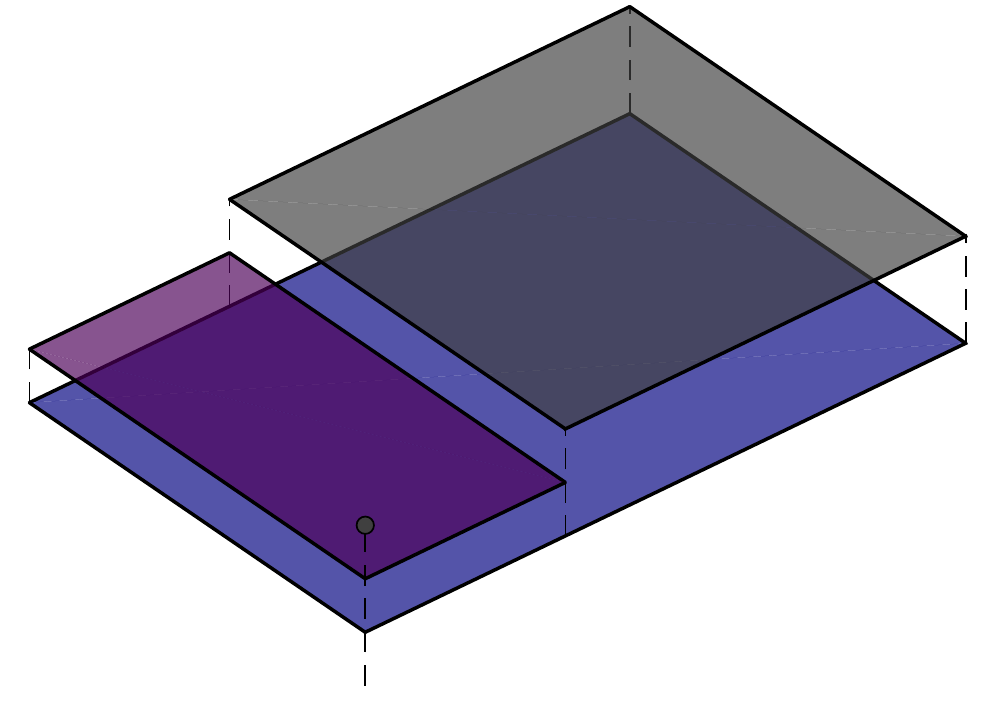}};


\draw[->] (-1.9,-.5) -- (-0.75,-.5);
\node[text width=0cm] at (-2.4,-.5) {$D_2$};

\draw[->] (1.7,-.4) -- (0.6,-.4);
\node[text width=0cm] at (1.8,-.4) {$D_1$};
\end{tikzpicture}
}
\subfigure[General exterior source point: \newline ``$D = D_1 - D_2 - D_3 + D_4$'']{
\begin{tikzpicture}[->,>=stealth',auto,node distance=2cm,
  main node/.style={circle,draw,font=\sffamily\Large\bfseries}]
\node[inner sep=0pt] (IG) at (0,0)
{\includegraphics[trim = 6cm 10.7cm 5.5cm 10.25cm, clip = true, height=3cm]{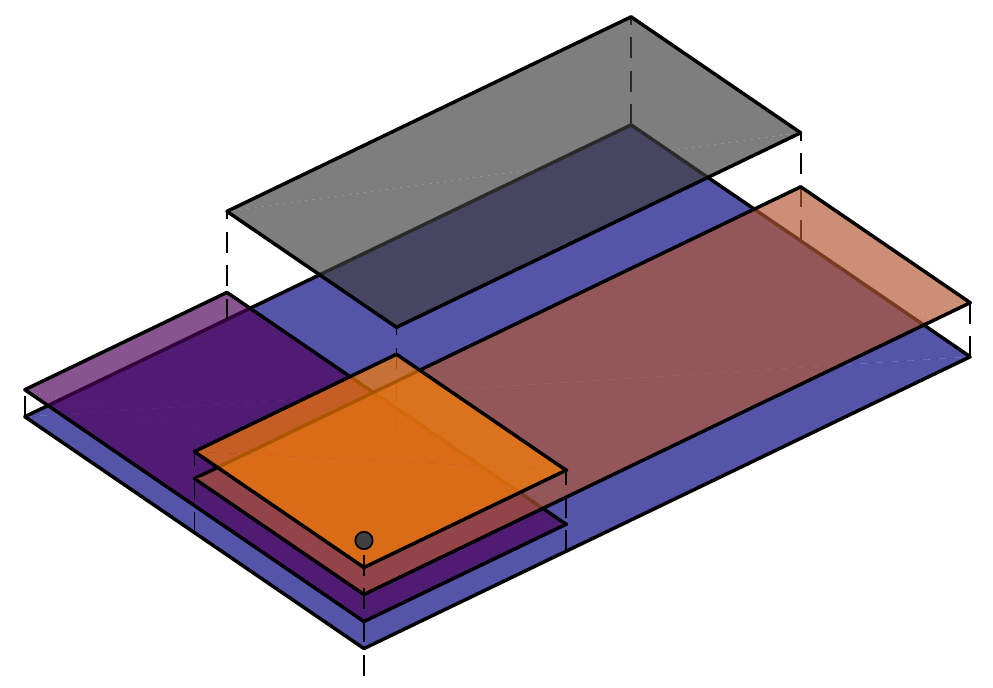}};

\draw[-] (-2.25,.35) -- (-0.84,.35);
\draw[-][dashed] (-0.8,.35) -- (0.1,.35);
\draw[->] (.15,.35) -- (1,.35);
\node[text width=0cm] at (-2.75,.35) {$D_3$};

\draw[->] (-2.25,-.25) -- (-1.5,-.25);
\node[text width=0cm] at (-2.75,-.25) {$D_2$};

\draw[->] (-2.25,-.75) -- (-0.5,-.75);
\node[text width=0cm] at (-2.75,-.75) {$D_4$};

\draw[->] (1.8,-.5) -- (.9,-.5);
\node[text width=0cm] at (1.9,-.5) {$D_1$};
\end{tikzpicture}
}

\caption{All possible positions of the rectangle $D$ (in gray) with respect to the source point $(0,0)$, depicted as a gray dot. Sub-rectangles are shown in different colours and for better visualization the overlapping ones are depicted at different vertical layers.}
\label{fig:Rxy_cases}
\end{figure}
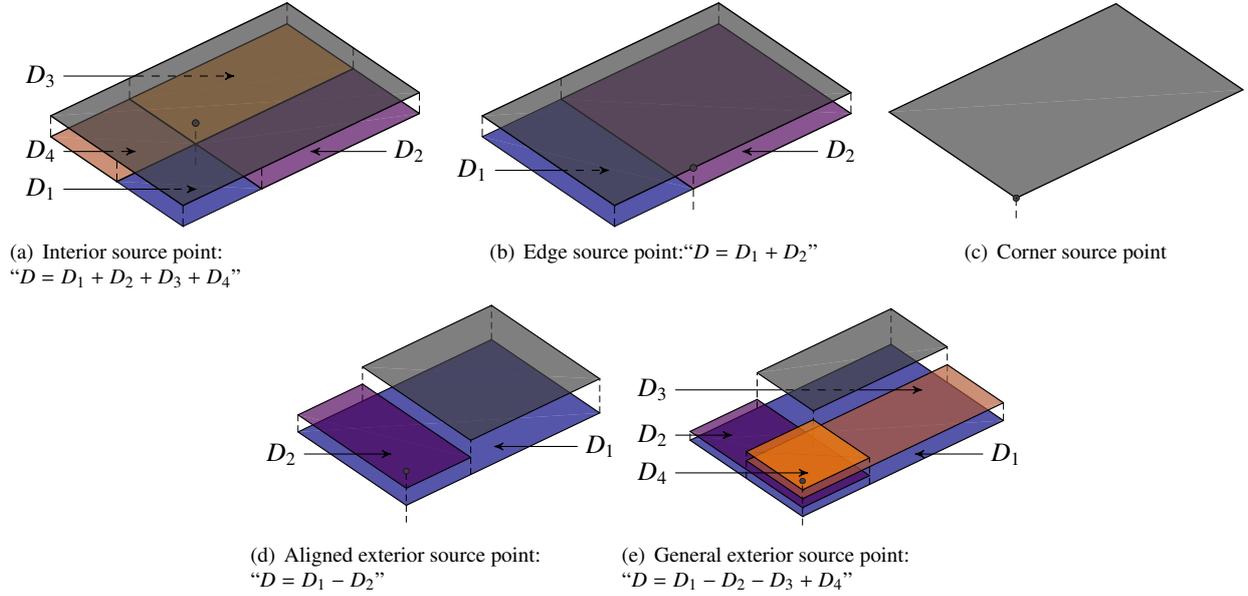

First, let us write $I$ as a linear combination of at most 4 integrals, called $I_0$, which contain the same integrand as $I$, but defined on rectangles with $(0,0)$ as one of their corners. In Figure~\ref{fig:Rxy_cases} we can see 5 different possible cases of the initial rectangle with respect to the source point $(0,0)$. In (a) and (b), the initial rectangle $D$ is split into 4 and 2 smaller rectangles, respectively. Case (c) is a special case, when $I=I_0$ and no treatment is needed. In (d),  the source point is outside the initial rectangle but one of rectangle's edge coordinate is 0 -- rectangle $D$ is extended to $D_1$ and then $D_2$ is subtracted from the augmented one. In (e), $D_2$ and $D_3$ are subtracted from $D_1$, and their intersection $D_4$ is added back.

For each $I_{0}$, defined on a rectangle $[0, x_2]\times[0, y_2]$, we apply a linear transformation to map  the rectangle to the unit square $[0,1]^2$. Then we divide $R$ by ${a}^{p/2}$. Therefore
\begin{align*}
I_0 =  x_2^{p+q+1} y_2^{r+1} a^{p/2} \int_{0}^{1} \int_{0}^{1} R \left(1, \frac{b}{a}\frac{y_2}{x_2}, \frac{c}{a} \frac{y_2^2}{x_2^2}, (x,y) \right)^p x^q y^r \dx \dy. 
\end{align*}

Observe that the integral in $I_0$ depends only on two continuous variables $\bar b:=({b}/{a}) ({y_2}/{x_2})$, $\bar c:=({c}/{a}) ({y_2}/{x_2})^2$ and on integers $p,q,r$. For every needed combination of integers $p,q,r$ we can beforehand compute integrals $\int_{0}^{1} \int_{0}^{1} R \left(1, \bar b, \bar c, (x,y) \right)^p x^q y^r \dx \dy$, using the formulae in Section~\ref{sec:recursiveFormulae}, for a dense enough grid of points $(\bar b,\bar c)$ and store them in a lookup table. Intermediate values, not stored in the table, can be computed via interpolation.

\subsubsection{Triangles}
\begin{figure}[t!]
\centering
{\subfigure[Interior source point: \newline ``$T = T_1 + T_2 + T_3$'']{
\begin{tikzpicture}[->,>=stealth',auto,node distance=2cm,
  main node/.style={circle,draw,font=\sffamily\Large\bfseries}]
\node[inner sep=0pt] (IG) at (0,0)
 {\includegraphics[trim = 6.5cm 10.5cm 5.5cm 12cm, clip = true, height=2.5cm]{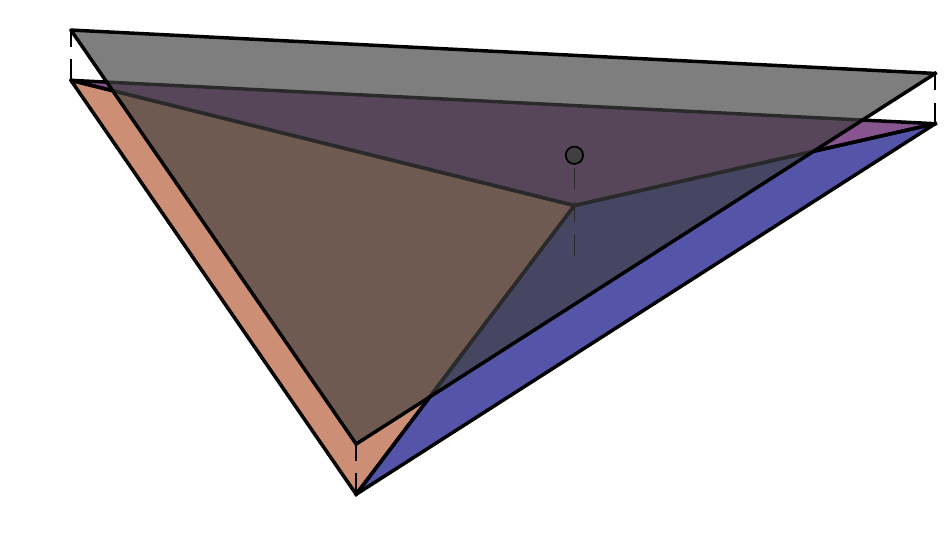}};

\draw[-] (-2.25,.5) -- (-1.45,.5);
\draw[->] [dashed](-1.45,.5) -- (0.25,.5);
\node[text width=0cm] at (-2.75,.5) {$T_2$};

\draw[-] (-2.25,-0) -- (-1.1,-0);
\draw[->][dashed] (-1.1,-0) -- (-0.5,-0);
\node[text width=0cm] at (-2.75,-0) {$T_3$};

\draw[->] (1.5,-.5) -- (0.1,-.5);
\node[text width=0cm] at (1.6,-.5) {$T_1$};
\end{tikzpicture}
}
}
{\subfigure[Edge source point: \newline ``$T = T_1 + T_2$'']{
\begin{tikzpicture}[->,>=stealth',auto,node distance=2cm,
  main node/.style={circle,draw,font=\sffamily\Large\bfseries}]
\node[inner sep=0pt] (IG) at (0,0)
{\includegraphics[trim = 6.5cm 10.5cm 5.5cm 12cm, clip = true, height=2.5cm]{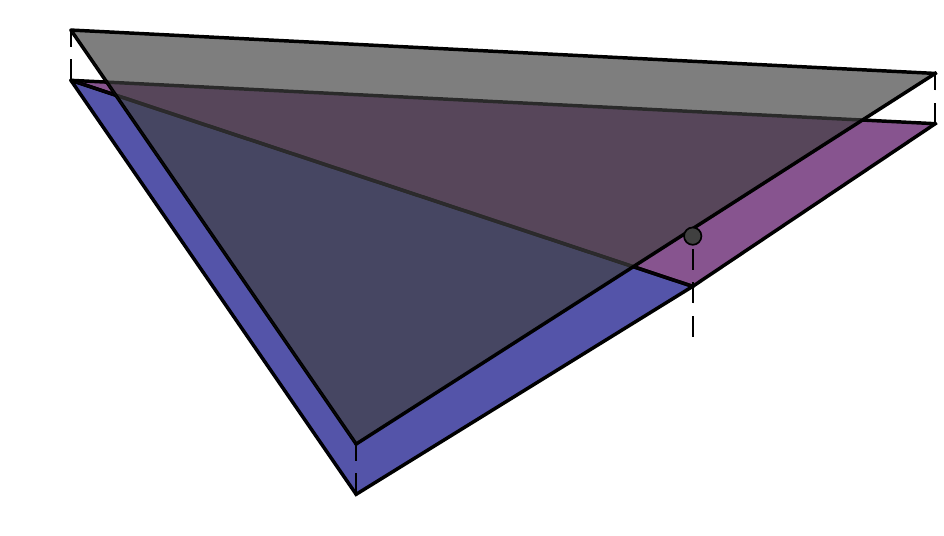}};

\draw[-] (-2.25,.5) -- (-1.45,.5);
\draw[->] [dashed](-1.45,.5) -- (0.25,.5);
\node[text width=0cm] at (-2.75,.5) {$T_2$};

\draw[-] (-2.25,-0) -- (-1.1,-0);
\draw[->][dashed] (-1.1,-0) -- (-0.5,-0);
\node[text width=0cm] at (-2.75,-0) {$T_1$};
\end{tikzpicture}
}
}

\subfigure[Corner source point]{
{\includegraphics[trim = 6.5cm 10.5cm 5.5cm 12cm, clip = true, height=2.5cm]{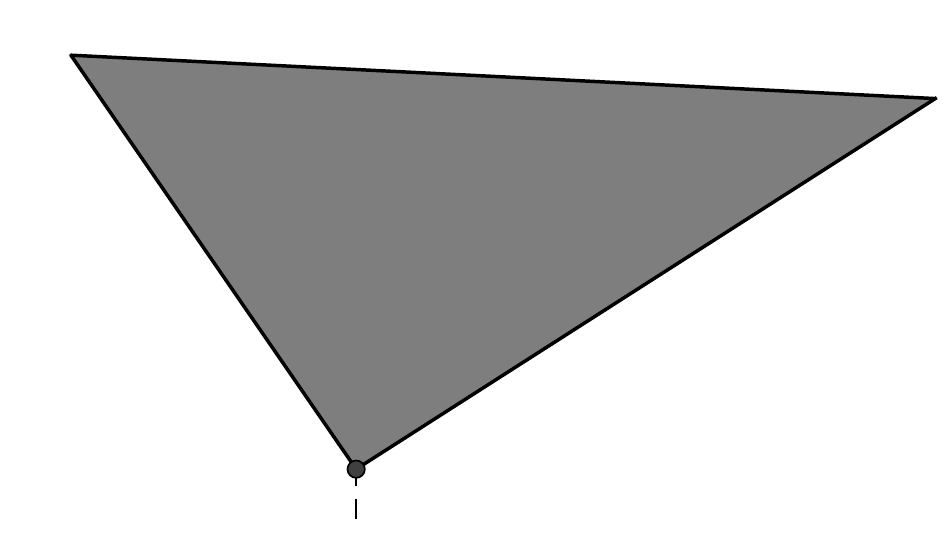}}
}
\subfigure[Exterior source point, close to an edge:  \newline ``$T = T_1 + T_2 - T_3$'']{
\begin{tikzpicture}[->,>=stealth',auto,node distance=2cm,
  main node/.style={circle,draw,font=\sffamily\Large\bfseries}]
\node[inner sep=0pt] (IG) at (0,0)
{\includegraphics[trim = 6.5cm 11.25cm 6cm 12cm, clip = true, height=2.5cm]{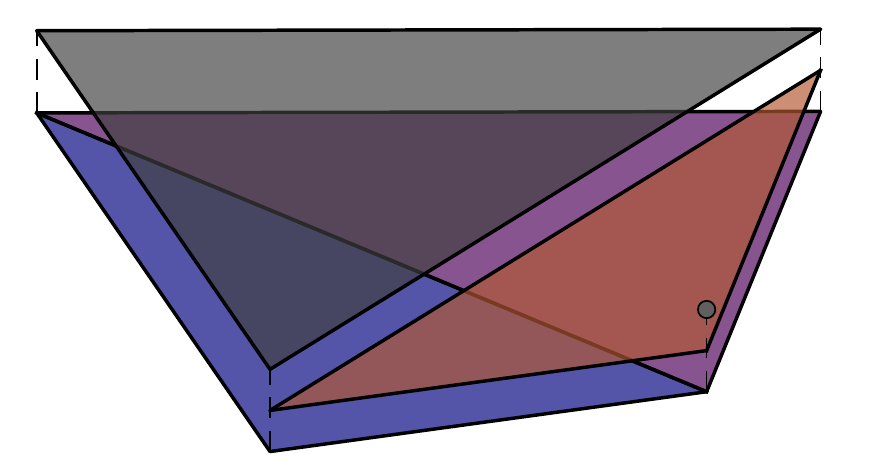}};

\draw[-] (-2.25,.5) -- (-1.75,.5);
\draw[->] [dashed] (-1.75,.5) -- (-0.9,.5);
\node[text width=0cm] at (-2.75,.5) {$T_2$};

\draw[-] (-2.25,.0) -- (-1.4,.0);
\draw[-] [dashed](-1.4,.0) -- (0.05,.0);
\draw[->] (0.05,.0) -- (1,.0);
\node[text width=0cm] at (-2.75,.0) {$T_3$};

\draw[->] (-2.25,-0.6) -- (-1.1,-0.6);
\node[text width=0cm] at (-2.75,-0.6) {$T_1$};

\end{tikzpicture}
}
\subfigure[Exterior source point, close to a corner:  \newline ``$T = T_1 - T_2 - T_3$'']{
\begin{tikzpicture}[->,>=stealth',auto,node distance=2cm,
  main node/.style={circle,draw,font=\sffamily\Large\bfseries}]
\node[inner sep=0pt] (IG) at (0,0)
{\includegraphics[trim = 4.75cm 11.75cm 6.5cm 10.5cm, clip = true, height=2.5cm]{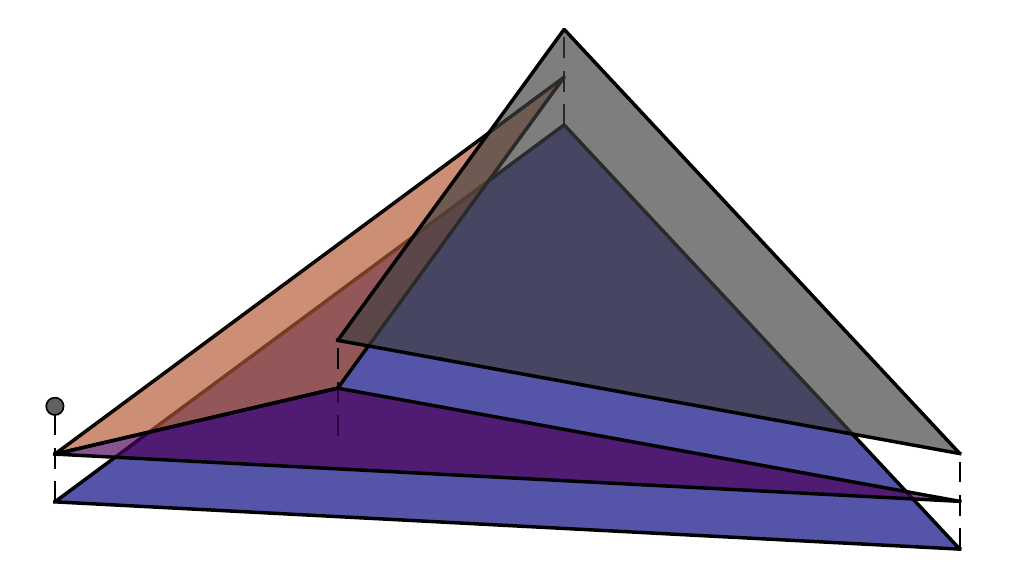}};

\draw[->] (-2.25,.0) -- (-0.8,.0);
\node[text width=0cm] at (-2.75,.0) {$T_3$};

\draw[->] (-2.25,-0.7) -- (-0.55,-0.7);
\node[text width=0cm] at (-2.75,-0.7) {$T_2$};

\draw[-] (1.5,.3) -- (1.1,.3);
\draw[->][dashed] (1.1,.3) -- (0.3,.3);
\node[text width=0cm] at (1.6,.3) {$T_1$};
\end{tikzpicture}
}

\caption{All possible positions of a triangle $T$ (in gray) with respect to the source point $(0,0)$, depicted as a black dot. Sub-triangles are shown in different colours and for better visualization the overlapping ones are depicted at different vertical layers.}
\label{fig:Rxy_cases_T}
\end{figure}
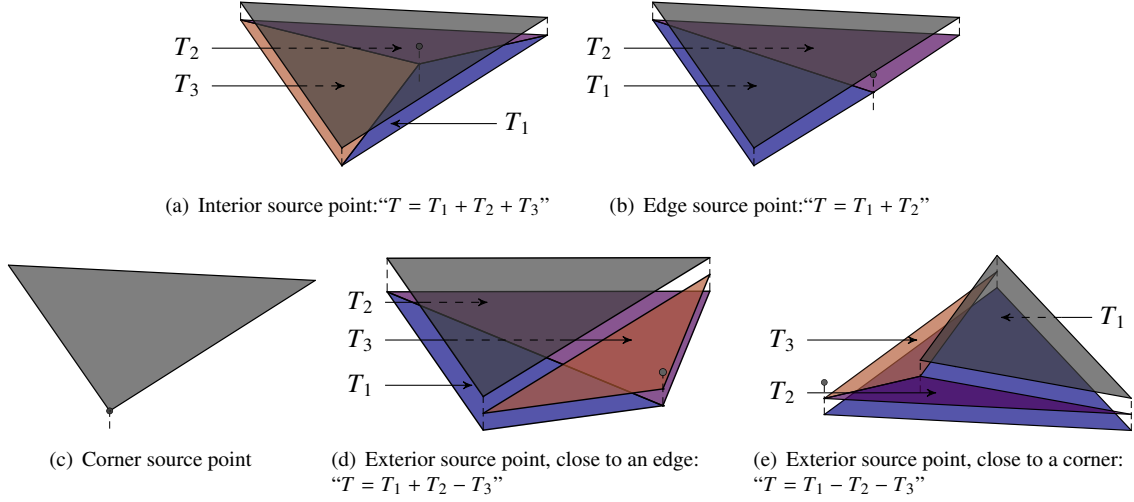

A simplification of the integration on a general triangle $T$ with vertices $(x_0, y_0), (x_1, y_1), (x_2,y_2)$ is more involved compared to the previous aligned rectangular domains. Here we demonstrate one way how to proceed for a general case as a proof of concept -- clearly it is not the optimal nor the most numerical stable approach for every $T$. Since integral
\begin{align*}
I := \int_{T} R^p x^q y^r \, dT
\end{align*}
depends on 9 ``continuous'' variables ($a,b,c,x_0,x_1,x_2,y_0,y_1,y_2$), we can write $I$ again as linear combination of at most 3 integrals, denoted by $I_0$. These integrals contain the same integrand as $I$, but are defined on triangles with $(0,0)$ as one of their corners. The cases are shown in Figure~\ref{fig:Rxy_cases_T}.  In (a) and (b), the initial rectangle $T$ is split into 3 and 2 smaller triangles, respectively. Case (c) is again a special case, when $I=I_0$ and no treatment is needed. In (d),  the source point is outside of $T$ and close to triangle's edge --  triangle $T$ is a union of $T_1$ and $T_2$, where we subtract triangle $T_3$. In (e), $T$ is extended to $T_1$ and then $T_2$ and $T_3$ are subtracted.

In the next step we focus on each of those newly constructed triangles. Let $T_0$ be the integration domain of $I_0$, hence a triangle with vertices $(0, 0), (x_1, y_1), (x_2,y_2)$. Linear transformation that maps a triangle with vertices $(0,0),(1,0), (0,1)$ to $T_0$ is represented by matrix
\begin{align*}
\begin{bmatrix}
x_1 & x_2\\
y_1 & y_2
\end{bmatrix}.
\end{align*}
Therefore, the integral $I_0$ can be written as
\begin{align}\label{eqn:triang_I_0}
I_0 := (x_1 y_2 - x_2 y_1) \int_{0}^1 \int_{0}^{1-y} R(\bar a, \bar b, \bar c, (x,y))^p (x_1 x + x_2 y)^q (y_1 x + y_2 y)^r \dx \dy,
\end{align}
where
\begin{align*}
\bar a &:= a x_1^2 + b x_1 y_1 + c y_1^2,\\
\bar b &:= 2 a x_1 x_2 + b x_2 y_1 + b x_1 y_2 + 2 c y_1 y_2,\\
\bar c &:= a x_2^2 + b x_2 y_2 + c y_2^2.
\end{align*}
The polynomial part in \eqref{eqn:triang_I_0} can be written as
\begin{align}\label{eqn:binom_product}
(x_1 x + x_2 y)^q (y_1 x + y_2 y)^r = \sum_{i=0}^{q+r} \alpha_i x^{q+r-i} y^i
\end{align}
for suitable coefficients $\alpha_i$. By exchanging the order of integration and summation and dividing $R$ by ${a}^{p/2}$, $I_0$ can be written as linear combination of functions
\begin{align*}
\int_{0}^1 \int_{y}^1  R \left(1, \frac{\bar b}{\bar a}, \frac{\bar c}{\bar a}, (x,y)  \right)^p x^{q+r-i} y^i \dx \dy.
\end{align*}

The latter integral depends only on two continuous variables $\bar{ \bar b} := {\bar b}/{\bar a}, \bar {\bar c} := {\bar c}/{\bar a}$ and integer powers $p, q+r-i, i$, and we can evaluate it for a dense enough grid of points $(\bar{ \bar b}, \bar{ \bar c})$ and store the values in a lookup table. At first glance it might seem that every $I_0$ in \eqref{eqn:triang_I_0} needs to {\cc be} computed as a sum of $q+r+1$ integrals due to formula \eqref{eqn:binom_product}. However, observe that the polynomial $P$ in \eqref{eqn:simplifiedInt} already contains different powers of monomials and thus we can group the integrals with the same powers of $x$ and $y$ to reduce the amount of computations.

\subsection{Smoothness of the regularized kernels}\label{sec:Kregularity}
The smoothness of the regularized kernel $\rho$ at $\s$ increases by one if $n$ is increased by one. First we consider the more common singularity subtraction technique, where the analysis is very straightforward.

\begin{thm}
For a singular kernel $\mk_\s$, defined in \eqref{eqn:kernel}, let $\rho =\mk_\s - \mk_{\s,n}$ be the regularized kernel with $m+ n \geq 1$.  Then $\rho$ is a $C^{m+n-1}$ smooth function at $\z = 0$ and $D^\alpha \rho = 0$ for $0 \leq |\alpha| \leq m+n-1$.
\end{thm}

\begin{pf}
Let us write $\mk_\s = \mk_{\s,n} + \tau_{n}$, where $\tau_{n}$ is the corresponding tail in the series expansion of $\mk_\s$. Then $\rho = \tau_{n}$ and from \eqref{eqn:kernel} it is clear that $\zeta(\rho) \geq  m+n$. Apply 
Proposition~\ref{pr:cont} and \ref{pr:smooth} and the proof is complete. 

\qed

\end{pf}

For the regularization by the division step we need to additionally assume $P_{m_2+3}^{[1]}(\z) \neq 0$ for $\z \in D_\varepsilon$. This condition is needed to not introduce new singularities in the regularized kernel $\rho$. In the next proposition we locally assume $P_{m_2+3}^{[1]} > 0$ but an analogous statement can be done for $P_{m_2+3}^{[1]} < 0$.

\begin{cor}\label{cor:dominant}
Let $P_{m_2+3}^{[1]} > 0$ for $\z \in D_\varepsilon$ for some $\varepsilon > 0$. Then  there exists $\varepsilon' > 0$ such that $\mk_{\s, n}(\z) > 0$ for every $\z \in D_{\varepsilon'}$.
\end{cor}

\begin{pf}
Write
\begin{align}\label{eqn:kernelFactorized}
\mk_{\s,n}  = 
\sum_{\ell=1}^n R^{-2 (\ell+m_1)+1} P_{3\ell +m_2}^{[\ell]} = 
R^{-2 (n+m_1)+1} \left ( R^{2n-2} P_{m_2+3}^{[1]} + \sum_{\ell=2}^n R^{2(n-\ell)} P_{3\ell +m_2}^{[\ell]} \right). 
\end{align}
Function $R^{-2 (n+m_1)+1}$ is positive, thus it is enough to analyze the remaining factor. Homogeneous polynomial $R^{2n-2} P_{m_2+3}^{[1]}$ is of degree $2n + m_2 +1$ and positive for $\z \in D_{\varepsilon}$. Each polynomial in the remaining sum is of degree greater than $2n + m_2 +1$, thus it decays to zero faster than $R^{2n-2} P_{m_2+3}^{[1]}$, when $\z$ goes to $\boldsymbol 0$. We can deduce that the kernel $\mk_{\s,n}$ is positive in $D_{\varepsilon'}$ for sufficiently small $\varepsilon'>0$.


\qed

\end{pf}

\begin{rmk}\label{rmk:correctionTerm}
The additional condition $P_{m_2+3}^{[1]} \neq 0$ is clearly an notable drawback of the regularization via division, since not all kernels are necessary of the same sign near the source point $\s$, e.g., kernel $\mh$. For these types of kernel the division approach cannot be directly applied in this simple form. 

The case when $\mk_\s (\z) \neq 0$ for all $\z \in D_\varepsilon$ but $\mk_{\s,n} (\z) = 0$ for some $\z \in D_\varepsilon$ is far less problematic. For example, by adding a simple correction term $\eta \| \z \|_2^2$ to $\mk_{\s,n} (\z)$ for a suitable $\eta$, we can enforce function $\mk_{\s,n} (\z) + \eta \| \z \|_2^2$ to not vanish inside $D_\varepsilon$. Namely, $\eta>0$ if $\mk_\s > 0$, and $\eta<0$ if $\mk_\s < 0$.
\end{rmk}

Despite the mentioned drawbacks, probably the main advantage of the division technique is the improved smoothness of $\rho$ when compared to the subtraction splitting for the same amount of regularization terms. Furthermore, the gap between the approaches is even bigger for kernels $\mk_\s$  with stronger singularity. Indeed, for the same number of terms $n$ the smoothness of $\rho$ is the same, regardless of the type of singularity, when the regularization via division step is used.

The following lemma is needed before proving the theorem for the smoothness of $\rho$.

\begin{lma}\label{lem:cont}
Let $\varepsilon >0$ and let
\begin{align*}
p_m(\z) := \sum_{|\alpha|=m} c_\alpha \z^{\alpha}
\end{align*}
be a bivariate homogeneous polynomial of degree $m$, such that $p_m(\z)>0$ for $\| \z \|_2 < \varepsilon$. Then
$
1/ ({R^\ell p_{m}})
$
is a continuous function at $\z = \boldsymbol 0$ if $\zeta ({R^\ell p_{m}})= \ell+m \leq -1$.
\end{lma}

\begin{pf}
The proof is similar to the proof of Proposition~\ref{pr:cont} but considering a negative parameter $\zeta$ instead of a positive one.

Fix a sufficiently small $\delta > 0$. Let us show that there exists a small enough $\varepsilon$, $0 < \varepsilon<1$, such that for $\| \z \|_2  \leq \varepsilon$ it follows $1/({R(\z)^{\ell} p_{m} (\z)}) \leq \delta$.

Let us again consider only the case $\z \in T_1$ (see \eqref{eqn:subdomains}), since the proof for the other subdomains follows similarly. We can write
$z_2= \gamma z_1$ for $0 \leq \gamma \leq 1$. 
From Lemma~\ref{lem:lemma1} we can infer that there exists a positive constant $C_{abc\ell,1}$ such that
$ C_{abc\ell,1} \leq ({a  + b \gamma + c \gamma^2})^{\ell/2}$ for all $\gamma$. 
From the assumptions we know that $p_m(\z) =  z_1^{m} \sum_{|\alpha|=m} c_\alpha \gamma^{\alpha_2}  > 0$ for $z_1>0$, thus there exists a positive constant $C$ such that $C \leq \sum_{|\alpha|=m} c_\alpha \gamma^{\alpha_2}$. By combining all the derived inequalities we can show that
\begin{align*}
R(\z)^{\ell} p_m (\z) = {(a  + b \gamma + c \gamma^2)}^{\ell/2} z_1^{\ell+m} \sum_{|\alpha|=m} c_\alpha \gamma^{\alpha_2} \geq C_{abc\ell,1}\, C\,  z_1^{\ell+m}
\end{align*}
and
\begin{align*}
\frac 1 {R(\z)^{\ell} p_m(\z)} \leq C_{abc\ell,1}^{-1}\, C^{-1}\,  z_1 \leq C_{abc\ell,1}^{-1}\, C^{-1}\,  \varepsilon.
\end{align*}
By setting $ \varepsilon \leq C_{abc\ell,1}\, C\, \delta$ the proof is complete.

\qed

\end{pf}

\begin{thm}

  Let $\mk_{\s,n}$ be an approximation of $\mk_\s$, defined in \eqref{eqn:kernel}, with $P_{m_2+3}^{[1]}(\z) \neq 0$ for $\z \in D_\varepsilon$. Let $\rho =\mk_\s / \mk_{\s,n}$ be the regularized kernel.  Then $\rho$ is a $C^{n-1}$ smooth function at $\z = \boldsymbol 0$ and $D^\alpha \rho(\boldsymbol 0) = 0$ for $1 \leq |\alpha| \leq n-1$.
\end{thm}

\begin{pf}
%
%
Let us write again $\mk_\s = \mk_{\s,n} + \tau_{n}$, where $\tau_{n}$ is the corresponding tail in the series expansion of $\mk_\s$. 

First, let us prove it for $n=1$, i.e., function $\rho = \mk_\s / \mk_{\s,1}$ is $C^0$ continuous. The function $\rho$ can be written as
\begin{align*}
\rho = \frac{\mk_\s}{\mk_{\s,1}} = 1+ \frac{\tau_1}{\mk_{\s,1}}
%
 = 1 + \frac{R^{-m-1} \tau_1}{R^{-m-1} \mk_{\s,1}}
 = 1 + \frac{R^{-m-1} \tau_1}{R^{-m -2m_1-2} P_{m_2+3}^{[1]}}.
\end{align*}
Since $\zeta(R^{-m-1} \tau_1) \geq (-m-1) + m+1 = 0$ function $R^{-m-1} \tau_1(\z)$ is bounded in $D_\varepsilon$ (see Proposition~\ref{pr:operations}) and 
we can find a positive constant $C_1$, independent of $\z$, such that $|R(\z)^{-m-1}\tau_1(\z)| \leq C_1$. For the denominator it holds $\zeta(R^{-m-1} \mk_{\s,1})=-1$, thus it is not bounded for $\| \z \|_2 \leq \varepsilon$. 
By defining $\rho(\boldsymbol 0) = 1$ we can see that
\begin{align*}
|\rho(\z) - \rho(\boldsymbol 0)| = |\rho(\z) - 1| = \left |\frac{R(\z)^{-m-1} \tau_1 (\z)}{R(\z)^{-m-1} \mk_{\s,1} (\z)}\right |  \leq  \frac{C_1}{|R(\z)^{-m-1} \mk_{\s,1} (\z)|}.
\end{align*}
From Lemma~\ref{lem:cont} it follows that the function $1 / | R^{-m-1} \mk_{\s,1}|$ is continuous at $\z = \boldsymbol 0$, hence so is $\rho$.

Let us prove the remaining part by induction. Assume $D^\alpha \rho$ is continuous and $D^\alpha \rho(\boldsymbol 0) = 0$ for $0 \leq |\alpha| \leq n-1$. Then we can define the following derivatives via recursion,
\begin{align}\label{eqn:pij}
\nonumber \rho = 1+ \frac{\tau_{n}}{\mk_{\s,{n}}} =: 1+ \frac{p_{0,0}}{q_{0,0}},\\
D^{i+1}_1 D^j_2 \rho = D^1_1 \frac{ p_{i,j}}{ q_{i,j}} = \frac{D^1_1 p_{i,j}\, q_{i,j} - p_{i,j}\, D^1_1 q_{i,j}}{(q_{i,j})^2}=:\frac{p_{i+1,j}}{q_{i+1,j}},\\
\nonumber D^{i}_1 D^{j+1}_2 \rho = D^1_2 \frac{ p_{i,j}}{ q_{i,j}} = \frac{D^1_2 p_{i,j}\, q_{i,j} - p_{i,j}\, D^2_1 q_{i,j}}{(q_{i,j})^2}=:\frac{p_{i,j+1}}{q_{i,j+1}},
\end{align}
for $i+j \leq n-1$.
Observe that $\zeta(p_{0,0}) = \zeta(\tau_{n}) \geq m+n$, while $\zeta(q_{0,0}) = \zeta(\mk_{\s,{n}}) = m$. From expressions \eqref{eqn:pij} and Proposition~\ref{pr:operations} we can infer that $\zeta(q_{i,j}) = 2^{i+j} m$ and $\zeta(p_{i,j}) \geq 2^{i+j}m + n-i-j$.

For $|\alpha| = i+j = n-1$ we can write
\begin{align}\label{eqn:contnDer}
\left| \frac{p_{i,j}(\z)}{q_{i,j}(\z)} - \boldsymbol 0 \right| = \frac{|R(\z)^{-2^{n-1} m - 1}\, p_{i,j} (\z)|}{R(\z)^{-2^{n-1} m - 1}\, |q_{i,j}(\z)|} \leq \frac{C_2}{R(\z)^{-2^{n-1} m - 1}\, |q_{i,j}(\z)|},
\end{align}
where $C_2$ is a positive constant, since $\zeta(R^{-2^{n-1}m - 1}\, p_{i,j}) \geq 0$. 
From the proof of Corollary~\ref{cor:dominant} it is easy to see that there exists a small enough $\varepsilon'>0$ such that
\begin{align*}
|q_{i,j}(\z)| = |\mk_{\s,{n}}^{2^{n-1}}(\z)| \geq \left| \left(\frac{1}{2} R(\z)^{-2 m_1-1} P_{m_2+3}^{[1]}(\z) \right)^{2^{n-1}} \right|
\end{align*}
for all $\z \in D_{\varepsilon'}$. Then for  $\z \in D_{\varepsilon'}$ we can further estimate the right-hand side of \eqref{eqn:contnDer},
\begin{align}\label{eqn:contnDer2}
 \frac{C_2}{R(\z)^{-2^{n-1} m - 1}\, |q_{i,j}(\z)|} \leq  \frac{2^{2^{n-1}} C_2}{R(\z)^{- 2^{n-1} (m_2+3)-1}\, |P_{m_2+3}^{[1]}(\z)^{2^{n-1}}|}.
\end{align}
The function on the right-hand side of \eqref{eqn:contnDer2} is continuous at $\z = \boldsymbol 0$ by Lemma~\ref{lem:cont}. Hence, $D^\alpha \rho = {p_{i,j}}/{q_{i,j}}$ is continuous for $|\alpha| = i+j = n-1$ and $D^\alpha \rho(\boldsymbol 0) = 0$.


\qed

\end{pf}

\section{Numerical tests}

In all experiments we test the accuracy of numerical integration for specific singular integrals with respect to the number of terms $n$ in the singular kernel series expansion; we vary $n=0,1,2,3$, whereas in the case $n=0$ we do not apply any regularization of the kernel. For regularization we considered the singularity subtraction and division (see Section~\ref{sec:num_int}).

\begin{figure}[ht!]
\centering
{\includegraphics[trim = 4.5cm 8.5cm 6cm 9.25cm, clip = true, height=5cm]{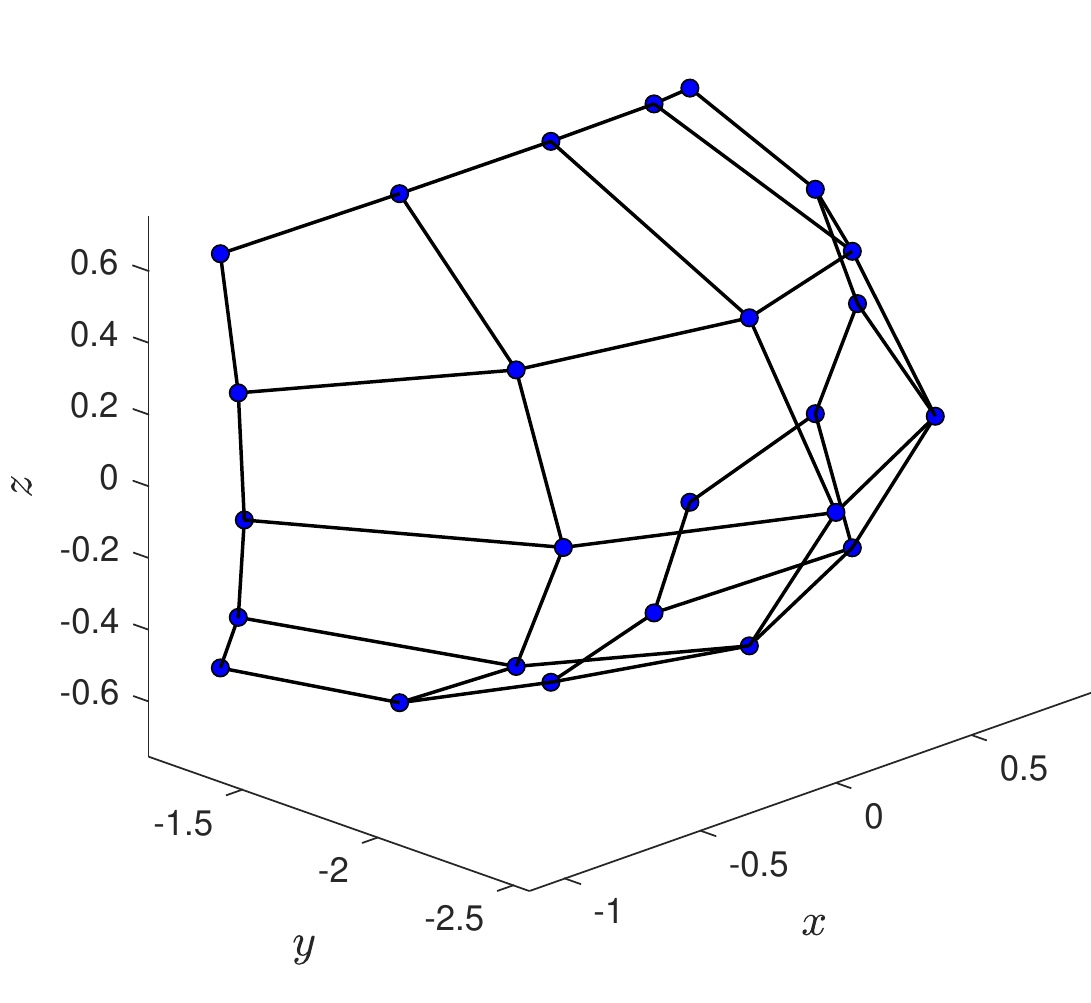}}
\hspace{.5cm}
{\includegraphics[trim = 5.5cm 8.5cm 6cm 9.25cm, clip = true, height=5cm]{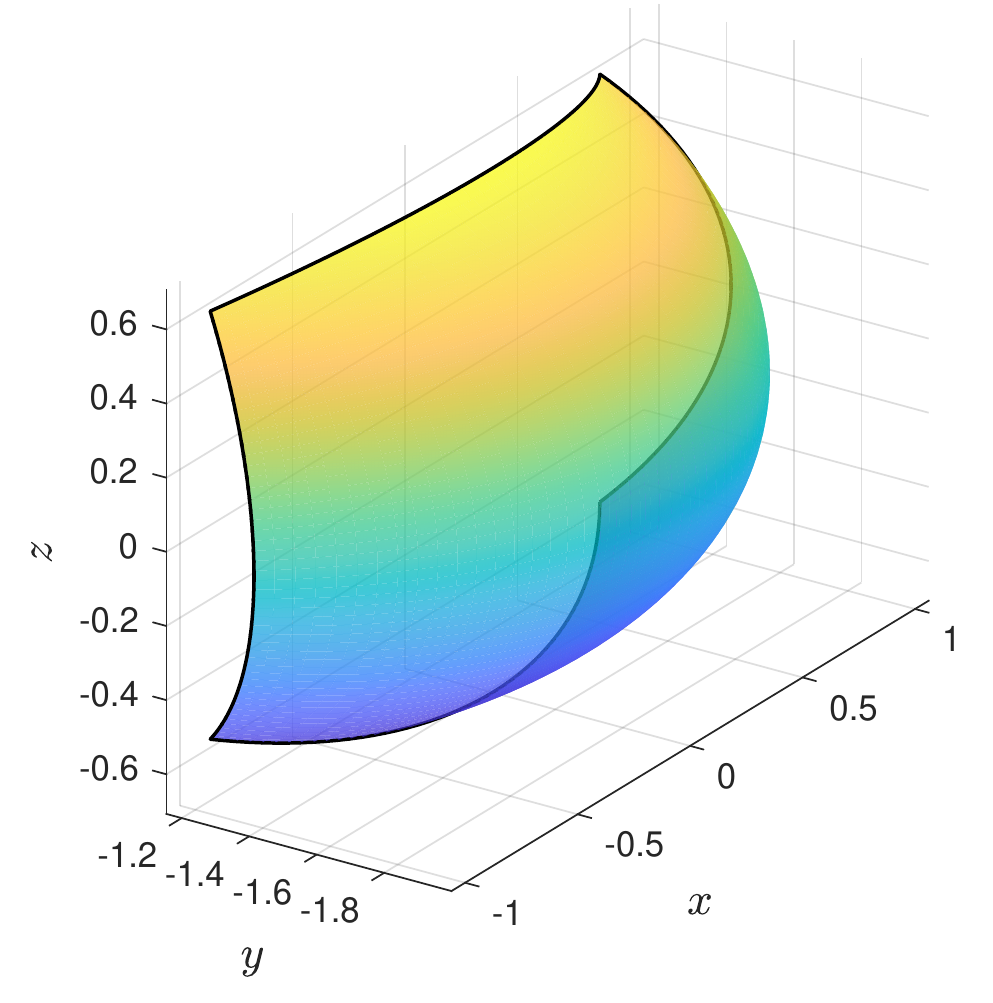}}

\caption{Control points and geometry of the sphereoid.}
\label{fig:spheroid_geom}
\end{figure}

For geometry let us consider a section of a spheroid, parameterized by quartic tensor product NURBS $\f: [0,1]^2 \to \RR^3$ (obtained by stretching a spherical section \cite{Cobb1994TilingTS}); see Fig.~\ref{fig:spheroid_geom}.
Its control points are $\boldsymbol c_{ij} = (C^x_{ij}, C^y_{ij}, C^z_{ij})/C^w_{ij}$ by defining the following matrices
\begin{align*}
C^x &= 1.5 
\begin{bmatrix}
-c_1 &-c_3 &-c_4 &-c_3 &-c_1\\
-c_2 &-c_5 &-c_6 &-c_5 &-c_2\\
0 &0 &0 &0 &0 \\
c_2& c_5 &c_6 &c_5 &c_2\\
c_1 &c_3 &c_4 &c_3& c_1
\end{bmatrix}, \;\;
C^y = -2
\begin{bmatrix}
c_1 &c_3 &c_4 &c_3 &c_1\\
c_3 &c_7 &c_8 &c_7 &c_3\\
c_4 &c_8 &c_9 &c_8 &c_4\\
c_3 &c_7 &c_8 &c_7 &c_3\\
c_1 &c_3 &c_4 &c_3 &c_1
\end{bmatrix}, \;\;
C^z = 
\begin{bmatrix}
-c_1 &-c_2 &0 & c_2&c_1\\
-c_3 &-c_5 &0 &c_5& c_3\\
-c_4 &-c_6 &0 &c_6 &c_4\\
-c_3 &-c_5 &0 &c_5 &c_3\\
-c_1 &-c_2& 0 &c_2 &c_1
\end{bmatrix}, \\
C_w &= 
\begin{bmatrix}
w_1 &w_2 &c_9 &w_2 &w_1\\
w_2 &c_7 &w_3 &c_7 &w_2\\
c_9 &w_3 &w_4 &w_3 &c_9\\
w_2 &c_7 &w_3 &c_7 &w_2\\
w_1 &w_2 &c_9 &w_2 &w_1
\end{bmatrix},
\end{align*}
where
\begin{align*}
c_1 &= 4(\sqrt{3}-1),
& c_2 &= \sqrt{2},
& c_3 &= \sqrt{2}(4-\sqrt{3}),
& c_4 &= 4(2\sqrt{3}-1)/3,
& c_5 &= (3\sqrt{3}-2)/2,\\
c_6 &= \sqrt{2}(7-2\sqrt{3})/3,
& c_7 &= (\sqrt{3}+6)/2,
& c_8 &= 5\sqrt{6}/3,
& c_9 &= 4(5-\sqrt{3})/3,\\
w_1  &= 4(3 - \sqrt{3}),
& w_2  &= \sqrt{2}(3\sqrt{3}-2),
& w_3  &= \sqrt{2}(\sqrt{3}+6)/3,
& w_4  &= 4(5\sqrt{3} - 1)/9.
\end{align*}

\subsection{N-refinement}

In the first example we consider a simplified version of a governing integral in BIE \eqref{eqn:BIEdirect} to test the effect of the regularization technique in the most direct manner. We test the accuracy of the numerical schemes in Section~\ref{sec:num_int} to compute
\begin{align}\label{eqn:test1}
\frac{1}{4 \pi}\int_{[0,1]^2} \mk(\s, \t) \dt,
\end{align}
where $\mk \in \{\mg, \bar \mh\}$ for 16 different source points $\s=(s_1,s_2)$ with $s_1, s_2  \in \{ 0.6, 0.7, 0.8,  0.9\}$ and for different number of quadrature points $N$. The integral can be considered as a simplified version of \eqref{eqn:initInt}, where $B_j \equiv 1$ and $v \equiv 1$.

\begin{figure}[ht!]
\centering
\subfigure[Subtraction, $n=1$]{
{\includegraphics[trim = 0cm .75cm 0.25cm 1cm, clip = true, height=4cm]{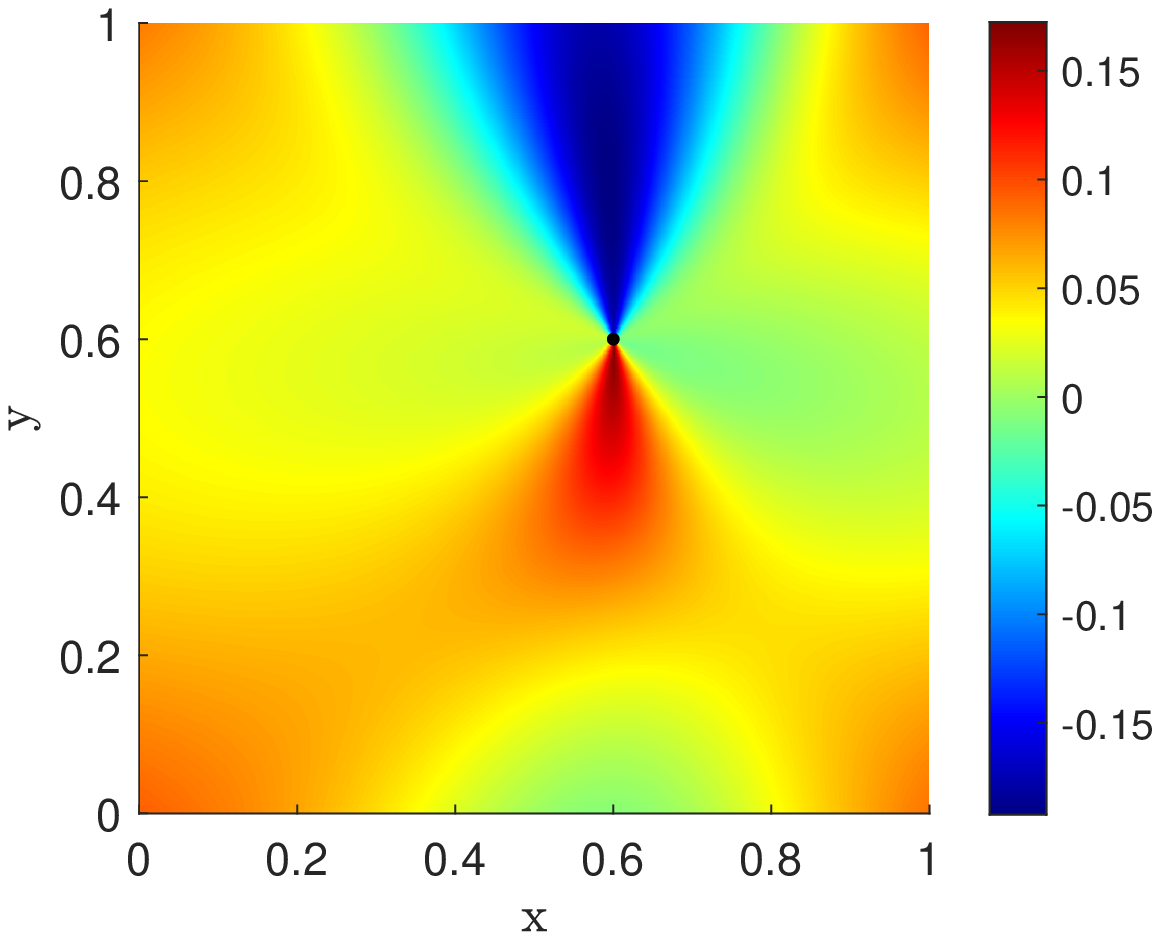}}
}
\subfigure[Subtraction, $n=2$]{
{\includegraphics[trim = 0cm .75cm 0.25cm 1cm, clip = true, height=4cm]{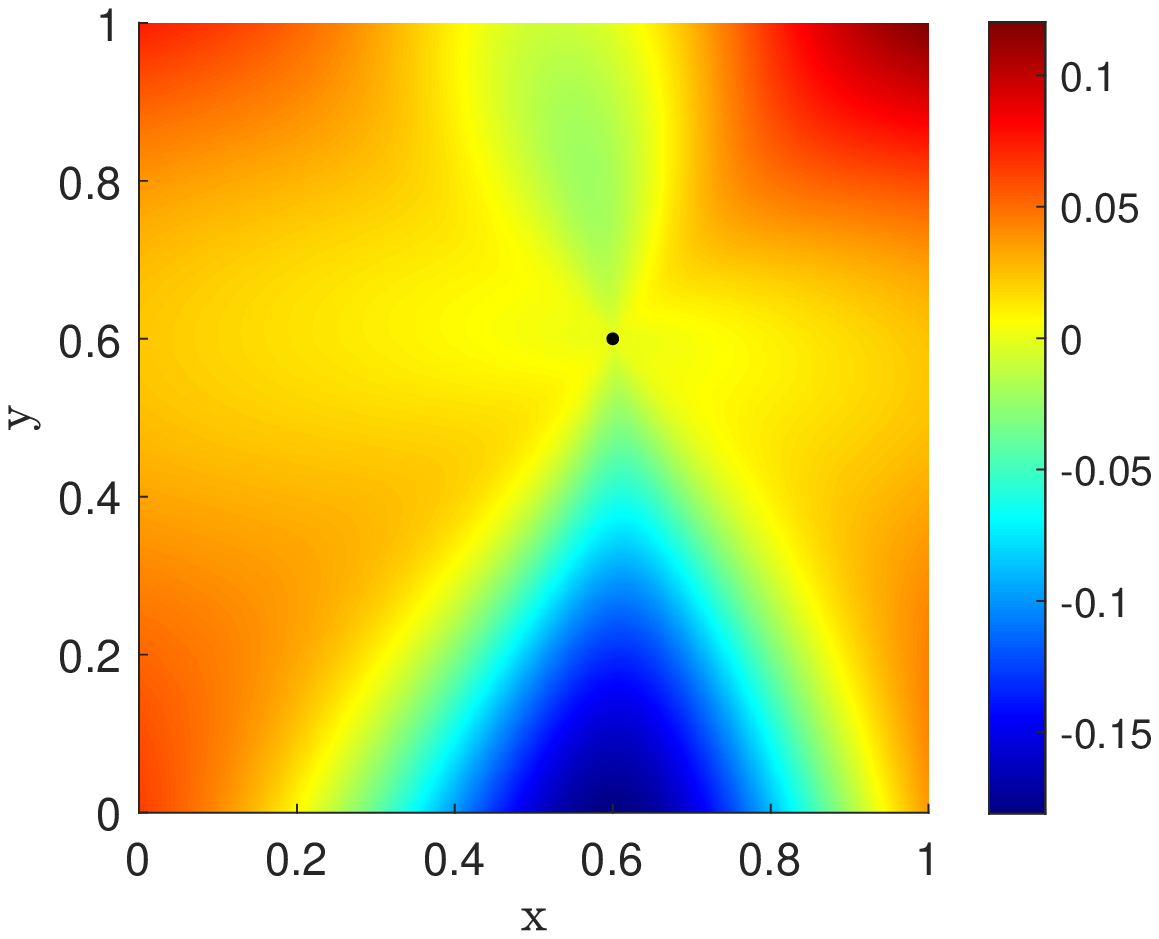}}
}
\subfigure[Subtraction, $n=3$]{
{\includegraphics[trim = 0cm .75cm 0.25cm 1cm, clip = true, height=4cm]{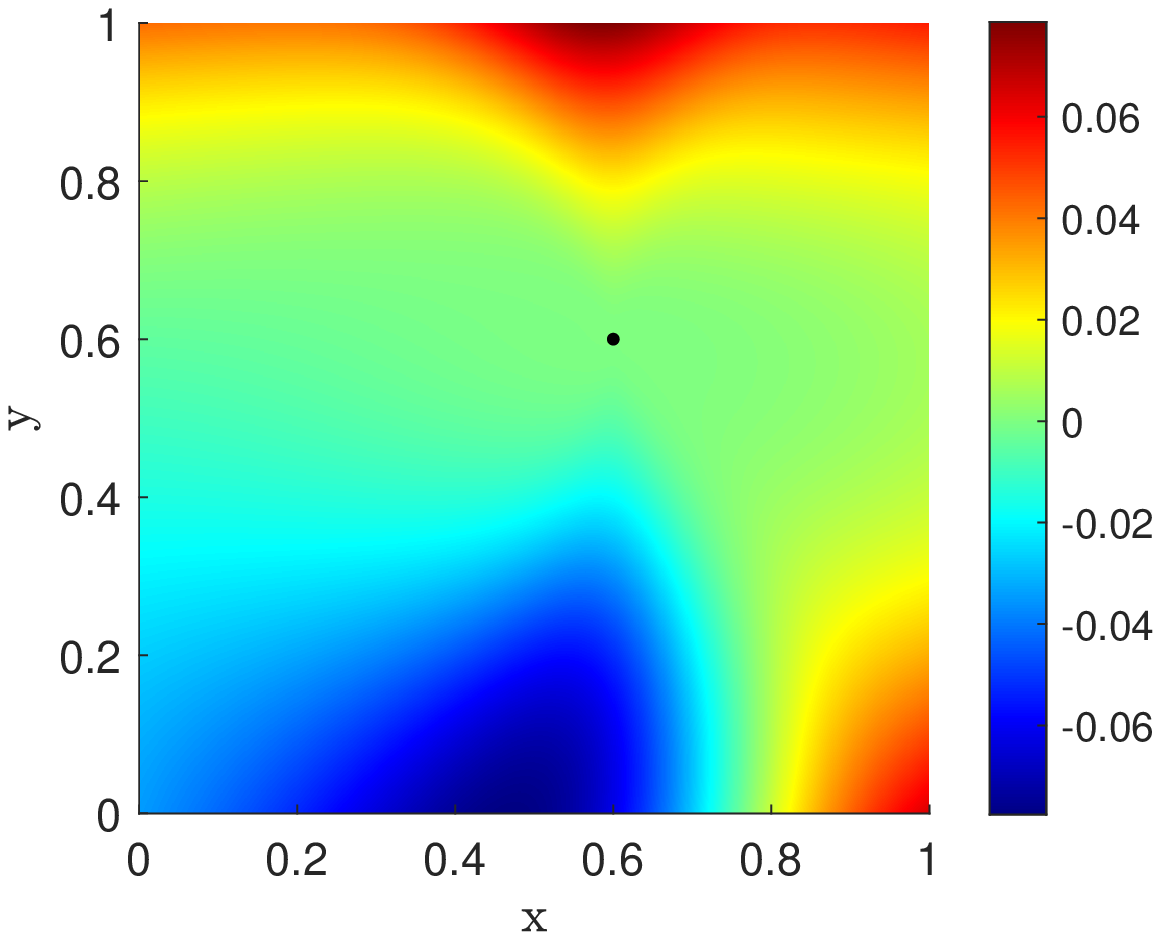}}
}

\subfigure[Division, $n=1$]{
{\includegraphics[trim = 0cm .75cm 0.25cm 1cm, clip = true, height=4cm]{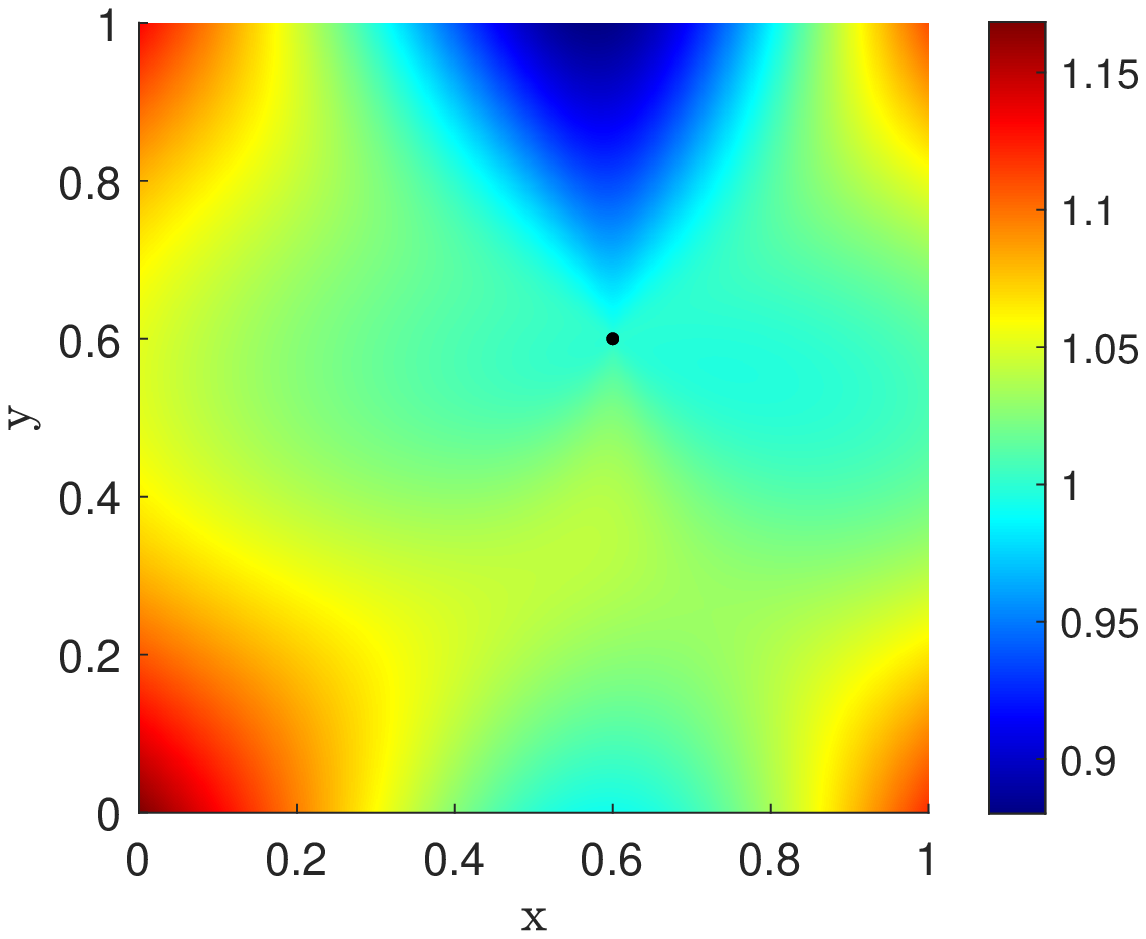}}
}
\subfigure[Division, $n=2$]{
{\includegraphics[trim = 0cm .75cm 0.25cm 1cm, clip = true, height=4cm]{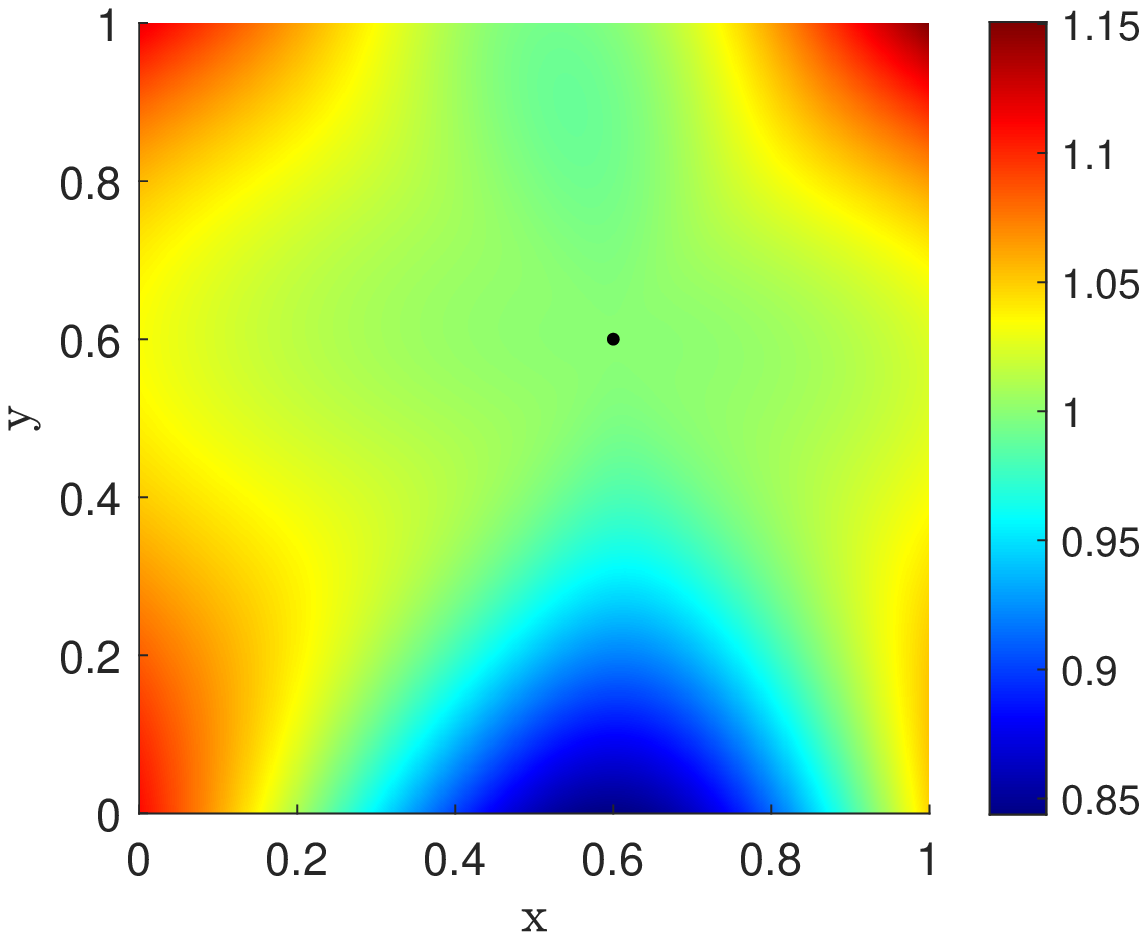}}
}
\subfigure[Division, $n=3$]{
{\includegraphics[trim = 0cm .75cm 0.25cm 1cm, clip = true, height=4cm]{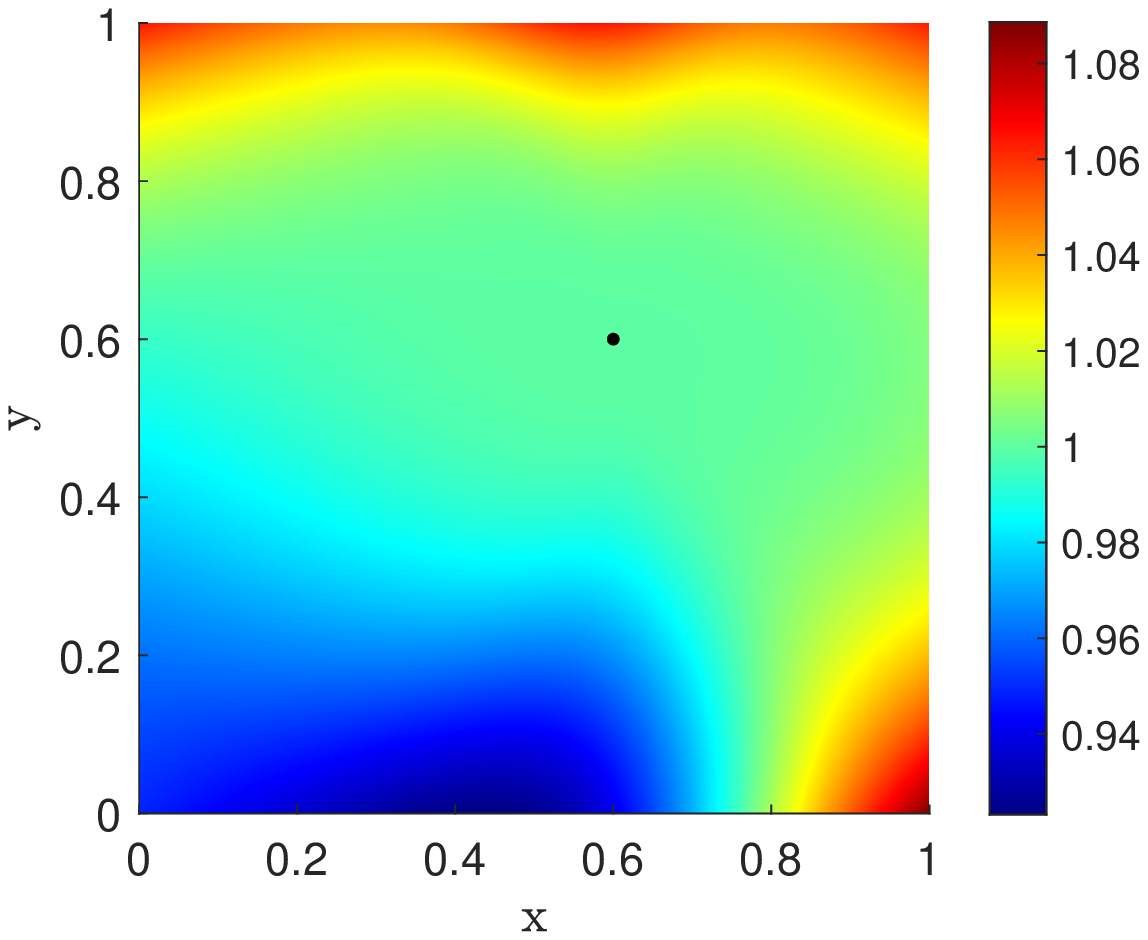}}
}

\caption{Regularized kernel $\rho$ for $\mg_\s$ with $\s=(0.6,0.6)$ marked as a black dot.}
\label{fig:regularizedKernelG}
\end{figure}

\begin{figure}[ht!]
\centering
\subfigure[Subtraction, $n=1$]{
{\includegraphics[trim = 0cm .75cm 0.25cm 1cm, clip = true, height=4cm]{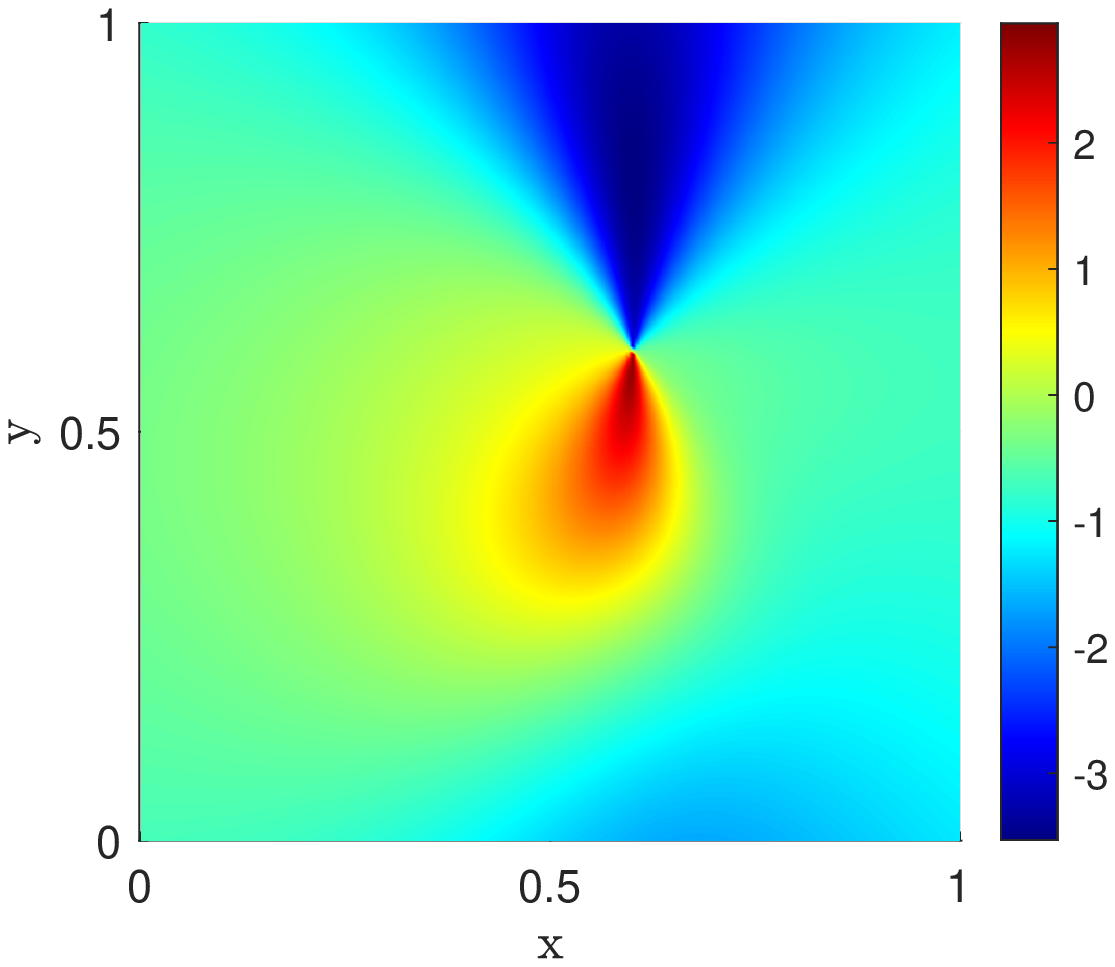}}
}
\subfigure[Subtraction, $n=2$]{
{\includegraphics[trim = 0cm .75cm 0.25cm 1cm, clip = true, height=4cm]{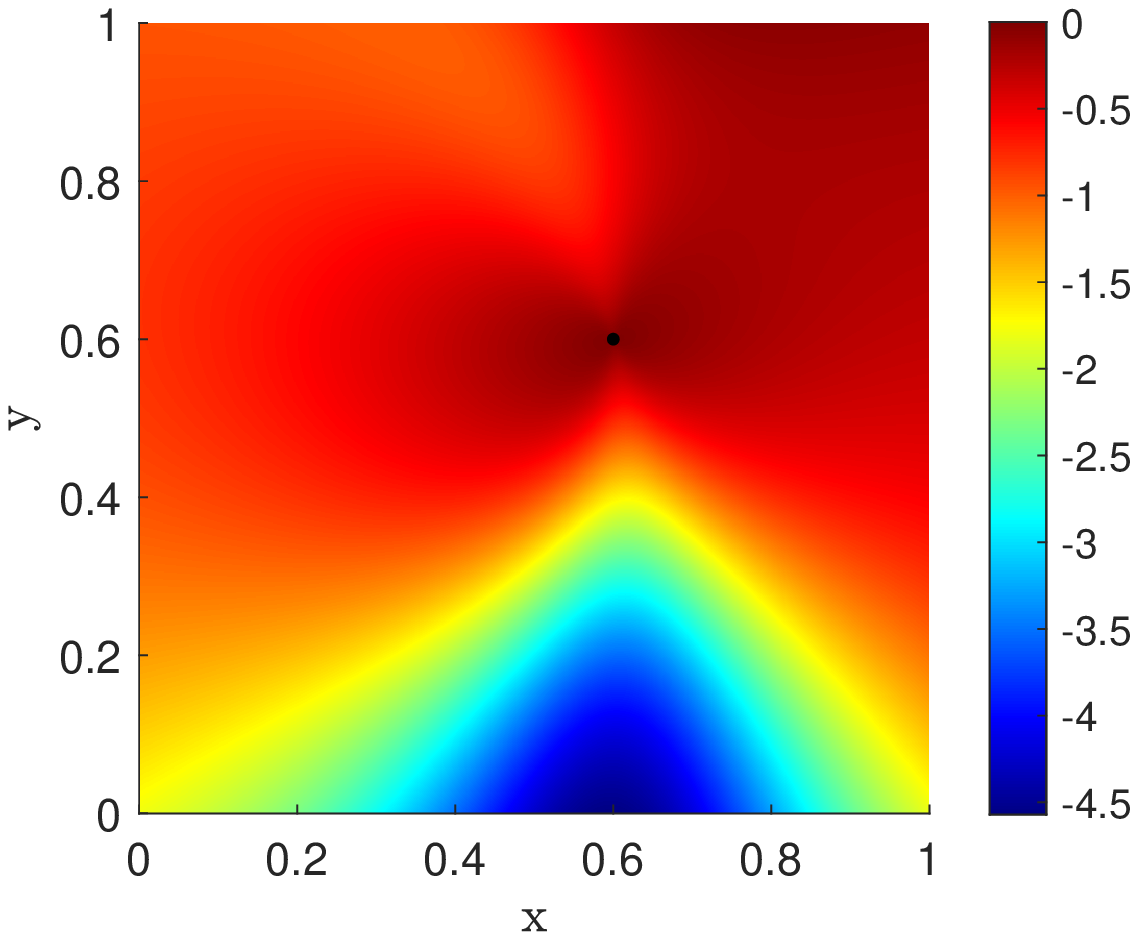}}
}
\subfigure[Subtraction, $n=3$]{
{\includegraphics[trim = 0cm .75cm 0.25cm 1cm, clip = true, height=4cm]{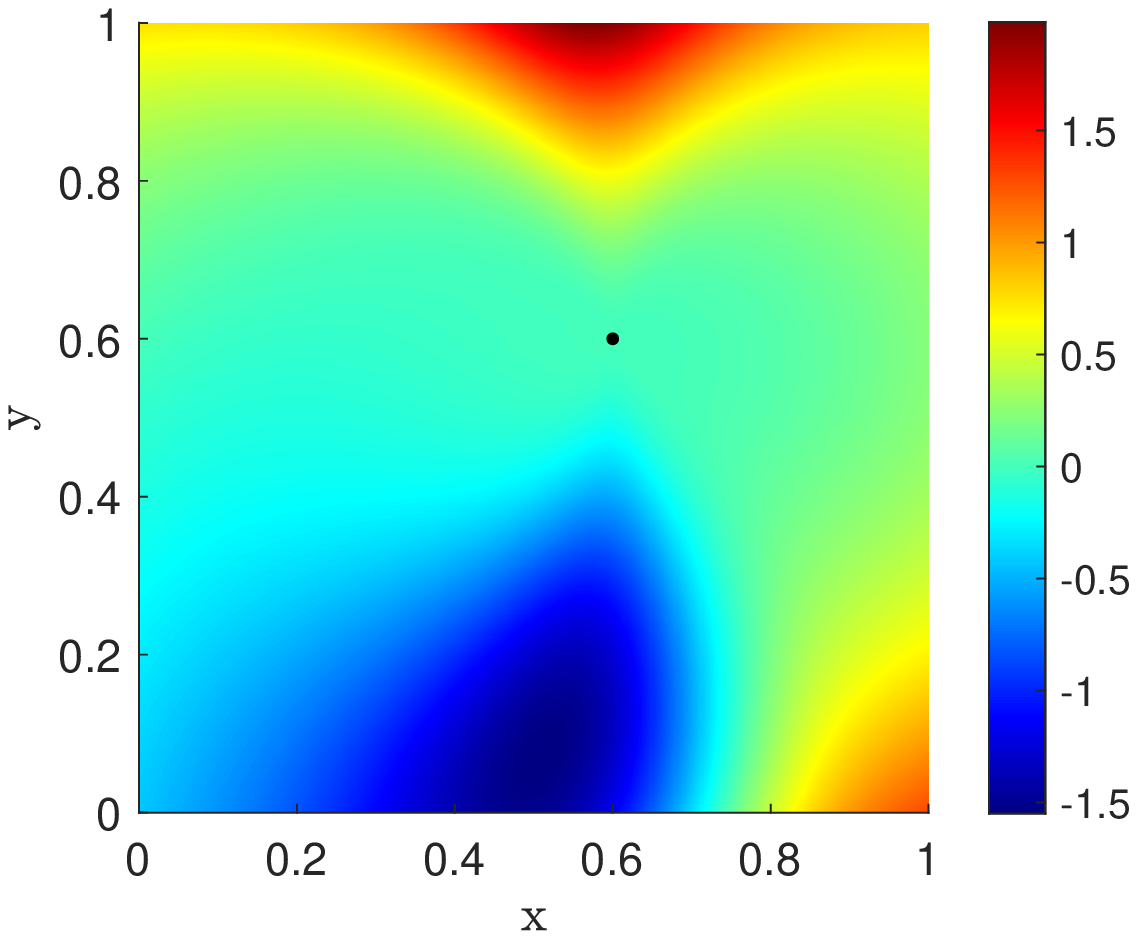}}
}

\subfigure[Division, $n=1$]{
{\includegraphics[trim = 0cm .75cm 0.25cm 1cm, clip = true, height=4cm]{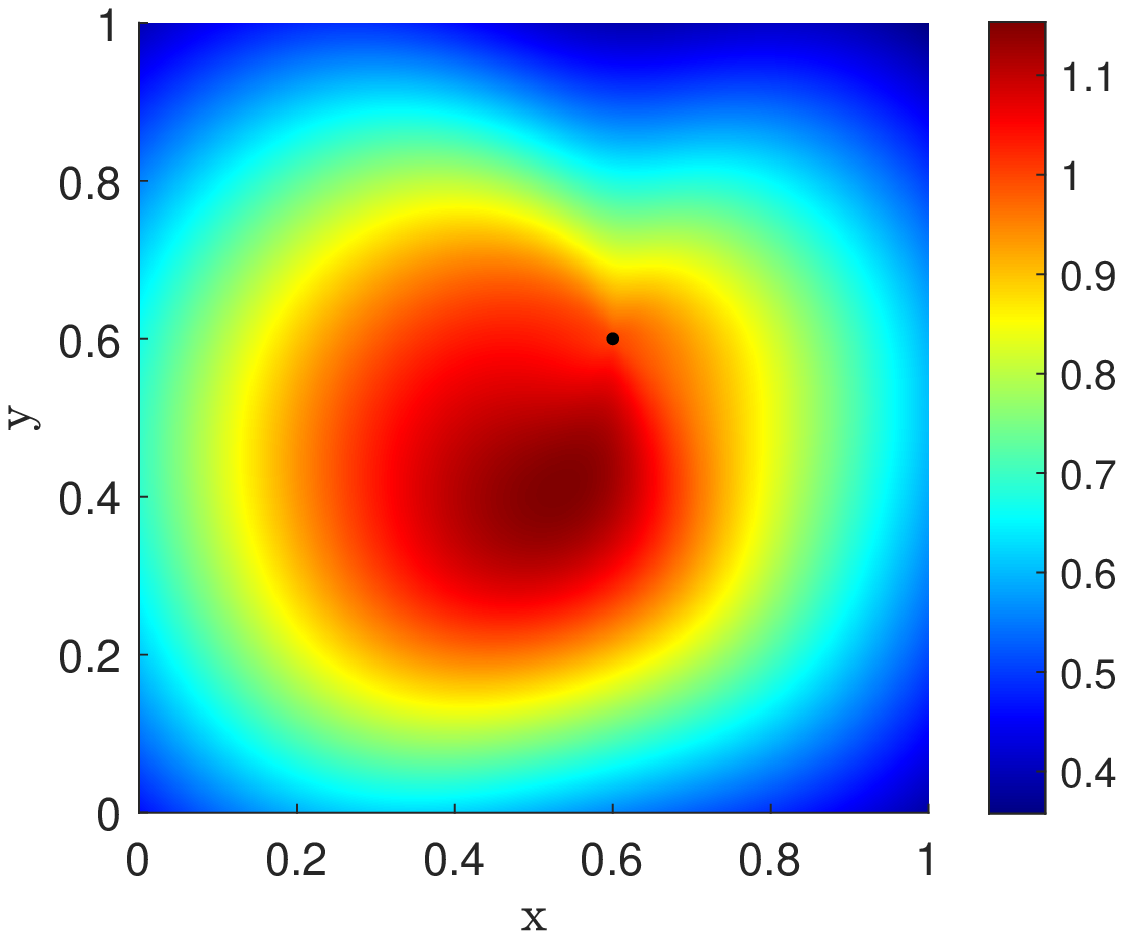}}
}
\subfigure[Division, $n=2$]{
{\includegraphics[trim = 0cm .75cm 0.25cm 1cm, clip = true, height=4cm]{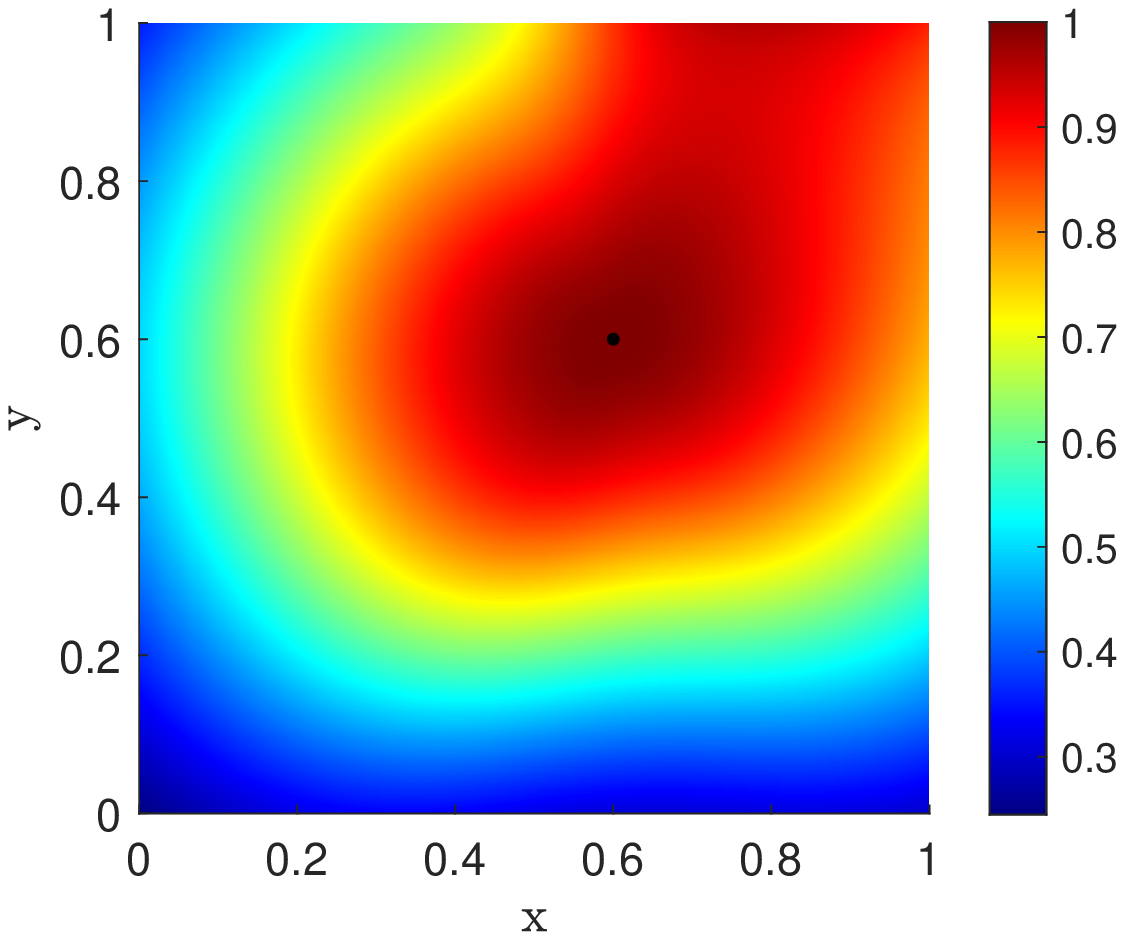}}
}
\subfigure[Division, $n=3$]{
{\includegraphics[trim = 0cm .75cm 0.25cm 1cm, clip = true, height=4cm]{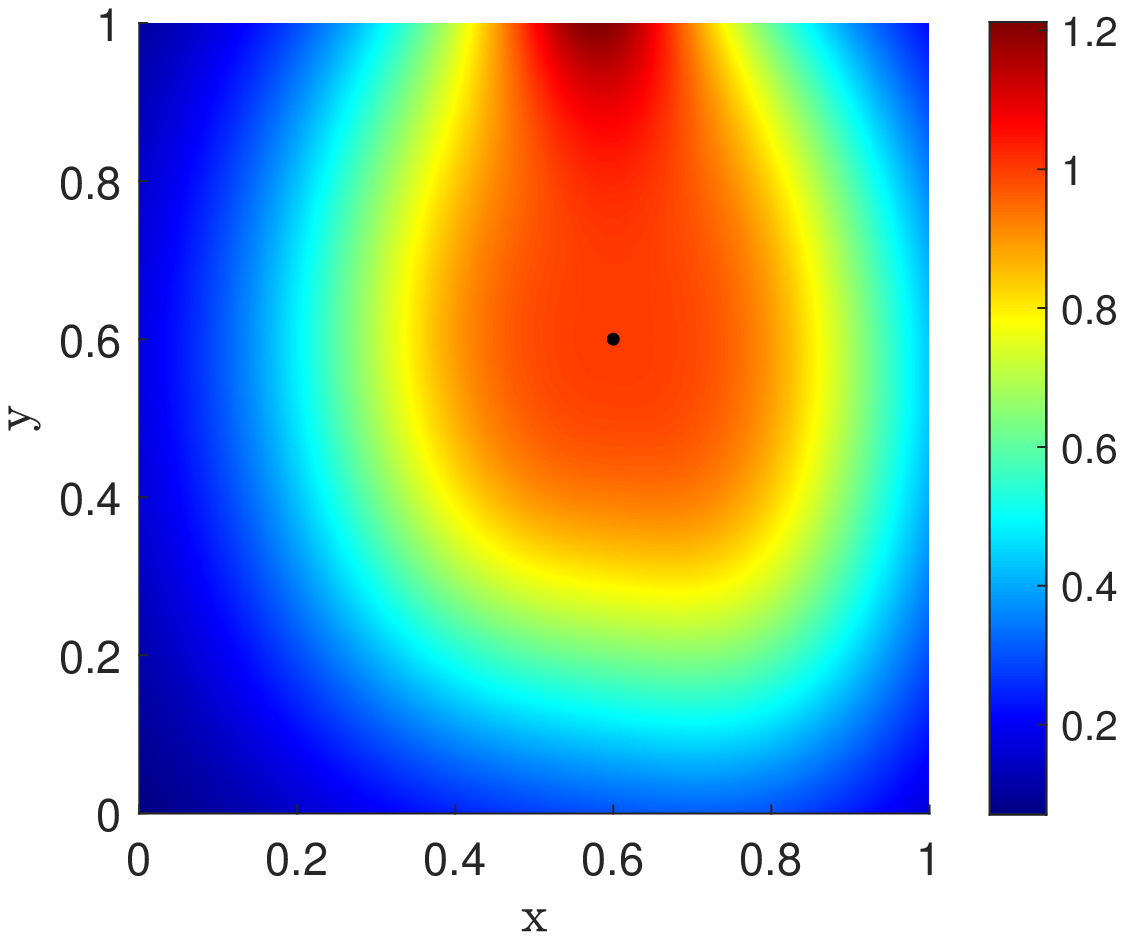}}
}

\caption{Regularized kernel $\rho$ for $\bar \mh_\s$ with $\s=(0.6,0.6)$ marked as a black dot.}
\label{fig:regularizedKernelH}
\end{figure}

In Fig.~\ref{fig:regularizedKernelG} and Fig.~\ref{fig:regularizedKernelH} we show the plots of the regularized kernels $\rho$ for the subtraction and the division regularization step via kernel $\mk_{\s,n}$ for $\cc \s=(0.6,0.6)$. As analyzed in Section~\ref{sec:Kregularity}, for the singularity subtraction the function $\rho$ is bounded, $C^0$ and $C^1$ for $n=1,2,3$, respectively. Instead, for the division step the regularity of $\rho$ increases by one. For the latter regularization for $n=3$ the function $\rho$ exhibits singularities near the edges of the domain $[0,1]^2$; the issue is resolved by adding the correction term with $\eta=10$ (see Remark~\ref{rmk:correctionTerm}).

\begin{figure}[ht!]
\centering
\subfigure[kernel $\mg_\s$]{
\includegraphics[trim = 0cm .0cm 0cm 0.0cm, clip = true, height=6cm]{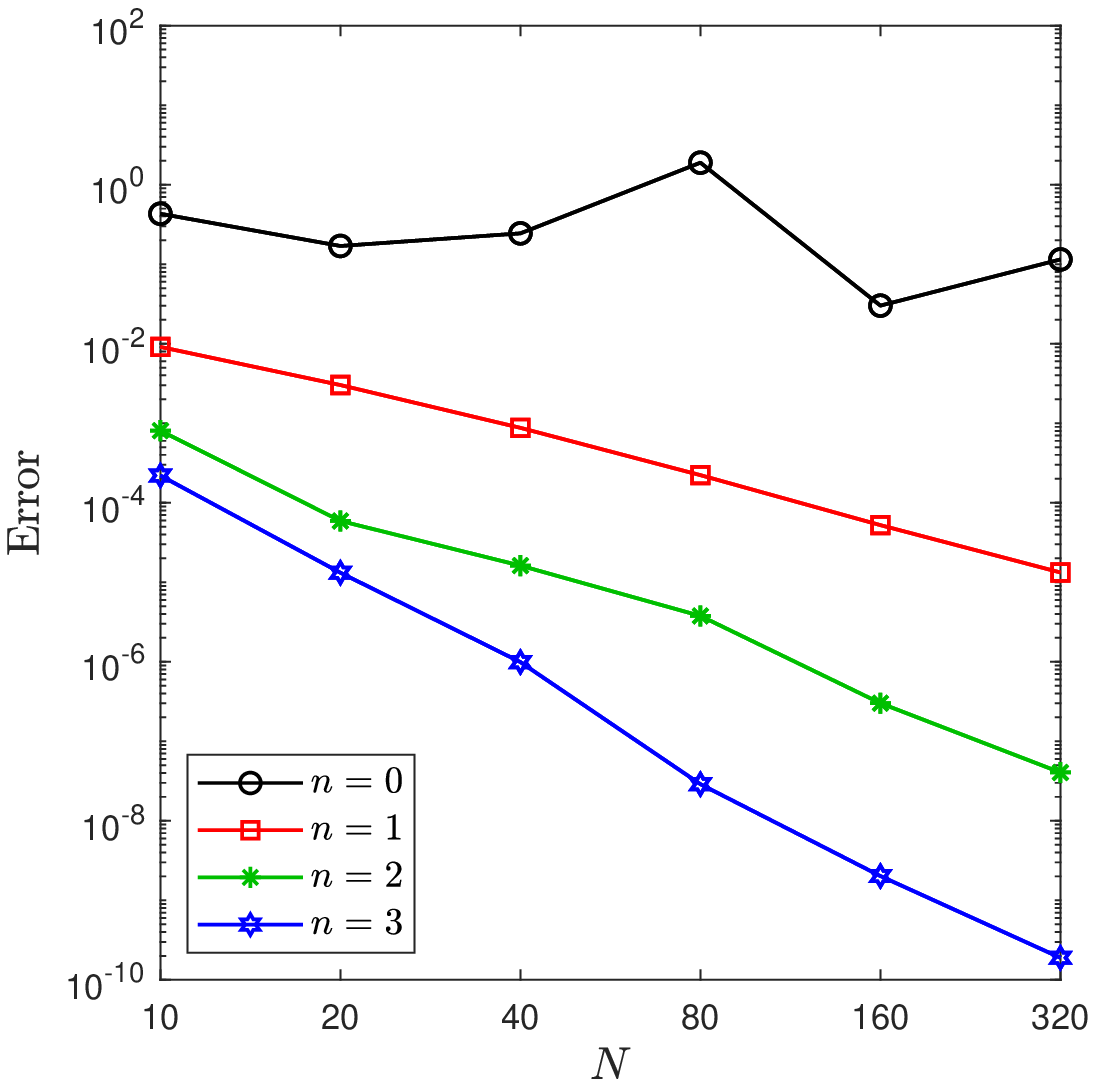} 
}
\subfigure[kernel $\bar \mh_\s$]{
\includegraphics[trim = 0cm .0cm 0cm 0.0cm, clip = true, height=6cm]{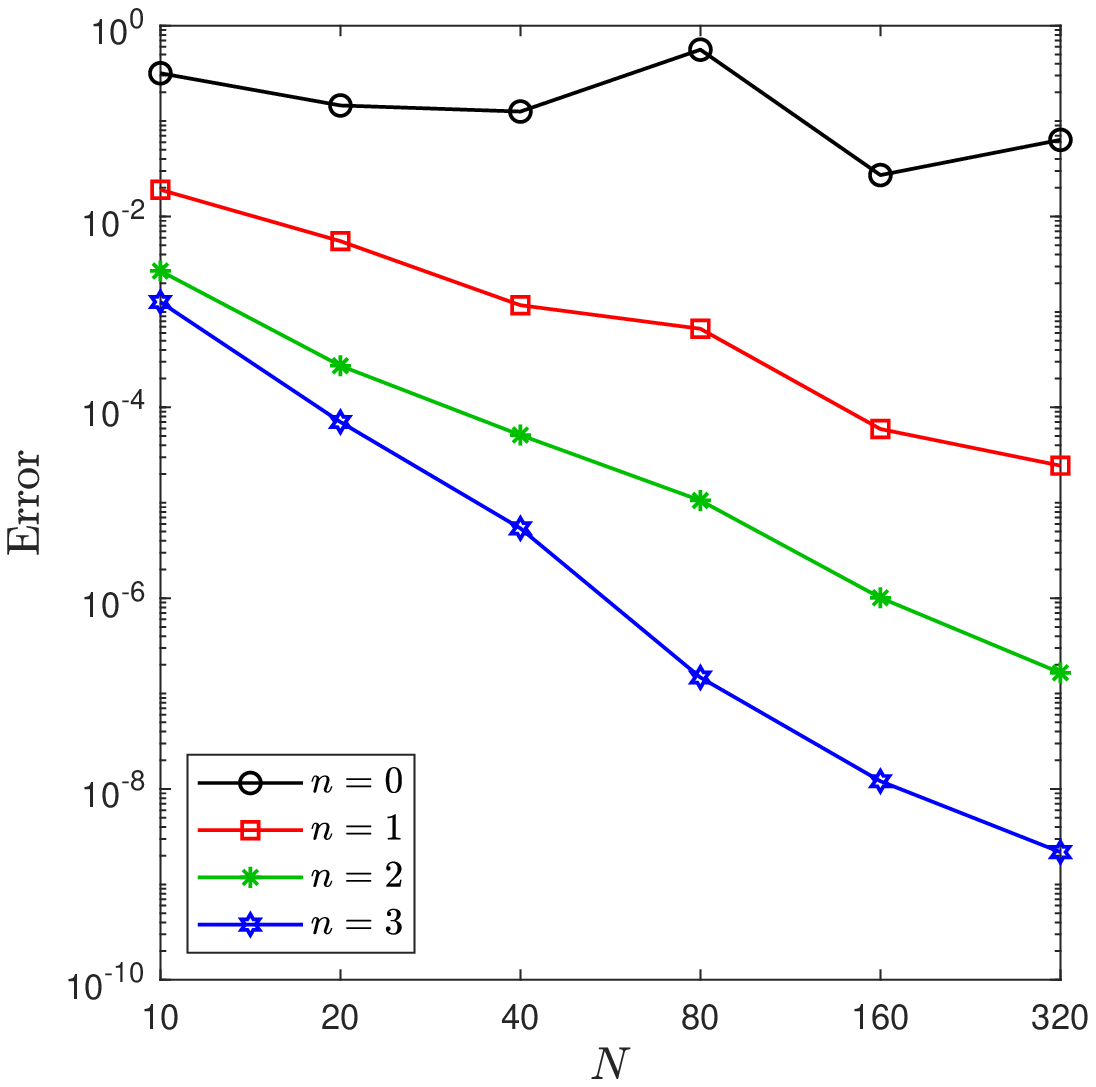} 
}
\caption{Maximum error plots for numerical integration involving kernels $\mg_\s$ and $\bar \mh_\s$ using singularity subtraction and $N$-refinement.}
\label{fig:test1Conv}
\end{figure}

After the subtraction of the singularity in \eqref{eqn:test1}, the singular part of the integral is computed analytically, using the recursive formulae in Section~\ref{sec:recursiveFormulae}. The regular part of the integral is computed numerically, using the tensor product Gauss-Legendre quadrature. For each iteration step we double the number of quadrature nodes $N$ in each direction. Since the initial $N=10$ is already sufficiently high, the accuracy of the integration is expected to be governed by the lower smoothness of the integrand at $\s$. The numerical values are compared against the values obtained by highly accurate numerical scheme that combines the Duffy transformation \cite{Duffy} and Matlab's adaptive quadrature routines. The regularization via division, however, cannot be used in this test in such a simple manner without a further approximation of the integrand and is thus excluded from the test. In Fig.~\ref{fig:test1Conv} by examining the error plots with respect to $N$ we can clearly see the benefits of using higher $n$ in the regularization. With increasing $n$, the accuracy and the convergence order increase for both type of kernels; for $n=1,2,3$ the error decays as $C N^{-(n+1)}$ for some constant $C>0$.

\subsection{h-refinement}

In this experiment we perform an h-refinement for the discretization space; we vary size $h$ of the support of basis functions $B$ in the integrals
\begin{align*}
\frac{1}{4\pi} \int_{[0, 1]^2}  \mk(\s, \t)  B(\t) v(\t) \dt,
\end{align*}
where again $\mk \in \{\mg, \bar \mh\}$. To analyze the integral of BIE with a constant boundary datum, let us fix $v=J$, when we test kernel $\mg$, and let $v\equiv 1$ for kernel $\bar \mh$ since $J$ is already included in the definition of the kernel.
Let us fix 16 source points $\s=(s_1,s_2)$ with $s_1, s_2  \in \{ 0.6, 0.7, 0.8,  0.9\}$. Let $B$ be a piece-wise constant basis function. To analyze only truly singular integrals, let us for each $\s$ consider ``active'' only the following 5 basis functions $B$, denoted as
\begin{align*}
B_\s^{(i)}(\t) :=  \left\{
\begin{array}{ll}
1, & \t\in D_\s^{(i)}\\
0, & \t\notin D_\s^{(i)}
\end{array}\right.,
\end{align*}
with supports on squares of sizes $h \times h$,
\begin{align*}
D_\s^{(1)} &:= [s_1-h/2, s_1+h/2] \times [s_2-h/2, s_2+h/2],                \\
D_\s^{(2)} &:= [s_1-2h/3, s_1+h/3] \times [s_2-2h/3, s_2+h/3],                \\
D_\s^{(3)} &:= [s_1-2h/3, s_1+h/3] \times [s_2-h/3, s_2+2h/3],                \\
D_\s^{(4)} &:= [s_1-h/3, s_1+2h/3] \times [s_2-2h/3, s_2+h/3],                \\
D_\s^{(5)} &:= [s_1-h/3, s_1+2h/3] \times [s_2-h/3, s_2+2h/3],                
\end{align*}
that are inside the domain $[0,1]^2$.

Following the itinerary in Section~\ref{sec:num_int}, by applying the subtraction of singularity, we derive to formula
\eqref{eqn:subInt}. The first integral is computed with Gauss-Legendre quadrature on domain $D_\s^{(i)}$ with $10 \times 10$ quadrature nodes. In the second one, the regular part of the integrand, i.e., $J$ is approximated with a polynomial of bi-degree $(4,4)$ using a least-squares method for the same $10 \times 10$ evaluation sites.
The second integral of \eqref{eqn:subInt} is thus approximated with an integral \eqref{eqn:simplifiedInt} and we integrate it analytically.

Similarly, for the singularity division, we transform the integral into \eqref{eqn:divInt}. The regular part $\rho$ of the integrand  is approximated with a polynomial using the least-squares approach for the same $10 \times 10$ evaluation sites. The obtained integral is of form \eqref{eqn:simplifiedInt} and we integrate it analytically.

\begin{figure}[ht!]
\centering
\subfigure[kernel $\mg_\s$]{
{\includegraphics[trim = 0cm .0cm 0cm 0.0cm, clip = true, height=6cm]{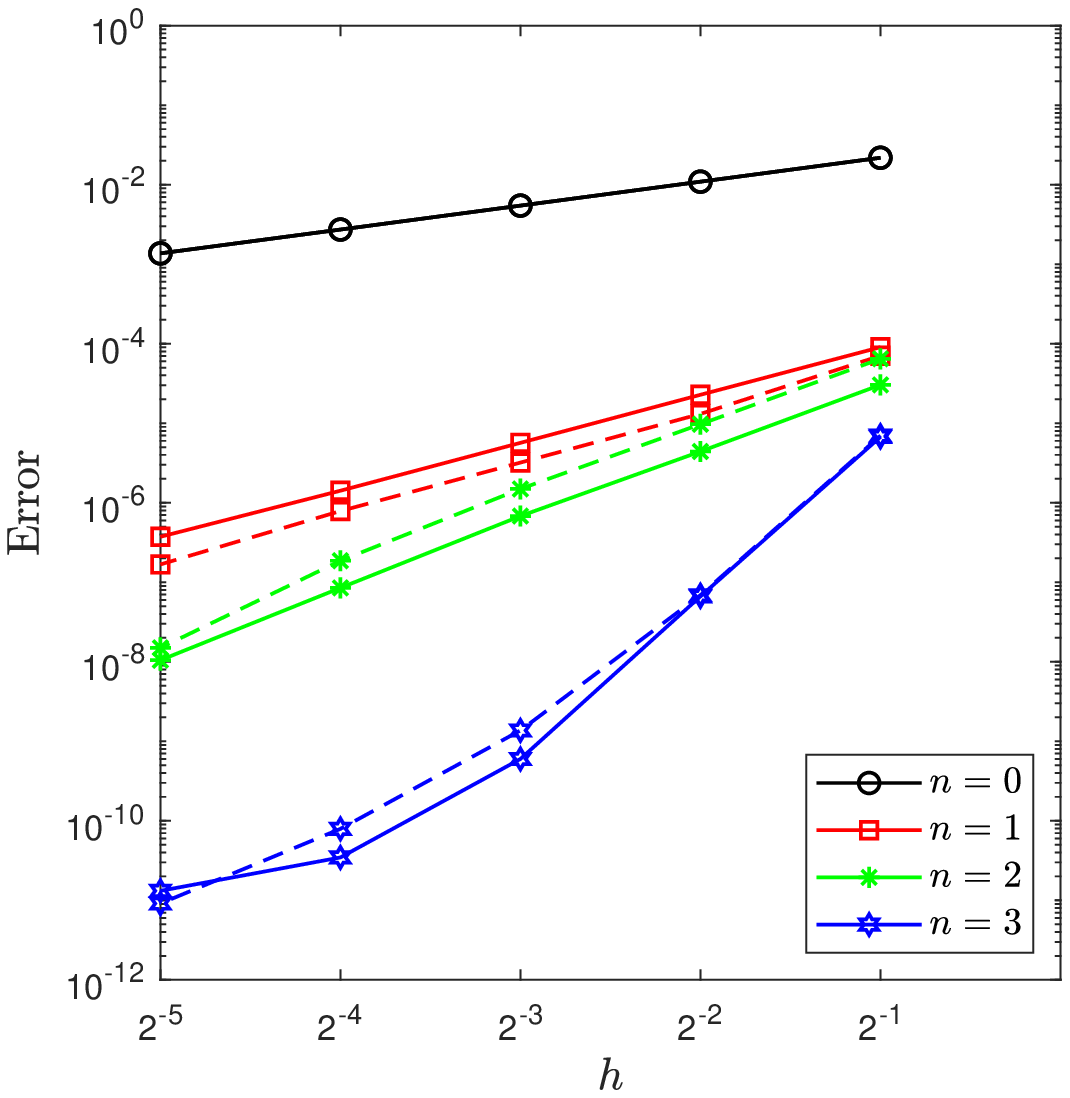} }
}
\subfigure[kernel $\bar \mh_\s$]{
{\includegraphics[trim = 0cm .0cm 0cm 0.0cm, clip = true, height=6cm]{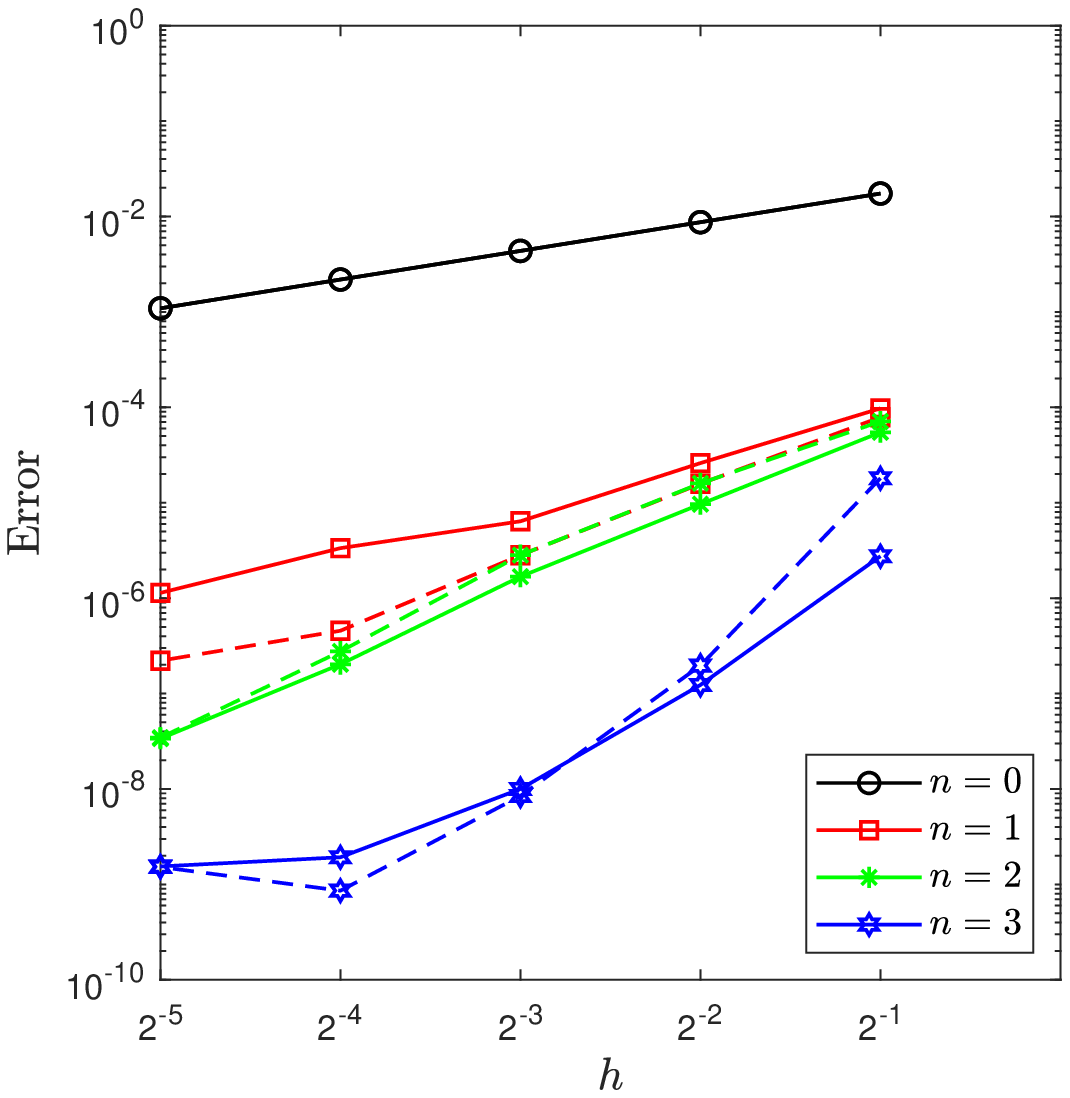} }
}
\caption{Maximum error plots for numerical integration involving kernels $\mg_\s$ and $\bar \mh_\s$ using singularity subtraction and division, and $h$-refinement. Plots with solid and dashed lines correspond to results using the subtraction and division procedure, respectively.}
\label{fig:test2Conv}
\end{figure}

As in the first experiment, we expect the accuracy of integration schemes to be subjected to reduced smoothness of $\rho$ at $\s$. In Fig.~\ref{fig:test2Conv} we can see the maximum error plots for all $\s$ and $D_\s^{(i)}$ for the numerical evaluation of integrals with respect to the support size $h$. Again, there is a clear benefit in using higher $n$ for both singularity subtraction and division. Both regularization techniques produce similar accuracy for both kernels and for all tested $h$.

\section*{Conclusion}

The singularity extraction is studied for singular integrals in BIE for 3D potential problems on smooth geometries, focusing on the single and and double layer operator. Integrals of truncated series expansions of singular kernels are computed via recurrence formulae -- to reduce the complexity of this process, a speedup via lookup tables that exploits correlated variables is provided. 

Numerical tests demonstrate that the presented extraction technique can be a useful prerequisite for a numerical quadrature to recover the optimal order of convergence of the approximate solution in BEM with a small number of quadrature nodes. Since the extraction acts directly on the starting parametric domain, no modification of integration domain is needed -- this is an advantageous property since the same quadrature nodes can be reused for several integrals involving neighbouring trial basis functions. The truncation error of the numerical integration rule can be attributed solely to the truncation error of approximating the regular part of the integrand; thus the accuracy of the rule can be effectively controlled by the number of terms in the series extraction and by the choice of the approximation operator for regular functions. Singularity subtraction and division give comparable results; the technique via division is arguably slightly more efficient since the governing integral does not need to be split into two integrals, however it cannot directly handle kernels that change signs. 

Interesting topics for future work include a study of other kernels with stronger types of singularities and of non-stationary nature. Efficient implementation for computing singular integrals via lookup tables is beyond the scope of this paper, but it should be investigated in future work. 
 Finding better ways to approximate kernels across non-smooth interfaces of the geometry, for examples on multipatch domains, could improve accuracy of the integration routines on more complex geometries.

\section*{References}
\bibliography{BEM3D_SE_references}

\begin{thebibliography}{54}
\expandafter\ifx\csname natexlab\endcsname\relax\def\natexlab#1{#1}\fi
\expandafter\ifx\csname url\endcsname\relax
  \def\url#1{\texttt{#1}}\fi
\expandafter\ifx\csname urlprefix\endcsname\relax\def\urlprefix{URL }\fi

\bibitem[{Aimi et~al.(2018)Aimi, Calabr\`o, Diligenti, Sampoli, Sangalli, and
  Sestini}]{ACDS3}
Aimi, A., Calabr\`o, F., Diligenti, M., Sampoli, M.~L., Sangalli, G., Sestini,
  A., 2018. Efficient assembly based on {B}-spline tailored quadrature rules
  for the {IgA-SGBEM}. Comput. Methods Appl. Mech. Engrg. 331, 327--342.

\bibitem[{Aimi et~al.(2020)Aimi, Calabr\`{o}, Falini, Sampoli, and
  Sestini}]{AIMI2020113441}
Aimi, A., Calabr\`{o}, F., Falini, A., Sampoli, M.~L., Sestini, A., 2020.
  Quadrature formulas based on spline quasi-interpolation for hypersingular
  integrals arising in {IgA-SGBEM}. Comput. Methods Appl. Mech. Engrg. 372,
  113441.

\bibitem[{Aimi and Diligenti(2002)}]{AD2002}
Aimi, A., Diligenti, M., 2002. Numerical integration in 3{D} {G}alerkin {BEM}
  solution of {HBIE}s. Comput. Mech. 28, 233--249.

\bibitem[{Atkinson(2009)}]{Atkinson_2009}
Atkinson, K.~E., 2009. The {N}umerical {S}olution of {I}ntegral {E}quations of
  the {S}econd {K}ind. Cambridge University Press.

\bibitem[{Beer et~al.(2015)Beer, Marussig, and Zechner}]{BEER2015776}
Beer, G., Marussig, B., Zechner, J., 2015. A simple approach to the numerical
  simulation with trimmed {CAD} surfaces. Comput. Methods Appl. Mech. Engrg.
  285, 776--790.

\bibitem[{Botha(2013)}]{Duffy4}
Botha, M.~M., 2013. A family of augmented {D}uffy transformations for
  near-singularity cancellation quadrature. IEEE Antennas Wirel. Propag. Lett.
  61~(6), 3123--3134.

\bibitem[{Calabr\`o et~al.(2018)Calabr\`o, Falini, Sampoli, and
  Sestini}]{CFSS18}
Calabr\`o, F., Falini, A., Sampoli, M.~L., Sestini, A., 2018. Efficient
  quadrature rules based on spline quasi-interpolation for application to
  {IgA-BEM}s. J. Comput. Appl. Math. 338, 153--167.

\bibitem[{Calabr\`o et~al.(2019)Calabr\`o, Loli, Sangalli, and
  Tani}]{IGA_quads}
Calabr\`o, F., Loli, G., Sangalli, G., Tani, M., 2019. Quadrature rules in the
  isogeometric {G}alerkin method: {S}tate of the art and an introduction to
  weighted quadrature. In: Giannelli, C., Speleers, H. (Eds.), Advanced Methods
  for Geometric Modeling and Numerical Simulation. Vol.~35 of Springer INdAM
  Series. Springer Cham, pp. 43--55.

\bibitem[{Chen and Hong(1999)}]{Chen_Hong_Review_dualBEM}
Chen, J.~T., Hong, H.-K., 1999. Review of dual boundary element methods with
  emphasis on hypersingular integrals and divergent series. Appl. Mech. Rev.
  52~(1), 17--33.

\bibitem[{Cobb(1988)}]{Cobb1994TilingTS}
Cobb, J.~E., 1988. Tiling the sphere with rational {B}\'{e}zier patches.
  Technical report {UUCS}-88-009, Computer Science, University of Utah.

\bibitem[{Costabel(1986)}]{costabel1986principles}
Costabel, M., 1986. Principles of boundary element methods. Techn. Hochsch.,
  Fachbereich Mathematik.

\bibitem[{Cottrell et~al.(2009)Cottrell, Hughes, and Bazilevs}]{LibroHughes}
Cottrell, J.~A., Hughes, T.~J.~R., Bazilevs, Y., 2009. Isogeometric analysis:
  toward integration of {CAD} and {FEA}. John Wiley \& Sons.

\bibitem[{de~Boor(2001)}]{librodeBoor01}
de~Boor, C., 2001. A practical guide to splines, revised Edition. Vol.~27 of
  Applied Mathematical Sciences. Springer-Verlag, New York.

\bibitem[{D\"olz et~al.(2018)D\"olz, Harbrecht, Kurz, Sch\"ops, and
  Wolf}]{Dolz18}
D\"olz, J., Harbrecht, H., Kurz, S., Sch\"ops, S., Wolf, F., 2018. A fast
  isogeometric {BEM} for the three dimensional {L}aplace and {H}elmholtz
  problems. Comput. Methods Appl. Mech. Engrg. 330, 83--101.

\bibitem[{Duffy(1982)}]{Duffy}
Duffy, M.~G., 1982. Quadrature over a pyramid or cube of integrands with a
  singularity at a vertex. SIAM J. Numer. Anal. 19~(6), 1260--1262.

\bibitem[{Falini et~al.(2022)Falini, Giannelli, Kandu{\v{c}}, Sampoli, and
  Sestini}]{FGKSS_2021}
Falini, A., Giannelli, C., Kandu{\v{c}}, T., Sampoli, M.~L., Sestini, A., 2022.
  A collocation {IGA-BEM} for {3D} potential problems on unbounded domains. In:
  Manni, C., Speleers, H. (Eds.), Springer INdAM volume ‘Geometric Challenges
  in Isogeometric Analysis’. Vol.~49 of Springer INdAM Series. Springer
  Cham., p. In press.

\bibitem[{Falini et~al.(2019)Falini, Giannelli, Kandu\v{c}, Sampoli, and
  Sestini}]{FGKSS_2019}
Falini, A., Giannelli, C., Kandu\v{c}, T., Sampoli, M.~L., Sestini, A., 2019.
  An adaptive {IgA-BEM} with hierarchical {B}-splines based on
  quasi-interpolation quadrature schemes. Int. J. Numer. Methods Engrg.
  117~(10), 1038--1058.

\bibitem[{Falini et~al.(2020)Falini, Kandu{\v{c}}, Sampoli, and
  Sestini}]{nash20}
Falini, A., Kandu{\v{c}}, T., Sampoli, M.~L., Sestini, A., 2020. Cubature rules
  based on bivariate spline quasi-interpolation for weakly singular integrals.
  In: Fasshauer, G.~E., Neamtu, M., Schumaker, L.~L. (Eds.), Approximation
  Theory XVI. AT 2019. Vol. 336 of Springer in Mathematics \& Statistics.
  Springer Cham, pp. 73--86.

\bibitem[{Falini and Kandu\v{c}(2019)}]{FKdreams}
Falini, A., Kandu\v{c}, T., 2019. A study on spline quasi-interpolation based
  quadrature rules for the isogeometric {G}alerkin {BEM}. In: Giannelli, C.,
  Speleers, H. (Eds.), Advanced Methods for Geometric Modeling and Numerical
  Simulation. Vol.~5 of Springer INdAM Series. Springer, pp. 193--227.

\bibitem[{Falini et~al.(2021)Falini, Kandu\v{c}, Sampoli, and
  Sestini}]{FKSS_2021}
Falini, A., Kandu\v{c}, T., Sampoli, M.~L., Sestini, A., 2021. Isogeometric
  {BEM} collocation for {3D} {L}aplace and {H}elmholtz problems. in
  preparation.

\bibitem[{Farin(2002)}]{Farin-02-CAGD}
Farin, G., 2002. Curves and surfaces for computer-aided geometric design, 5th
  Edition. Computer Graphics and Geometric Modeling. Academic Press Inc., San
  Diego, CA.

\bibitem[{Fata(2009)}]{Fata2009ExplicitEF}
Fata, S.~N., 2009. Explicit expressions for {3D} boundary integrals in
  potential theory. Int. J. Numer. Meth. Eng. 78, 32--47.

\bibitem[{Feischl et~al.(2017)Feischl, Gantner, Haberl, and Praetorius}]{Dirk1}
Feischl, M., Gantner, G., Haberl, A., Praetorius, D., 2017. Optimal convergence
  for adaptive {IGA} boundary element methods for weakly-singular integral
  equations. Numer. Math. 136, 147--182.

\bibitem[{Gantner and Praetorius(2020)}]{Dirk2}
Gantner, G., Praetorius, D., 2020. Adaptive {BEM} for elliptic {PDE} systems,
  part {I}: abstract framework, for weakly-singular integral equations. Appl.
  Anal., 1--34.

\bibitem[{Gao(2002)}]{GAO2002905}
Gao, X.-W., 2002. The radial integration method for evaluation of domain
  integrals with boundary-only discretization. Eng. Anal. Bound. Elem. 26~(10),
  905--916.

\bibitem[{Gao(2010)}]{gao_cmame}
Gao, X.-W., 2010. An effective method for numerical evaluation of general {2D}
  and {3D} high order singular boundary integrals. Comput. Methods Appl. Mech.
  Engrg. 199~(45), 2856--2864.

\bibitem[{Ginnis et~al.(2014)Ginnis, Kostas, Politis, Kaklis, Belibassakis,
  Gerostathis, Scott, and Hughes}]{GINNIS2014425}
Ginnis, A., Kostas, K., Politis, C., Kaklis, P., Belibassakis, K., Gerostathis,
  T., Scott, M., Hughes, T., 2014. Isogeometric boundary-element analysis for
  the wave-resistance problem using {T}-splines. Comput. Methods Appl. Mech.
  Engrg. 279, 425--439.

\bibitem[{Gradshteyn and Ryzhik(2007)}]{integralsBook}
Gradshteyn, I.~S., Ryzhik, I.~M., 2007. Table of integrals, series, and
  products, 7th Edition. Academic Press, Elsevier.

\bibitem[{Guiggiani et~al.(1992)Guiggiani, Krishnasamy, Rudolphi, and
  Rizzo}]{Guiggiani_generalAlg_HBIE}
Guiggiani, M., Krishnasamy, G., Rudolphi, T.~J., Rizzo, F.~J., 09 1992. A
  general algorithm for the numerical solution of hypersingular boundary
  integral equations. J. Appl. Mech. 59~(3), 604--614.

\bibitem[{Heltai et~al.(2014)Heltai, Arroyo, and
  DeSimone}]{heltai2014nonsingular}
Heltai, L., Arroyo, M., DeSimone, A., 2014. Nonsingular isogeometric boundary
  element method for {S}tokes flows in {3D}. Comput. Methods Appl. Mech. Engrg.
  268, 514--539.

\bibitem[{Heltai et~al.(2017)Heltai, Kiendl, DeSimone, and
  Reali}]{HELTAI2017522}
Heltai, L., Kiendl, J., DeSimone, A., Reali, A., 2017. A natural framework for
  isogeometric fluid–structure interaction based on {BEM}–shell coupling.
  Comput. Methods Appl. Mech. Engrg. 316, 522--546.

\bibitem[{Hughes et~al.(2005)Hughes, Cottrell, and Bazilevs}]{Hughes_2005}
Hughes, T.~J.~R., Cottrell, J.~A., Bazilevs, Y., 2005. {Isogeometric analysis:
  CAD, finite elements, NURBS, exact geometry and mesh refinement}. Comput.
  Methods Appl. Mech. Engrg. 194~(39-41), 4135--4195.

\bibitem[{J\"{a}rvenp\"{a}\"{a} et~al.(2003)J\"{a}rvenp\"{a}\"{a}, Taskinen,
  and Yl\"{a}-Oijala}]{jarvenpaa}
J\"{a}rvenp\"{a}\"{a}, S., Taskinen, M., Yl\"{a}-Oijala, P., 2003. Singularity
  extraction technique for integral equation methods with higher order basis
  functions on plane triangles and tetrahedra. Int. J. Numer. Meth. Eng.
  58~(8), 1149--1165.

\bibitem[{Johnston et~al.(2013)Johnston, Johnston, and
  Elliott}]{JOHNSTON2013148}
Johnston, B.~M., Johnston, P.~R., Elliott, D., 2013. A new method for the
  numerical evaluation of nearly singular integrals on triangular elements in
  the {3D} boundary element method. J. Comput. Appl. Math. 245, 148--161.

\bibitem[{Khayat et~al.(2008)Khayat, Wilton, and Fink}]{Khayat_polarSingInteg}
Khayat, M.~A., Wilton, D.~R., Fink, P.~W., 2008. An improved transformation and
  optimized sampling scheme for the numerical evaluation of singular and
  near-singular potentials. IEEE Antennas Wirel. Propag. Lett. 7, 377--380.

\bibitem[{Klaseboer et~al.(2009)Klaseboer, Fernandez, and
  Khoo}]{KLASEBOER2009796}
Klaseboer, E., Fernandez, C.~R., Khoo, B.~C., 2009. A note on true
  desingularisation of boundary integral methods for three-dimensional
  potential problems. Eng. Anal. Bound. Elem. 33~(6), 796--801.

\bibitem[{Kostas et~al.(2018)Kostas, Fyrillas, Politis, Ginnis, and
  Kaklis}]{KOSTAS2018600}
Kostas, K., Fyrillas, M., Politis, C., Ginnis, A., Kaklis, P., 2018. Shape
  optimization of conductive-media interfaces using an {IGA-BEM} solver.
  Comput. Methods Appl. Mech. Engrg. 340, 600--614.

\bibitem[{Mazzia and Sestini(2009)}]{MSbit09}
Mazzia, F., Sestini, A., 2009. The {BS} class of {H}ermite spline
  quasi-interpolants on nonuniform knot distributions. BIT 49~(3), 611--628.

\bibitem[{Mazzia and Sestini(2012)}]{MSJcam12}
Mazzia, F., Sestini, A., 2012. Quadrature formulas descending from {BS}
  {H}ermite spline quasi-interpolation. J. Comput. Appl. Math. 236, 4105--4118.

\bibitem[{Monegato and Sloan(1997)}]{Monegato_airfol}
Monegato, G., Sloan, I.~H., 1997. Numerical solution of the generalized airfoil
  equation for an airfoil with a flap. SIAM J. Numer. Anal. 34~(6), 2288--2305.

\bibitem[{M{\o}rken(1991)}]{Morken91}
M{\o}rken, K., 1991. Some identities for products and degree raising of
  splines. Constr. Approx. 7, 195--208.

\bibitem[{Mousavi and Sukumar(2010)}]{Duffy2}
Mousavi, S.~E., Sukumar, N., 2010. Generalized {D}uffy transformation for
  integrating vertex singularities. Comput. Mech. 45, 127--140.

\bibitem[{Peng et~al.(2017)Peng, Atroshchenko, Kerfriden, and
  Bordas}]{PENG2017151}
Peng, X., Atroshchenko, E., Kerfriden, P., Bordas, S., 2017. Isogeometric
  boundary element methods for three dimensional static fracture and fatigue
  crack growth. Comput. Methods Appl. Mech. Engrg. 316, 151--185, special Issue
  on Isogeometric Analysis: Progress and Challenges.

\bibitem[{Reid et~al.(2015)Reid, White, and Johnson}]{Duffy3}
Reid, M. T.~H., White, J.~K., Johnson, S.~G., 2015. Generalized
  {T}aylor-{D}uffy method for efficient evaluation of {G}alerkin integrals in
  boundary-element method computations. IEEE Antennas Wirel. Propag. Lett.
  63~(1), 195--209.

\bibitem[{Rong et~al.(2014)Rong, Wen, and Xiao}]{Rong_polarTransform}
Rong, J., Wen, L., Xiao, J., 2014. Efficiency improvement of the polar
  coordinate transformation for evaluating {BEM} singular integrals on curved
  elements. Eng. Anal. Bound. Elem. 38, 83--93.

\bibitem[{Sauter and Schwab(2011)}]{BEMbook}
Sauter, S.~A., Schwab, C., 2011. Boundary element methods. Vol.~39 of Springer
  Series in Computational Mathematics. Springer-Verlag, Berlin, Heidelberg.

\bibitem[{Schumaker(2007)}]{schumaker2007spline}
Schumaker, L.~L., 2007. Spline functions: basic theory, 3rd Edition. Cambridge
  University Press.

\bibitem[{Scuderi(2009)}]{SCUDERI2009406}
Scuderi, L., 2009. A new smoothing strategy for computing nearly singular
  integrals in {3D} {G}alerkin {BEM}. J. Comput. Appl. Math. 225~(2), 406--427.

\bibitem[{Simpson et~al.(2014)Simpson, Scott, Taus, Thomas, and
  Lian}]{SimpsonScott}
Simpson, R., Scott, M., Taus, M., Thomas, D., Lian, H., 2014. Acoustic
  isogeometric boundary element analysis. Comput. Methods Appl. Mech. Engrg.
  269, 265--290.

\bibitem[{Taus et~al.(2019)Taus, Rodin, Hughes, and Scott}]{TAUS2019112591}
Taus, M., Rodin, G.~J., Hughes, T.~J., Scott, M.~A., 2019. Isogeometric
  boundary element methods and patch tests for linear elastic problems:
  Formulation, numerical integration, and applications. Comput. Methods Appl.
  Mech. Engrg. 357, 112591.

\bibitem[{Taus et~al.(2016)Taus, Rodin, and Hughes}]{TauRodHug}
Taus, M., Rodin, G.~J., Hughes, T. J.~R., 2016. Isogeometric analysis of
  boundary integral equations: {H}igh-order collocation methods for the
  singular and hyper-singular equations. Math. Models and Methods in Appl. Sci.
  26~(8), 1447--1480.

\bibitem[{Telles(1987)}]{telles1987self}
Telles, J. C.~F., 1987. A self-adaptive co-ordinate transformation for
  efficient numerical evaluation of general boundary element integrals. Int. J.
  Numer. Meth. Eng. 24~(5), 959--973.

\bibitem[{Ven{\aa}s and Kvamsdal(2020)}]{VENAS2020112670}
Ven{\aa}s, J.~V., Kvamsdal, T., 2020. Isogeometric boundary element method for
  acoustic scattering by a submarine. Comput. Methods Appl. Mech. Engrg. 359,
  112670.

\bibitem[{Zhang et~al.(2015)Zhang, Li, Sladek, Sladek, and Gao}]{ZHANG201557}
Zhang, Y., Li, X., Sladek, V., Sladek, J., Gao, X., 2015. A new method for
  numerical evaluation of nearly singular integrals over high-order geometry
  elements in {3D} {BEM}. J. Comput. Appl. Math. 277, 57--72.

\end{thebibliography}

\appendix
\section{The underlying integrals}\label{basicInt}

In this section we analyze some basic properties of the integrand functions for the underlying integrals in intrinsic coordinates. Then we provide recursive formulae to evaluate these kind of integrals analytically. 
Following the isoparametric paradigm, the most common surface integrals are defined either on rectangular or triangular domains. The definite integrals are computed as a linear combination of indefinite ones, evaluated at appropriate boundary points. 

Suppose $f$ is a general integrable bivariate function we would like to integrate analytically. For both type of domains we need a primitive function of $f$,  $F(x,y) := \int \int f(x,y) \dx \dy$. Since indefinite integrals are not uniquely defined, we can neglect all byproducts that do not contribute to the final values of the definite integrals, e.g., functions $g$ that satisfy $\partial^2 g  / (\mathop{\partial x}  \mathop{\partial y} ) \equiv 0$. For the triangular domain let us also define $\cc{F^{[\rm t]}(y)} := \int (\int f(x,y) \dx )|_{x=1-y}\dy$.
The surface integral on rectangle $[x_0, x_1] \times [y_0, y_1]$ can be simply computed by evaluating the indefinite integral at the 4 corners,
\begin{align*}
\int_{y_0}^{y_1} \int_{x_0}^{x_1} f(x,y) \dx \dy = F(x_1,y_1) - F(x_0,y_1) - F(x_1,y_0) + F(x_0,y_0).
\end{align*}
To compute the surface integral of $\tilde f$ on a general triangle $T$, a preliminary step is to apply an affine mapping $\varphi$ to map, for example, a triangle with vertices $(0,0),(1,0), (0,1)$ to $T$. Then we need to again insert appropriate boundary values into the two indefinite integrals,
\begin{align*}
\int_T \tilde f(\tilde x, \tilde y) \, dT =  
\int_{0}^{1} \int_{0}^{1-y} f(x,y) \dx \dy = F(0,0) -F(0,1)  + F^{[\rm t]}(1) - F^{[\rm t]}(0),
\end{align*}
where $f := (\tilde f \circ \varphi) \cdot J$ and constant $J$ is the area of the infinitesimal surface element.

The two derived formulae allow us to analytically compute surface integrals on a rectangle and triangle, if the primitive function of $f$ is known. In the remaining part of the section we focus only on a special family of these functions that are the building blocks for the singularity extraction technique, and provide formulae for their indefinite integrals.

\subsection{Main properties}

Our goal in this section is to derive recursive formulae to compute integrals of a type
\begin{align}\label{eqn:fullInt}
I_{p,q,r} = I_{p,q,r}(x,y):=\int \int R^{p} x^q y^r \dx \dy,
\end{align}
where $R = R(a,b,c,(x,y)) := \sqrt{a x^2 + b x y + c y^2}$ is a function, corresponding to the first fundamental form of a smooth surface. Here we assume $a,c>0$ and $b\in\RR$, such that $4ac - b^2 > 0$, which ensures the quadratic function to be positive and hence $R > 0$. Then, $p$ is a negative odd integer, while $q,r$ are non-negative integers. Double integrals \eqref{eqn:fullInt} can be computed via analytical formulae for single inner and outer integrals,
\begin{align}\label{eqn:general_intdx}
I_{p,q}^{(1)} = I_{p,q}^{(1)}(x):=\int R^{p} x^q \dx,
\qquad 
I_{p,r}^{(2)} = I_{p,r}^{(2)}(y):=\int R^{p} y^r \dy,
\end{align} namely $I_{p,q,r} = \int y^r I_{p,q}^{(1)}  \dy.$
To compute integrals on a triangle, we also need formulae to compute
\begin{align}\label{eqn:fullInt_t}
I_{p,q,r}^{[\rm t]} = I_{p,q,r}^{[\rm t]} (y):= \int  y^r (I_{p,q}^{(1)}|_{x=1-y}) \dy.
\end{align}

\begin{rmk}\label{rmk:exchange}
Since $x$ and $y$ appear symmetrically in \eqref{eqn:fullInt} it is sufficient to study analytical formuale for $I_{p,q}^{(1)}$ only. Formulae for $I_{p,r}^{(2)}$ are therefore simply transformed from the derived expressions for $I_{p,q}^{(1)}$  by exchanging variables $x$ and $y$, and parameters $a$ and $c$, and $q$ and $r$.
\end{rmk}

Smoothness (regularity) of the integrand $R^p x^q y^r $ in \eqref{eqn:fullInt} at $(x,y)=(0,0)$ is easily characterized by the parameter $\zeta$,
\begin{align}\label{eqn:zeta}
\zeta(R^p x^q y^r):=p+q+r.
\end{align}
For the study of the smoothness we relax the condition on $p$ to be an arbitrary integer. When we consider a sum of integrands of this type, it is also convenient to use the following definition
\begin{align*}
\zeta \left( \sum_\ell c_\ell R^{p_\ell} x^{q_\ell} y^{r_\ell} \right) := \min_\ell \zeta(R^{p_\ell} x^{q_\ell} y^{r_\ell}),
\end{align*}
where we assume that the sum is written in irreducible form for some coefficients $c_\ell$; it is also meaningful to set $\zeta = \infty$ if we sum over an empty set of indices. 

Without a proof we mention two simple properties of the parameter $\zeta$.
\begin{prn}\label{pr:operations}
Let $f,g$ be functions of a type $R^p x^q y^r $.
Then 
\begin{itemize}
\item $\zeta(f g) = \zeta(f) + \zeta(g)$,
\item if $\zeta(f)=m$, then $n$-th derivative of $f$ can be written as ${\partial^n f}/{(\mathop{\partial x}^i \mathop{\partial y}^{n-i})} = \sum_{\ell=0}^n f_\ell$, $i=0,1,\dots,n$, where $f_\ell$ are also of a type $R^p x^q y^r $ and $\zeta(f_\ell) = m-n$ for every $\ell$.
\end{itemize}
\end{prn}

Before analyzing the smoothness of the integrands, we need to prove the following lemma.

\begin{lma}\label{lem:lemma1}
Function $ \gamma \mapsto (a  + b \gamma + c \gamma^2)^{p/2}$ is positive and bounded for $|\gamma| \leq 1$. Namely, there exist two positive constants $ C_{abcp,1},  C_{abcp,2}$, independent of $\gamma$, such that
\begin{align*}
 C_{abcp,1} \leq (a  + b \gamma + c \gamma^2)^{p/2} \leq  C_{abcp,2}.
\end{align*}
\end{lma}

\begin{pf}
There exist two positive constants $\tilde C_{abc,1},\tilde C_{abc,2}$, independent of $\gamma$, such that
\begin{align}\label{eqn:lemma1}
\tilde C_{abc,1} \leq {a  + b \gamma + c \gamma^2} \leq \tilde C_{abc,2}.
\end{align}
The upper bound $\tilde C_{abc,2}$ exists since $\gamma \mapsto a  + b \gamma + c \gamma^2$ is continuous and thus bounded on $|\gamma|\leq1$. The existence of the lower bound $\tilde C_{abc,1}$ follows from assumption $0< 4ac - b^2$.
By raising \eqref{eqn:lemma1} to the power $p/2$ we get the sought expression.  \qed
\end{pf}

\begin{prn}\label{pr:bounded}
Function $R^p x^q y^r$ is bounded if $\zeta(R^p x^q y^r) = 0$.
\end{prn}

\begin{pf}
Let us split $\RR^2$ into eight subdomains, formed by cutting the plane with lines $x=0$, $y=0$, $y=x$ and $y=-x$. Let us prove the proposition only for the subdomain
\begin{align}\label{eqn:subdomains}
T_1:=\{(x,y)\in\RR^2: 0  \leq y \leq x\}, 
\end{align}
since for the other 7 cases the proof is similar.

Let us choose a point $(x,y) \in T_1$. By writing $y= \gamma x$ for $0 \leq \gamma \leq 1$ we get
\begin{align*}
R^{p} x^q y^r  = {(a  + b \gamma + c \gamma^2)}^{p/2}  \gamma^r  x^{p+q+r}.
\end{align*}
By considering $p+q+r=0$, from Lemma~\ref{lem:lemma1} there exist positive constants $C_{abcp,1}, C_{abcp,2}$ such that
\begin{align*}
0 \leq C_{abcp,1} \gamma^r \leq {(a  + b \gamma + c \gamma^2)}^{p/2}\,  \gamma^r \leq C_{abcp,2} \gamma^r \leq C_{abcp,2}
\end{align*}
for all $0 \leq \gamma \leq 1$. Since this estimate holds for every $(x,y) \in T_1$, the function $R^{p} x^q y^r$ is bounded on $T_1$.


\qed
\end{pf}

\begin{prn}\label{pr:cont}
Let  $f(x,y):=R^p x^q y^r$  with $\zeta(f)\geq 1$. Then $f$ is $C^0$ continuous and $f(0,0) =0$.
\end{prn}

\begin{pf}
The function is analytical on $\RR^2 \backslash \{(0,0)\}$, hence it is enough to check it is continuous at $(0,0)$. Fix sufficiently small $0 < \delta$. Let us show that there exists sufficiently small $\varepsilon$, $0 \leq \varepsilon \leq 1$, such that from $\| (x, y) \|_2 \leq \varepsilon$ it follows $|R^p x^q y^r| \leq \delta$.

By splitting $\RR^2$ again into eight subdomains and assuming $(x,y)\in T_1$ (see \eqref{eqn:subdomains}), we can write
$y= \gamma x$ for $0 \leq \gamma \leq 1$. Then we can again write the function as
\begin{align*}
R^{p} x^q y^r  = {(a  + b \gamma + c \gamma^2)}^{p/2}  \gamma^r  x^{p+q+r}.
\end{align*}
From Lemma~\ref{lem:lemma1} there exists a positive constant $C_{abcp,2}$, independent of $\gamma$, and we can take $\varepsilon$ small enough, so that
\begin{align*}
R^{p} x^q y^r  \leq C_{abcp,2}\, \varepsilon^{p+q+r} \leq C_{abcp,2} \, \varepsilon \leq \delta.
\end{align*}

A similar proof can be asserted for the other 7 subdomains.

\qed
\end{pf}

\begin{prn}\label{pr:smooth}
Let $f(x,y):=R^p x^q y^r $ with $\zeta(f) \geq n+1$.  Then $f$ is $C^n$ continuous and ${\partial^n f}/{(\mathop{\partial x}^i {\mathop{\partial y}}^{n-i})} (0,0) = 0$ for $i=0,1,\dots,n$.
\end{prn}

\begin{pf}


Let $\zeta(f) = p+q+r \geq 2$. It is straightforward to check that
\begin{align}\label{eqn:der_Rxy}
\frac{\partial}{\partial x} (R^p x^q y^r) = p a R^{p-2} x^{q+1}y^r + \frac{pb}{2} R^{p-2} x^q y^{r+1}.
\end{align}
From Proposition \ref{pr:cont} it follows that both functions on the right-hand side of \eqref{eqn:der_Rxy} are continuous. A similar argument can be stated for $\partial (R^p x^q y^r) / {\partial y}$. Thus $R^p x^q y^r$ is a $C^1$ smooth function.

If $p+q+r \geq n+1$, we can apply formula \eqref{eqn:der_Rxy}, or its variation for the derivative in variable $y$, $n$ times and get
\begin{align}\label{eqn:der_Rxy2}
\frac{\partial^n}{\mathop{\partial x}^i \mathop{\partial y}^{n-i}} (R^p x^q y^r) = \sum_{j=0}^n c_j R^{p-2n} x^{q+j} y^{r+n-j}
\end{align}
for suitable coefficients $c_j$ and $i=0,1,\dots,n$. Since $\zeta(R^{p-2n} x^{q+j} y^{r+n-j}) \geq 1 $ and by applying Proposition~\ref{pr:cont}, functions on the right-hand side of \eqref{eqn:der_Rxy2} are continuous for all mixed $n$-th derivatives, hence $R^p x^q y^r$ is a $C^n$ smooth function. From Proposition~\ref{pr:cont} also follows directly the condition for the all $n$-th derivatives to be zero at $(0,0)$.  \qed
\end{pf}

\subsection{Recursive formulae for $I_{p,q,r}$ and $I_{p,q,r}^{[\rm t]}$}
\label{sec:recursiveFormulae}

Formulae for double integrals are derived in two steps. First, we focus on the {\cc inner integral} of \eqref{eqn:fullInt}, $I_{p,q}^{(1)}$. The {\cc outer integral} is computed similarly, by using formulae for the {\c inner integral}, and by suitably replacing the parameters $a,c$ and $q,r$ and variables $x,y$.

First, we examine a border case of $I_{p,q}^{(1)}$ in \eqref{eqn:general_intdx} when $p=-1$. 
\begin{lma}
Let $q \geq 0$. Then
\begin{align}\label{eqn:Rm1xq_F}
I_{-1,q}^{(1)} = R P_{q-1} + C_{q-1} y^q I_{-1,0}^{(1)},
\end{align}
where 
\begin{align}\label{eqn:Rm1AndRm1x}
I_{-1,0}^{(1)} = \frac{1}{\sqrt{a}}\log(2 \sqrt a R + 2ax + by).
\end{align}
The bivariate homogeneous polynomials $P_{q} = P_q(x,y)$ of degree $q$ and coefficients $C_q$ are defined as
\begin{align*}
P_q=
\left\{\begin{array}{ll}
0, & q=-1\\
\frac{1}{a}, & q=0\\
\frac{1}{a(q+1)} \left(x^{q} - \frac{b(2q +1)}{2}y P_{q-1} - c q y^2 P_{q-2}\right),	& q \geq 1
\end{array}\right.,
\qquad
C_q=
\left\{\begin{array}{ll}
1, & q=-1\\
-\frac{b}{2a}, & q=0\\
\frac{1}{a(q+1)} \left(- \frac{b(2q + 1)}{2}C_{q-1} - c q C_{q-2}\right), & q\geq1
\end{array}\right.
.
\end{align*}
\end{lma}

\begin{pf}
We refer to \cite{integralsBook}, Section 2.26 for the special case $q=0$, $I_{-1,0}^{(1)} = \int R^{-1} \dx$. 
%
%
From the same reference we also know that \begin{align*}
I_{-1,1}^{(1)} = \frac{1}{a}R - \frac{b}{2a} y I_{1,0}^{(1)}.
\end{align*}

When $q\geq 2$ the considered integral can be rewritten in the following way
\begin{align}\label{eqn:Rm1xq}
a \int R^{-1} {x^q} \dx & = \int R^{-1} x^{q-2} (a x^2 + (bxy+cy^2) - (bxy + cy^2)) \dx  \nonumber \\ 
& = \int R x^{q-2} \dx - \int R^{-1} x^{q-2} (bxy + cy^2) \dx.
\end{align}
For the first integral in \eqref{eqn:Rm1xq} we apply integration by parts,
\begin{align}\label{eqn:R_pp}
I_{1,q-2}^{(1)} = \int R x^{q-2} \dx = \frac{1}{q-1}R x^{q-1} - \frac{1}{2(q-1)} \int R^{-1}x^{q-1}(2ax+by) \dx.
\end{align}
By combining recursive formulae \eqref{eqn:Rm1xq} and \eqref{eqn:R_pp} we obtain the formula for $q\geq 2$,
\begin{align*}
I_{-1,q}^{(1)} =  \frac{1}{aq} \left( R x^{q-1} - \frac{b(2q - 1)}{2} y I_{-1,q-1}^{(1)} - c(q-1) y^2 I_{-1,q-2}^{(1)} \right).
\end{align*}
By considering two stopping cases of the recursion, $q=0,1$, after some simplification we can express $I_{-1,q}^{(1)} $ with formula \eqref{eqn:Rm1xq_F}.

\qed
\end{pf}

%


Let us analyze another border case of $I_{p,q}^{(1)}$, when  $p \leq -3$ and $q=0,1$.
\begin{lma}
Let $p\leq -3$. Then
\begin{align}\label{eqn:Rmp}
I_{p,0}^{(1)} 
= \sum_{\ell=p+2,p+4,\dots,-1} \bigg(\prod_{\ell_2 = p,p+2,\dots,\ell - 4} D_{\ell_2} \bigg) \bigg(\prod_{\ell_2 = p,p+2,\dots,\ell-2} E_{\ell_2} \bigg) R^{\ell} (2 a x + b y) y^{p-\ell},
\end{align}
where
\begin{align*}
D_\ell  = -a(2\ell + 6), \qquad E_\ell = -\frac 2 {(4 a c - b^2)(\ell+2)},
\end{align*}
and
\begin{align}\label{eqn:Rmpx}
I_{p,1}^{(1)} = \frac{1}{a (p+2)} R^{p+2} - \frac{b}{2a} y I_{p,0}^{(1)}.
\end{align}
\end{lma}

\begin{pf}
For $p \leq -3$ we integrate $I_{p+2,0}^{(1)}$ by parts,
\begin{align*}
2a \int R^{p+2} \dx = R^{p+2} (2ax+by) - \frac {p+2} 2 \int R^{p} (2ax+by)^2 \dx.
\end{align*}
The expression can be rewritten as
\begin{align*}
0 &= R^{p+2} (2ax+by) - \frac {p+2} 2 \int R^{p} (2ax+by)^2 \dx - 2a \int R^{p+2} \dx\\
  &= R^{p+2} (2ax+by)  - \frac {p+2} 2 \int R^{p} \left ( (2ax+by)^2 - 4a  R^{2} \right) \dx - a(2(p+2)+2) \int  R^{p+2} \dx\\
  &= R^{p+2} (2ax+by)  + \frac{ ( 4ac - b^2 )y^2(p+2)}{2} \int R^{p} \dx - a(2p+6) \int  R^{p+2} \dx,
\end{align*}
hence
\begin{align*}
I_{p,0}^{(1)} &= \int R^{p} \dx = \frac{2}{(4ac - b^2) (p+2)} \left ( a(2p+6)  y^{-2} I_{p+2,0}^{(1)} - R^{p+2}(2ax+by) y^{-2} \right ).
\end{align*}
By continuing with the recursion for higher $p$ we obtain after some simplifications the expression \eqref{eqn:Rmp}.


From
\begin{align*}
R^{p+2} = \int \frac{\partial}{\partial x} R^{p+2} \dx = \int \frac{p+2}{2} R^{p}(2 a x + b y) \dx = a(p+2) I_{p,1}^{(1)} + \frac{(p+2)b}{2} y  I_{p,0}^{(1)}
\end{align*}
we get formula \eqref{eqn:Rmpx}.


\qed
\end{pf}

\begin{figure}[ht!]
\centering
\begin{tikzpicture}[->,>=stealth',auto,node distance=2cm,
  main node/.style={circle,draw,font=\sffamily\Large\bfseries}]
\node[inner sep=0pt] (IG) at (0,0)
{\includegraphics[trim = 2.5cm 0.0cm -5cm 0.25cm, clip = true, height=6cm]{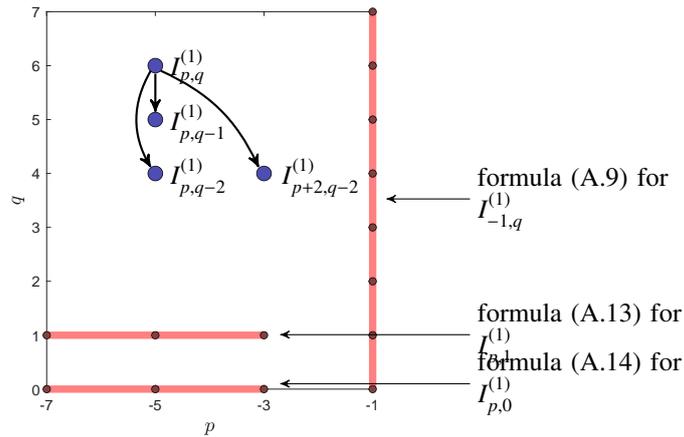}};


\node[text width=1cm] at (-1.95,1.95) {$I_{p,q}^{(1)}$};
\node[text width=1cm] at (-1.95,1.2) {$I_{p,q-1}^{(1)}$};
\node[text width=1cm] at (-1.95,0.48) {$I_{p,q-2}^{(1)}$};
\node[text width=1cm] at (-0.48,0.48) {$I_{p+2,q-2}^{(1)}$};

\draw[->] (1.5,0.2) -- (0.4,0.2);
\node[text width=3cm] at (3.1,0.2) {formula \eqref{eqn:Rm1xq_F} for $I_{-1,q}^{(1)}$};

\draw[->] (1.5,-1.6) -- (-1,-1.6);
\node[text width=3cm] at (3.1,-1.6) {formula \eqref{eqn:Rmp}  for $I_{p,1}^{(1)} $};

\draw[->] (1.5,-2.25) -- (-1,-2.25);
\node[text width=3cm] at (3.1,-2.25) {formula \eqref{eqn:Rmpx}  for $I_{p,0}^{(1)} $};

  \path[every node/.style={font=\sffamily\small}, thick] 
    (-2.64,1.86) edge[] node [left] {} (-2.64,1.35);
    
      \path[every node/.style={font=\sffamily\small}, thick] 
    (-2.57,1.9) edge[bend left=20] node [left] {} (-1.28,0.62);
    
      \path[every node/.style={font=\sffamily\small}, thick] 
    (-2.7,1.9) edge[bend right] node [left] {} (-2.7,0.62);

\end{tikzpicture}

%
%
\caption{Illustration of the recursion for $I_{p,q}^{(1)}$ in \eqref{eqn:recurrence_Ipq}.}
\label{fig:recursionScheme}
\end{figure}

When $p\leq -3$ and $2 \leq q$ the integral $ I_{p,q}^{(1)}$ is computed recursively from integrals of higher and lower powers $p$ and $q$, respectively, by applying formula $a R^{p} x^{q} = R^{p} x^{q-2}(a x^2 + (bxy+cy^2) - (bxy+cy^2))$,
\begin{align} \label{eqn:recurrence_Ipq}
I_{p,q}^{(1)}  = \frac{1}{a}I_{p+2,q-2}^{(1)} - \frac{b}{a} y I_{p,q-1}^{(1)} - \frac{c}{a}  y^2 I_{p,q-2}^{(1)},
\end{align}
see Fig.~\ref{fig:recursionScheme}. The recursion terminates when either the case $p=-1$ is reached and we use the formula \eqref{eqn:Rm1xq_F}, or the case $p \leq -3, q=0,1$ is encountered and we use the expressions \eqref{eqn:Rmp} and \eqref{eqn:Rmpx}, respectively.

We are left to show how to compute the outer integration of $I_{p,q,r}$, using {\cc similar formulae to} the inner integration.
%
%
\begin{prn}\label{prn:I_pqr}
Integral $I_{p,q,r}$ is a linear combination of integrals  $y^{\zeta_0+2}   I_{-1,0}^{(1)}$ and 
$x^{\zeta_0-\hat p -\hat r +1} I_{\hat p,\hat r}^{(2)}$ for $\hat p= p+2,p+4,\dots,-1$, $\hat r=0,1,\dots,\zeta_0-\hat p +1$, where $\zeta_0 :=p +q +r$ .

\end{prn}

\begin{pf}
Using recursion \eqref{eqn:recurrence_Ipq}  on $y^r I_{p,q}^{(1)} $ we obtain
\begin{align} \label{eqn:Ipq_final}
I_{p,q,r} = \int  y^r I_{p,q}^{(1)}  \dy = \sum_{\hat q=0,1,\dots,q} \alpha_{\hat q} \int y^{\zeta_0-\hat q+1} I_{-1,\hat q}^{(1)} \dy + \sum_{\hat p = p,p+2,\dots,-1} \int \beta_{\hat p} y^{\zeta_0-\hat p} I_{\hat p,0}^{(1)} +  \gamma_{\hat p} y^{\zeta_0-\hat p-1} I_{\hat p,1}^{(1)} \dy,
\end{align}
for suitable coefficients $\alpha_{\hat p}, \beta_{\hat p}, \gamma_{\hat q}$. Observe that $\zeta$ parameter is preserved (equal to $\zeta_0$) for all integrands inside double integrals on the right-hand side of \eqref{eqn:Ipq_final}.

By applying formulae \eqref{eqn:Rmp} and \eqref{eqn:Rmpx} on the second term of \eqref{eqn:Ipq_final}, the term writes as
\begin{align}\label{eqn:compact1}
\sum_{\hat p= p+2,p+4,\dots,-1} \int R^{\hat p} P_{\zeta_0-\hat p+1}^{[1]} \dy
\end{align}
for some bivariate homogeneous polynomials $P_{\zeta_0-\hat p+1}^{[1]}$ of degree $\zeta_0-\hat p+1$, hence the sum is a linear combination of $x^{\zeta_0-\hat p -\hat r +1} I_{\hat p,\hat r}^{(2)}$ for $\hat p= p+2,p+4,\dots,-1$ and $\hat r=0,1,\dots,\zeta_0-\hat p +1$. 

On each integral of the first term of \eqref{eqn:Ipq_final} we apply formula \eqref{eqn:Rm1xq_F}. Each non-logarithmic term can be rewritten as
\begin{align*}
\int R P_{\zeta_0}^{[2]} \dy = \int R^{-1} (R^{2}P_{\zeta_0}^{[2]}) \dy =: \int R^{-1} P_{\zeta_0+2}^{[3]} \dy,
\end{align*}
for suitable polynomials $P_{\zeta_0}^{[2]}, P_{\zeta_0+2}^{[3]}$, hence it can be treated in the same manner as the addend in \eqref{eqn:compact1} for $\hat p=-1$. The remaining integrals of \eqref{eqn:Ipq_final} are integrated by parts using
\begin{align*}
\int  y^{\zeta_0+1}  I_{-1,0}^{(1)} \dy
=\frac{y^{\zeta_0+2}}{\zeta_0 + 2}  I_{-1,0}^{(1)}  - \frac{1}{\zeta_0 + 2} \int \left(\frac{1}{\sqrt{a}}-xR^{-1} \right) y^{\zeta_0+1} \dy
=\frac{1}{\zeta_0 + 2}\left( y^{\zeta_0+2} I_{-1,0}^{(1)}  - \frac{y^{\zeta_0+2}}{\sqrt{a} (\zeta_0+2)}   + x I_{-1,\zeta_0+1}^{(2)} \right).
\end{align*}
We can ignore the second addend since it does not depend on variable $x$.

\qed
\end{pf}

The integral $I_{p,q,r}^{[\rm t]}$ from \eqref{eqn:fullInt_t} is computed similarly.
\begin{cor}
Let  $\hat a=a$, $\hat b:=(-2a+b)$, $\hat c:=a+b+c$. Then integral $I_{p,q,r}^{(t)}$  is a linear combination of integrals  $y^{\zeta_0+2}   I_{-1,0}^{(1)}|_{x=1-y}$ and 
$I_{\hat p,\hat r}^{(2)}$, involving integrand function $R(\hat a, \hat b, \hat c, (1,y))$, $\hat p= p+2,p+4,\dots,-1$, $\hat r=0,1,\dots,\zeta_0-\hat p +1$ and $\zeta_0 := p +q +r$.
\end{cor}
\begin{pf}
By applying similar arguments as	in the proof of Proposition~\ref{prn:I_pqr}, 
$I_{p,q,r}^{(t)} = \int  y^r I_{p,q}^{(1)}|_{x=1-y} \dy$ is computed from relation \eqref{eqn:Ipq_final} and by using $x=1-y$. 
Instead of \eqref{eqn:compact1} we get
\begin{align}\label{eqn:nekaj}
\sum_{\hat p= p+2,p+4,\dots,-1} \int (R^{\hat p} P_{\zeta_0-\hat p+1}^{[1]})|_{x=1-y} \dy =
%
\sum_{\hat p= p+2,p+4,\dots,-1} \int {R(\hat a, \hat b, \hat c, (1,y))}^{\hat p} P_{\zeta_0-\hat p+1}^{[1]}|_{x=1-y} \dy,
\end{align}
Thus we can write each term in the sum \eqref{eqn:nekaj} as a linear combination of  $I_{\hat p,\hat r}^{(2)}$, for $\hat p= p+2,p+4,\dots,-1$, $\hat r = 0,1,\dots,\zeta_0-\hat p+1$ and for the modified parameters $\hat a, \hat b, \hat c$.

For the remaining integral it holds
\begin{align*} 
\int y^{\zeta_0+1} I_{-1,0}^{(1)}|_{x=1-y}  \dy = \frac{y^{\zeta_0+2}}{\zeta_0+2} I_{-1,0}^{(1)}|_{x=1-y}   -  \frac{1}{\zeta_0+2} \int \left( \frac 1 {\sqrt a}  - R(\hat a, \hat b, \hat c, (1,y))^{-1} \right) y^{\zeta_0+1} \dy
\end{align*}
and
\begin{align*} 
I_{-1,0}^{(1)}|_{x=1-y} = \frac{1}{\sqrt{a}}\log(2 \sqrt a R(\hat a, \hat b, \hat c, (1,y)) + ( - 2a+b)y +2a).
\end{align*}

\qed
\end{pf}

\end{document}